%% file: incsamp_techreport.tex
\documentclass[10pt]{article}

\RequirePackage{amsmath}
\RequirePackage{amsfonts}
\RequirePackage{amssymb}
\RequirePackage{hyperref}
\RequirePackage{multirow}
\RequirePackage{algorithm}
\RequirePackage{algpseudocode}
\RequirePackage{graphicx}
\RequirePackage{xcolor}
\RequirePackage{amsthm}
\RequirePackage{amsmath}
\RequirePackage{amsfonts}
\RequirePackage{amssymb}
\RequirePackage{hyperref}
\RequirePackage{multirow}
\RequirePackage{algorithm}
\RequirePackage{algpseudocode}
\RequirePackage{graphicx}
\RequirePackage{xcolor}
\allowdisplaybreaks

\newif\iftechreport
\techreporttrue

\usepackage{dpr_fec,ifthen,dsfont}

\newtheorem{definition}{Definition}[section]
\newtheorem{lemma}{Lemma}[section]

\newtheorem{theorem}{Theorem}[section]

\usepackage{isetechreport}
\def\reportyear{25T}
\def\reportno{017}
\def\revisionno{0}
\def\originaldate{\today}
\def\revisiondate{\today}
\coraltrue
\cvcrfalse

\begin{document}

\title{\input{title}}

\author{Frank E.~Curtis\thanks{E-mail: frank.e.curtis@lehigh.edu}}
\author{Lingjun Guo\thanks{E-mail: lig423@lehigh.edu}}
\author{Daniel P.~Robinson\thanks{E-mail: daniel.p.robinson@lehigh.edu}}
\affil{Department of Industrial and Systems Engineering, Lehigh University}
\titlepage

\maketitle

\begin{abstract}
  \input{abstract}
\end{abstract}

\input{body}

\bibliographystyle{plain}
\bibliography{references}

\appendix
\input{appendix}

\end{document}

%% file: title.tex
Progressively Sampled Equality-Constrained Optimization

%% file: abstract.tex
An algorithm is proposed, analyzed, and tested for solving continuous nonlinear-equality-constrained optimization problems where the objective and constraint functions are defined by expectations or averages over large, finite numbers of terms.  The main idea of the algorithm is to solve a sequence of related problems, each involving finite samples of objective- and constraint-function terms, over which the sample sets grow progressively.  Under assumptions about the problem functions and their first- and second-order derivatives that are reasonable in real-world settings of interest, it is shown that---with sufficiently large initial sample sizes---solving a sequence of problems defined through progressive sampling yields a better worst-case sample complexity bound compared to solving a single problem with the full sets of samples. The results of numerical experiments with a set of test problems demonstrate that the proposed approach can be effective in practice.

%% file: body.tex
\section{Introduction}\label{sec.introduction}

We propose, analyze, and test the performance of an algorithm for solving nonlinear-equality-constrained continuous optimization problems.  Our particular setting of interest is when the objective and equality-constraint functions are defined by expectations or averages over large, finite numbers of terms.  A broad setting in which such a problem arises is the context of least-squares regression when one knows, in advance, a target residual error, at least for a subset of the regression terms.  Another class of problems that fit our setting include some arising from physics-informed supervised machine learning where the constraints state that the expected residual of a differential equation should be zero.  

The inspiration for our proposed algorithm is the method for solving unconstrained continuous optimization problems proposed in \cite{mokhtari2019efficient}.  In that work, the aim is to minimize an objective function, which may be nonconvex, when it is defined by an expectation or an average over a large finite number of terms.  The proposed algorithm involves solving a sequence of subproblems, where in each subproblem the objective function is defined by an average over a sample of the original objective terms.  It is shown that, under certain conditions about the original objective function and its sample terms, one can obtain an improved worst-case gradient-sample-complexity bound by solving a sequence of subproblems over which the sample sets grow progressively.  That is, each sample set is a superset of the previous one.  The main idea is to exploit the fact that, under their stated conditions, a minimizer of a sampled objective is close to a minimizer of an objective corresponding to an augmented sample set.

Our aim in this paper is to consider a similar algorithmic approach, except that instead of sampling only objective function terms, we sample terms defining the equality constraints as well.  Our algorithm employs some similar ideas as proposed in \cite{mokhtari2019efficient}, but our setting is quite distinct for multiple reasons.  First, our setting necessitates consideration of \emph{constrained} subproblems and worst-case complexity bounds for finding second-order stationary points of such problems, as opposed to consideration of only the \emph{unconstrained} setting as in \cite{mokhtari2019efficient}.  With respect to our constrained setting, we show that one can rely, e.g., on recent complexity bounds developed in \cite{GoyeEfteBoum2024,DimoONei2026}.  Second, the presence of constraints means that our method needs to rely on approximate stationarity conditions for constrained optimization when determining when one subproblem has been solved to sufficient accuracy and the next in the sequence should be solved.  With respect to this, we rely on properties of least-squares multipliers for ensuring our overall worst-case sample complexity guarantees.  Finally, due to fundamental changes in an optimization problem and its solutions that result when a constraint function is perturbed, our overall worst-case sample complexity bound relies on theory pertaining to \emph{acute perturbations}; see \cite{Stew1973,Stew1977,Stew1979}. For the unconstrained setting in \cite{mokhtari2019efficient}, a fundamental role is played by an assumption that the objective function is \emph{strongly Morse}.  Our analysis also involves a similar assumption (with respect to a Lagrangian function), but this alone is not sufficient.  In the context of equality constraints, the sample sizes need to be sufficiently large in order to ensure that the constraint Jacobians are acute perturbations of each other, as revealed in our analysis.

We also include the results of numerical experiments to demonstrate that computational gains can be achieved through the use of our proposed algorithm.  Our main experiment focuses on a physics-informed supervised learning problem, where we show that progressive sampling can offer reduced computational cost---in terms of accesses to data points---compared to solving a full-sample problem directly.

\subsection{Notation}\label{sec.notation}

We use $\R{}$ to denote the set of real numbers, $\R{}_{\geq r}$ (resp.,~$\R{}_{>r}$) to denote the set of real numbers greater than or equal to (resp.,~greater than) $r \in \R{}$, $\R{n}$ to denote the set of $n$-dimensional real vectors, and $\R{m \times n}$ to denote the set of $m$-by-$n$-dimensional real matrices.  We define $\N{} := \{0,1,2,\dots\}$, and, for any positive integer $N \geq 1$, we use $[N]$ to denote the set $\{1, \dots, N\}$. For any finite set $\Scal$, we use $|\Scal|$ to denote its cardinality. For vectors we define $\|\cdot\| := \|\cdot\|_2$, or else specify the norm explicitly.  We use $\|\cdot\|$ to denote the spectral norm of any matrix.

For any matrix $A \in \R{m \times n}$, we use $\sigma_i(A)$ for each $i \in [\min\{m,n\}]$ to denote its $i$th largest singular value.  Given any such $A$, we use $\Null(A)$ to denote its null space, i.e., $\{d \in \R{n} : Ad = 0\}$.  Assuming that $B \in \R{n \times m}$ has full-column rank, its Moore-Penrose pseudoinverse is $B^\dag := (B^TB)^{-1}B^T$.  Note that in this case of $B$ having full-column rank, the pseudoinverse $B^\dag$ is a left inverse of $B$ in the sense that $B^\dag B = I$.  For any subspace $\Xcal \subseteq \R{n}$ and point $x \in \R{n}$, we denote the Euclidean projection of $x$ onto $\Xcal$ as $\proj_\Xcal(x) := \arg\min_{\xbar \in \Xcal} \|\xbar - x\|$.  Given $B \in \R{n \times m}$ with full column rank, we use $\Rcal(B) := BB^\dag$ and $\Ncal(B^T) = I - \Rcal(B)$ to denote orthogonal projection matrices onto the span of the columns of $B$ and the null space of $B^T$, respectively.  That such $\Rcal(B)$ and $\Ncal(B^T)$ are orthogonal projection matrices follows since each of them are both symmetric and idempotent.

\subsection{Organization}

We state our problem class of interest and our proposed algorithm in \S\ref{sec.algorithm}.  Our analysis of our algorithm is given in \S\ref{sec.analysis}.  The results of our numerical experiments are provided in \S\ref{sec.numerical} and a conclusion is provided in \S\ref{sec.conclusion}.  \iftechreport Finally, a proof that a certain subproblem solver yields the worst-case complexity properties required by our analysis in \S\ref{sec.analysis} is provided in an appendix. \fi

\section{Algorithm}\label{sec.algorithm}

Our proposed algorithm is designed to solve a sample average approximation (SAA) of the continuous nonlinear-equality-constrained problem
\bequation\label{prob.expected}
  \min_{x \in \R{n}}\ \fbar(x)\ \st\ \cbar(x) = 0,
\end{equation}
where the objective and constraint functions, i.e., $\fbar : \R{n} \to \R{}$ and $\cbar : \R{n} \to \R{m}$ (with $m \leq n$), respectively, are twice-continuously differentiable and defined by expectations.  Formally, with respect to random variables $\omega_f$ and $\omega_c$ defined respectively by probability spaces $(\Omega_f,\Fcal_f,\P_f)$ and $(\Omega_c,\Fcal_c,\P_c)$, expectations $\E_f$ and $\E_c$ defined respectively by $\P_f$ and $\P_c$, and $\Fbar : \R{n} \times \Omega_f \to \R{}$ and $\Cbar : \R{n} \times \Omega_c \to \R{m}$, the functions $\fbar$ and $\cbar$ are defined by $\fbar(x) = \E_f[\Fbar(x,\omega_f)]$ and $\cbar(x) = \E_c[\Cbar(x,\omega_c)]$ for all $x \in \R{n}$.

The SAA of problem~\eqref{prob.expected} that our algorithm is designed to solve is defined with respect to a (large) sample of $N_f \in \N{}$ realizations of $\omega_f$ and a (large) sample of $N_c \in \N{}$ realizations of $\omega_c$, say, $\{\omega_{f,i}\}_{i\in[N_f]}$ and $\{\omega_{c,i}\}_{i\in[N_c]}$.  Defining the SAA functions $f : \R{n} \to \R{}$ and $c : \R{n} \to \R{m}$ for all $x \in \R{n}$ by
\begin{align*}
  f(x) = \tfrac{1}{N_f} \sum_{i \in [N_f]} f_i(x),\ \ &\text{where}\ \ f_i(x) \equiv \Fbar(x,\omega_{f,i})\ \ \text{for all}\ \ i \in [N_f] \\ \text{and}\ \ 
  c(x) = \tfrac{1}{N_c} \sum_{i \in [N_c]} c_i(x),\ \ &\text{where}\ \ c_i(x) \equiv \Cbar(x,\omega_{c,i})\ \ \text{for all}\ \ i \in [N_c],
\end{align*}
the problem that our algorithm is designed to solve is that given by
\bequation\label{prob.opt.N}
  \min_{x \in \R{n}}\ f(x)\ \st\ c(x) = 0.
\eequation
Under reasonable assumptions about the expectation functions $\fbar$ and $\cbar$, and an assumption that both $N_f$ and $N_c$ are sufficiently large, differences between values of $\fbar$ and $\cbar$ in \eqref{prob.expected} with $f$ and $c$ in \eqref{prob.opt.N}, as well as values of their first- and second-order derivatives at any $x \in \R{n}$, can be bounded with high probability \cite{BoucLugoMass2013,MeiBaiMont2018,Geer2000}. These bounds can, in turn, be used to relate approximate stationary points of~\eqref{prob.opt.N} with those of~\eqref{prob.expected}, again with high probability.  Thus, for our purposes, we focus on our proposed algorithm and our analysis of it for solving problem~\eqref{prob.opt.N}.

The main idea of our proposed algorithm for solving~\eqref{prob.opt.N} is to generate a sequence of iterates, each of which is a stationary point (at least approximately) with respect to a sampled problem involving only nonempty subsets $\Scal_f \subseteq [N_f]$ and $\Scal_c \subseteq [N_c]$ of objective and constraint function terms, respectively.  For any such $\Scal := (\Scal_f,\Scal_c)$, we denote the approximate objective as $f_\Scal : \R{n} \to \R{}$ and the approximate constraint function as $c_\Scal : \R{n} \to \R{m}$, and state the approximation of problem~\eqref{prob.opt.N} as
\bequation\label{prob.opt.S}
  \baligned
    \min_{x \in \R{n}}\ f_\Scal(x)\ &\st\ c_\Scal(x) = 0, \\ \text{where}\ \ f_\Scal(x) = \tfrac{1}{|\Scal_f|} \sum_{i\in{\Scal_f}} f_i(x) \ \ &\text{and}\ \ c_\Scal(x) = \tfrac{1}{|\Scal_c|} \sum_{i\in{\Scal_c}} c_i(x).
  \ealigned
\end{equation}
We use the subscript $\Scal$ for both $f$ and $c$, even though they only rely on $\Scal_f$ and~$\Scal_c$, respectively, to simplify our notation throughout the remainder of the paper. Observe that with this notation the choice $\Scal = ([N_f],[N_c])$ gives the objective and constraint functions in~\eqref{prob.opt.N}.  The primary benefit of considering \eqref{prob.opt.S} for the pair of samples $\Scal = (\Scal_f,\Scal_c)$, rather than \eqref{prob.opt.N} directly, is that any evaluation of the objective function, its gradient, or its Hessian requires a sum of only $|\Scal_f| \leq N_f$ terms, as opposed to $N_f$ terms, and any evaluation of the constraint function, its Jacobian, or any constraint Hessian function requires a sum of only $|\Scal_c| \leq N_c$ terms, as opposed to $N_c$ terms.  Also, under reasonable assumptions about the objective and constraint functions, we show in this paper that, by starting with an approximate stationary point for~\eqref{prob.opt.S} and aiming to solve a subsequent instance of~\eqref{prob.opt.S} with respect to supersets $\overline\Scal_f \supseteq \Scal_f$ and $\overline\Scal_c \supseteq \Scal_c$, our proposed algorithm can obtain an approximate stationary point for the subsequent instance with lower sample complexity than if the problem with the larger sample set were solved directly from an arbitrary starting point.  Overall, we show that for sufficiently large sample sets relative to~$(N_f,N_c)$, a sufficiently approximate stationary point of~\eqref{prob.opt.N} can be obtained more efficiently through progressive sampling than by solving the full-sample problem directly.

Now let us introduce stationarity conditions for~\eqref{prob.opt.S} for any pair of nonempty sample sets $\Scal = (\Scal_f,\Scal_c) \subseteq [N_f] \times [N_c]$, which in particular also represent stationarity conditions for~\eqref{prob.opt.N} when one considers the case that $\Scal_f = [N_f]$ and $\Scal_c = [N_c]$.  Let us assume that $f_{\Scal}$ and $c_{\Scal}$ are twice-continuously differentiable for any pair of nonempty sets $\Scal \equiv (\Scal_f,\Scal_c) \subseteq [N_f] \times [N_c]$.  Let the Lagrangian of problem~\eqref{prob.opt.S} be denoted by $L_\Scal : \R{n} \times \R{m} \to \R{}$, defined for all $(x,y) \in \R{n} \times \R{m}$ by
\bequationNN
  L_\Scal(x,y) = f_\Scal(x) + c_\Scal(x)^Ty = \tfrac{1}{|\Scal_f|} \sum_{i \in \Scal_f} f_i(x) + \tfrac{1}{|\Scal_c|} \sum_{i \in \Scal_c} c_i(x)^Ty,
\eequationNN
where $y \in \R{m}$ is a vector of Lagrange multiplier estimates (also known as dual variables).  When $\Scal = (\Scal_f,\Scal_c) = ([N_f],[N_c])$, we denote the Lagrangian simply as $L$.  Second-order necessary conditions for optimality for \eqref{prob.opt.S} can then be stated as
\bequation\label{eq.soc.S.1}
  \nabla L_\Scal(x,y) \equiv \bbmatrix \nabla_x L_\Scal(x,y) \\ \nabla_y L_\Scal(x,y) \ebmatrix \equiv \bbmatrix \nabla f_\Scal(x) + \nabla c_\Scal(x)y \\ c_\Scal(x) \ebmatrix = 0
\eequation
and, with $[c_\Scal]_j$ denoting the $j$th component of the constraint function $c_\Scal$,
\begin{align}
  d^T \nabla_{xx}^2 L_\Scal(x,y)d \equiv d^T \( \nabla^2 f_\Scal(x) + \sum_{j \in [m]} \nabla^2 [c_\Scal]_j(x) y_j \) d &\geq 0 \label{eq.soc.S.2} \\
  \text{for all}\ \ d &\in \Null(\nabla c_{\Scal}(x)^T). \nonumber
\end{align}
We refer to any $(x,y) \in \R{n} \times \R{m}$ satisfying \eqref{eq.soc.S.1} as a first-order stationary point with respect to~\eqref{prob.opt.S}, and we refer to any such point satisfying both \eqref{eq.soc.S.1} and \eqref{eq.soc.S.2} as a second-order stationary point.  Also, consistent with the literature on worst-case complexity for constrained optimization (e.g., see \cite{DimoONei2026}), we say for $(\epsilon,\zeta) \in \R{}_{>0} \times \R{}_{>0}$ that $(x,y) \in \R{n} \times \R{m}$ is $(\epsilon,\zeta)$-stationary with respect to~\eqref{prob.opt.S} if and only if
\bsubequations\label{eq.soc.S.approx}
  \begin{align}
    \| \nabla L_\Scal(x,y) \| &\leq \epsilon \label{eq.soc.S.approx.1} \\ \text{and}\ \ 
    d^T \nabla_{xx}^2 L_\Scal(x,y)d &\geq -\zeta\|d\|_2^2\ \ \text{for all}\ \ d \in \Null(\nabla c_{\Scal}(x)^T). \label{eq.soc.S.approx.2}
  \end{align}
\esubequations

Generally speaking, an algorithm for solving \eqref{prob.opt.S} can be a \emph{primal} method that only generates a sequence of primal iterates $\{x_k\}$, or it can be a \emph{primal-dual} method that generates a sequence of primal and dual iterate pairs $\{(x_k,y_k)\}$.  For an application of our proposed algorithm, either type of method can be employed, but for certain results in our analysis, we use properties of \emph{least-square multipliers} corresponding to a given primal point $x \in \R{n}$.  Assuming that the Jacobian of~$c_\Scal$ at~$x$ (i.e., $\nabla c_\Scal(x)^T$) has full-row rank, the least-squares multipliers with respect to a primal point $x$ are given by $y_\Scal(x) \in \R{m}$ that minimizes $\|\nabla_x L_\Scal(x,\cdot)\|^2$, which is denoted by
\begin{equation}\label{eq.lsm}
  y_\Scal(x) = - (\nabla c_\Scal(x)^T \nabla c_\Scal(x))^{-1} \nabla c_\Scal(x)^T \nabla f_\Scal(x) = - \nabla c_\Scal(x)^\dag \nabla f_\Scal(x).
\end{equation}
The computation of the least-squares multipliers is negligible relative to the overall cost of our method.  Moreover, in practice, they need not be computed by the subproblem solver during every iteration.  A reasonable approach would be to employ a subproblem solver until the prescribed accuracy has been obtained with some multipliers, then compute least-squares multipliers.  The first-order conditions will remain satisfied to the prescribed accuracy, so one only needs to confirm whether the second-order conditions remain satisfied as well, or whether higher accuracy is needed.

Our proposed method is stated as Algorithm~\ref{alg.psm} below.  Our analysis in the next section formalizes assumptions under which Algorithm~\ref{alg.psm} is well defined and yields our claimed worst-case sample complexity guarantees.

\balgorithm
  \caption{Progressive Sampling Algorithm for \eqref{prob.opt.N}}
  \label{alg.psm}
  \balgorithmic[1]
    \Require Initial sample sizes $(p_{f,1},p_{c,1}) \in [N_f] \times [N_c]$, sample increase factors $(\theta_f,\theta_c) \in (1,\infty) \times (1,\infty)$, initial point $x_0 \in \R{n}$, iteration limit $K = \max \{\lceil \log_{\theta_f} \tfrac{N_f}{p_{f,1}} \rceil, \lceil \log_{\theta_c} \tfrac{N_c}{p_{c,1}} \rceil\}$, and tolerances $\{(\epsilon_k,\zeta_k)\}_{k=1}^K \subset \R{}_{>0} \times \R{}_{>0}$
    \State set $\Scal_{f,0} \gets \emptyset$ and $\Scal_{c,0} \gets \emptyset$
    \For {$k \in [K]$}
      \State choose $\Scal_{f,k} \supseteq \Scal_{f,k-1}$ with $|\Scal_{f,k}| = p_{f,k}$ and $\Scal_{c,k} \supseteq \Scal_{c,k-1}$ with $|\Scal_{c,k}| = p_{c,k}$
      \State using $x_{k-1}$ as a starting point, employ an algorithm to solve \eqref{prob.opt.S}, terminating once a primal iterate $x_k$ has been obtained such that $(x_k,y_{\Scal_k}(x_k))$ (see \eqref{eq.lsm}) is $(\epsilon_k,\zeta_k)$-stationary with respect to subproblem~\eqref{prob.opt.S} for $\Scal = \Scal_k := (\Scal_{f,k},\Scal_{c,k})$
      \State set $p_{f,k+1} \gets \min\{\theta_f p_{f,k}, N_f\}$ and $p_{c,k+1} \gets \min\{\theta_c p_{c,k}, N_c\}$
    \EndFor
    \State \Return $(x_K,y(x_K))$, which is $(\epsilon_K,\zeta_K)$-stationary with respect to \eqref{prob.opt.N}
  \ealgorithmic
\ealgorithm

\section{Analysis}\label{sec.analysis}

Our analysis of Algorithm~\ref{alg.psm} focuses on generic worst-case iteration and sample complexity bounds under broad assumptions about the objective and equality-constraint functions, as well as on the worst-case complexity properties of the algorithm that is employed for solving the arising sampled subproblems.  We state these bounds in generic terms since they hold for any subproblem solver that possesses the stated worst-case complexity properties. We have also proved that a specific subproblem solver based on minimizing Fletcher's augmented Lagrangian function possesses the worst-case complexity properties that lead to desirable worst-case sample complexity bounds. \iftechreport See Appendix~\ref{sec.specific}. \fi

Algorithm~\ref{alg.psm} is written in full generality wherein the sequences of sample sizes $\{|\Scal_{f,k}|\}$ and $\{|\Scal_{c,k}|\}$ may be initialized differently and increase with different factors, meaning that the number of outer iterations (over $k$) required until $|\Scal_{f,k}| = N_f$ or $|\Scal_{c,k}| = N_c$ may be different.  This is reflected in the $\max$ in the definition of $K$.  To simplify our notation throughout our analysis, let us make the mild simplification that $N_f = N_c$, $p_{f,1} = p_{c,1}$, and $\theta_f = \theta_c$.  In this manner, we can refer simply to
\bequationNN
  N := N_f = N_c,\ \ p_k := p_{f,k} = p_{c,k}\ \ \text{for all}\ \ k \in [K],\ \ \text{and}\ \ \theta := \theta_f = \theta_c.
\eequationNN
Let us also abuse notation and refer to the sample sets and their cardinalities as
\bequationNN
  \Scal_k := \Scal_{f,k} = \Scal_{c,k}\ \ \text{and}\ \ |\Scal_k| := |\Scal_{f,k}| = |\Scal_{c,k}|\ \ \text{for all}\ \ k \in [K].
\eequationNN
All of the results that are proved in this section extend readily to the more general setting as stated in Algorithm~\ref{alg.psm}, only with more cumbersome index notation.

Let us begin by stating the basic assumptions under which we prove our theoretical guarantees in this subsection.  Assumption~\ref{ass.subproblems} below ensures that any minimizer of each subproblem is a second-order stationary point and that one can expect an algorithm that is employed to solve each subproblem to find a sufficiently approximate second-order stationary point.  It would be possible to prove reasonable convergence guarantees for Algorithm~\ref{alg.psm} under looser assumptions.  For example, if an algorithm employed to solve \eqref{prob.opt.S} for some nonempty sample set~$\Scal$ were to encounter an (approximate) infeasible stationary point, then it would be reasonable to terminate the subproblem solver and either terminate Algorithm~\ref{alg.psm} in its entirety or move on to solve the next subproblem (with a larger sample set).  However, since consideration of such scenarios would distract from the essential properties of our algorithm, we make Assumption~\ref{ass.subproblems}.  Our remarks after the assumption justify it further.

\bassumption\label{ass.subproblems}
  The functions $f : \R{n} \to \R{}$ and $c_i : \R{n} \to \R{m}$ for each $i \in [N]$ are twice-continuously differentiable.  In addition, the following hold for~\eqref{prob.opt.N}, instances of~\eqref{prob.opt.S}, and the algorithm employed to solve them.
  \benumerate
    \item[(a)] There exists $\sigma_{\min} \in \R{}_{>0}$ such that, for all $x \in \R{n}$ and $\Scal \subseteq [N]$ with $|\Scal| \geq p_1$, the constraint Jacobian has $\sigma_m(\nabla c_\Scal(x)^T) \geq \sigma_{\min}$.
    \item[(b)] For all $\Scal \subseteq [N]$ with $|\Scal| \geq p_1$ and any initial point, the algorithm employed to solve~\eqref{prob.opt.S} is guaranteed to generate a sequence of iterates for which a limit point is a second-order stationary point satisfying \eqref{eq.soc.S.1} and \eqref{eq.soc.S.2}.
    \item[(c)] For some $(\alpha,\beta) \in \R{}_{>0} \times \R{}_{>0}$, problem~\eqref{prob.opt.N} is $(\alpha,\beta)$-strongly Morse in the sense that, for any $x \in \R{n}$, if $(x,y(x))$ satisfies $\| \nabla L(x,y(x)) \| \leq \alpha$, then
    \bequationNN
      |d^T \nabla_{xx}^2 L(x,y(x)) d| \geq \beta \|d\|_2^2\ \ \text{for all}\ \ d \in \Null(\nabla c(x)^T).
    \eequationNN
  \eenumerate
\eassumption

Part~(a) of Assumption~\ref{ass.subproblems} guarantees that the algorithm employed to solve~\eqref{prob.opt.S} will not, e.g., generate an infeasible stationary point. Since the algorithm involves sampling of the objective and constraint functions, one could loosen this assumption to state that the bound on the singular values holds only with high probability, in which case our subsequent complexity bound would hold with high probability. However, since the core of our subsequent analysis remains the same in either case, we make part~(a) of the assumption for ease of our subsequent discussions.  Part~(b) of Assumption~\ref{ass.subproblems} implicitly requires one of various types of conditions in the literature on equality-constrained optimization that guarantee convergence of an algorithm to a second-order stationary point.  (This assumption will be strengthened in our final theorem of this section.)  Part~(c) of Assumption~\ref{ass.subproblems} is essential for proving our desired worst-case sample complexity properties of Algorithm~\ref{alg.psm}.  Related to~(b) and~(c), the following comments are important, so we emphasize them as formal remarks.

\bremark
  If one were to change the requirement for the subproblem solver and only require that it produces an approximate first-order stationary point, rather than an approximate second-order stationary point, then under Assumption~\ref{ass.subproblems} (with ``second-order'' replaced by ``first-order'' in part (b) and part (c) removed) the algorithm would be well defined and, by its construction, would guarantee convergence to an approximate first-order stationary point of problem~\eqref{prob.opt.N}.  Therefore, in practice, Algorithm~\ref{alg.psm} might be run with only approximate first-order stationarity requirements.  However, such a set-up would not allow us to prove our specific desired strong worst-case sample complexity guarantee, which focuses on attainment of approximate second-order stationarity.  Therefore, for the purposes of our analysis, we state Algorithm~\ref{alg.psm} and Assumption~\ref{ass.subproblems} as they are given.
\eremark

\bremark
  Part~(c) of Assumption~\ref{ass.subproblems} is motivated by the concept of a Morse function~\cite{milnor1963morse}, i.e., one that is twice-continuously differentiable and has a nonsingular Hessian at all first-order stationary points. The concept of a Morse function was extended to that of a Morse program for constrained optimization in~\cite{fujiwara1982morse}, wherein relationships between second-order sufficient conditions and stationary points of Morse programs are provided. Problem \eqref{prob.opt.N} is a Morse program if and only if
  \begin{align*}
    \nabla L(x,y(x)) = 0\ \ \implies\ \ Z(x)^T\nabla_{xx}^2 L(x,y(x)) Z(x)\ \ \text{is nonsingular},
  \end{align*}
  where $Z(x)$ is any basis for $\Null(\nabla c(x)^T)$. From the fact that stationary points of Morse programs are isolated~\cite{guillemin1975topology}, one can conclude that for any Morse program and any bounded region there exists a pair $(\alpha,\beta) \in \R{}_{>0} \times \R{}_{>0}$ such that the Morse program is $(\alpha,\beta)$-strongly Morse in the region. Hence, part~(c) of Assumption~\ref{ass.subproblems} is relatively mild in the context of constrained problems with isolated local minima.
\eremark

Our next assumption articulates bounds on derivatives of the objective and constraint functions corresponding to the full-sample problem~\eqref{prob.opt.N}.  An assumption of this type is typical in the literature on continuous constrained optimization.

\bassumption\label{ass.boundness}
  There exists $(\kappa_{\nabla f}, \kappa_{\nabla c}, \kappa_{\nabla^2 f}, \kappa_{\nabla^2 c}) \in \R{}_{>0} \times \R{}_{>0} \times \R{}_{>0} \times \R{}_{>0}$ such that, for all $x \in \R{n}$ and $j \in [m]$, one has $\| \nabla f(x) \| \leq \kappa_{\nabla f}$, $\| \nabla c(x) \| \leq \kappa_{\nabla c}$, $\| \nabla^2 f(x) \| \leq \kappa_{\nabla^2 f}$, and $\| \nabla^2 [c]_j (x) \| \leq \kappa_{\nabla^2 c}$, where $[c]_j$ is the $j$th component of $c$.
\eassumption

Our next assumption introduces constants that bound discrepancies between single-sample functions and the full-sample functions in \eqref{prob.opt.N}, as well as between derivatives of these functions.  The assumption is quite loose in general.  Indeed, the right-hand sides of each inequality are either (large) positive constants or a product of a (large) positive constant with a norm of a constraint derivative that, under Assumption~\ref{ass.subproblems}, is bounded away from zero.  We note that Assumptions~\ref{ass.boundness} and \ref{ass.bounded.distribute} complement each other in the sense that, together, they allow us to prove similar bounds for general nonempty sample sets $\Scal \subseteq [N]$; see upcoming Lemma~\ref{lemma.sample.average.result}.  The unique form of the bound for the constraint derivatives is needed for our subsequent results about acute perturbations; see upcoming Lemma~\ref{lemma.acute.perturb}.

\bassumption\label{ass.bounded.distribute}
  There exists $(\gamma_f, \gamma_c, \gamma_{\nabla f}, \gamma_{\nabla c}, \gamma_{\nabla^2 f}, \gamma_{\nabla^2 c}) \in \R{}_{>0} \times \R{}_{>0} \times \R{}_{>0} \times \R{}_{>0} \times \R{}_{>0} \times \R{}_{>0}$ such that, for all $(x,j) \in \R{n} \times [m]$, one has
  \begin{align*}
    \tfrac{1}{N} \sum_{i \in [N]} |f_i(x) - f(x) |^2 &\leq \gamma_f, \\
    \tfrac{1}{N} \sum_{i \in [N]} \| c_i(x) - c(x) \|^2 &\leq \gamma_c, \\
    \tfrac{1}{N} \sum_{i \in [N]} |\nabla f_i(x) - \nabla f(x) |^2 &\leq \gamma_{\nabla f}, \\
    \tfrac{1}{N} \sum_{i \in [N]} \| \nabla c_i(x)-\nabla c(x) \|^2 &\leq \gamma_{\nabla c} \|\nabla c (x) \|^2, \\
    \tfrac{1}{N} \sum_{i \in [N]} \| \nabla^2 f_i(x) - \nabla^2 f(x) \|^2 &\leq \gamma_{\nabla^2 f} \\
    \text{and}\ \ 
    \tfrac{1}{N} \sum_{i \in [N]} \| \nabla^2 [c_i]_j(x) - \nabla^2 [c]_j(x) \|^2 &\leq \gamma_{\nabla^2 c}.
  \end{align*}
\eassumption

Like for Assumption~\ref{ass.subproblems}(a), one could loosen the bounds in Assumption~\ref{ass.bounded.distribute} to say that they only hold with high probability. However, again, since the core of our analysis remains the same in either case, we employ Assumption~\ref{ass.bounded.distribute} as stated.

Our first lemma employs bounds introduced in Assumptions~\ref{ass.subproblems} and \ref{ass.boundness} to offer bounds on the pseudoinverse of the constraint Jacobian, the least-squares multipliers, and the Hessian of the Lagrangian at all points in the domain of problem~\eqref{prob.opt.N}.

\blemma\label{lemma.bound.c.y.L}
  For all $x \in \R{n}$, one has $($recall \eqref{eq.lsm}$)$ that
  \bsubequations
    \begin{align}
      \|\nabla c(x)^\dag\| &\leq \sigma_{\min}^{-1}, \label{lemma.bound.nabla_c_dag} \\
      \|y(x)\| &\leq \kappa_{\nabla f} \sigma_{\min}^{-1}, \label{lemma.bound.y_N} \\ \text{and}\ \ 
      \|\nabla^2_{xx} L(x,y(x))\| &\leq \kappa_{\nabla^2 f} + \sqrt{m} \kappa_{\nabla f} \kappa_{\nabla^2 c} \sigma_{\min}^{-1}. \label{lemma.bound.L_Hess}
    \end{align}
  \esubequations
\elemma

\bproof
  First, consider \eqref{lemma.bound.nabla_c_dag}.  As is well known, the spectral norm of the pseudoinverse of a matrix is equal to the reciprocal of its smallest singular value; see, e.g., \cite[Chapter 21]{gallier2019linear}.  Therefore, since under Assumption~\ref{ass.subproblems} the smallest singular value of $\nabla c(x)^T$ is bounded below by $\sigma_{\min}$ for all $x \in \R{n}$, the bound in \eqref{lemma.bound.nabla_c_dag} follows.  Second, consider \eqref{lemma.bound.y_N}.  By \eqref{eq.lsm}, submultiplicity of the matrix 2-norm, Assumption~\ref{ass.boundness}, and~\eqref{lemma.bound.nabla_c_dag} the desired conclusion follows since, for all $x \in \R{n}$, one finds
  \bequationNN
    \|y(x)\| = \|\nabla c(x)^\dag\nabla f(x)\| \leq \|\nabla c(x)^\dag\| \|\nabla f(x)\| \leq \kappa_{\nabla f} \sigma_{\min}^{-1}.
  \eequationNN
  Finally, consider \eqref{lemma.bound.L_Hess}.  By the triangle inequality and absolute homogeneity of matrix norms, Assumption~\ref{ass.boundness}, the fact that for any vector $y \in \R{m}$ one has $\|y\|_1 \leq \sqrt{m} \|y\|$, and \eqref{lemma.bound.y_N}, one has for all $x \in \R{n}$ that
  \begin{align*}
    \|\nabla^2_{xx}L(x,y(x))\|
      &= \bigg\|\nabla^2 f(x) + \sum_{j\in[m]} \nabla^2 [c]_j(x) [y(x)]_j \bigg\| \\
      &\leq \|\nabla^2 f(x)\| + \sum_{j\in[m]} \|\nabla^2 [c]_j(x) \| \cdot |[y(x)]_j| \\
      &\leq \kappa_{\nabla^2 f} + \kappa_{\nabla^2 c}  \|y(x)\|_1 \\
      &\leq \kappa_{\nabla^2 f} + \sqrt{m} \kappa_{\nabla^2 c} \|y(x)\| \leq \kappa_{\nabla^2 f} + \sqrt{m}\kappa_{\nabla f}\kappa_{\nabla^2 c} \sigma_{\min}^{-1},
  \end{align*}
  which is the desired conclusion.
\eproof

Our second lemma leverages Assumption~\ref{ass.bounded.distribute} in order to provide bounds that are similar to those in the assumption, except that they are with respect to any nonempty sample set $\Scal \subseteq [N]$.  The resulting bounds depend on the sample size $|\Scal|$.  For convenience here and throughout the rest of the paper, let us define
\bequation\label{eq.key_value}
  \xi_\Scal := \sqrt{\tfrac{N(N-|\Scal|)}{|\Scal|^2}} \in [0,\sqrt{N(N-1)}] \ \ \text{for all nonempty $\Scal\subseteq [N]$.}
\eequation
One could prove a stronger result than upcoming Lemma~\ref{lemma.sample.average.result} if one were to make stronger assumptions about the distribution of the random variables $\omega_f$ and~$\omega_c$. We address this possibility in further detail in upcoming Remark~\ref{rem.loose}.
  
\blemma\label{lemma.sample.average.result}
  For all $x \in \R{n}$, nonempty $\Scal \subseteq [N]$, and $j \in [m]$, one has
  \bsubequations
  \begin{align}
    | f_\Scal(x) - f(x) |^2 &\leq \xi_\Scal^2 \gamma_f, \label{lemma.saa.f} \\
    \| c_\Scal(x) - c(x) \|^2 &\leq \xi_\Scal^2 \gamma_c, \label{lemma.saa.c} \\
    \| \nabla f_\Scal(x) - \nabla f(x) \|^2 &\leq \xi_\Scal^2 \gamma_{\nabla f}, \label{lemma.saa.nabla_f} \\
    \| \nabla c_\Scal(x) - \nabla c(x) \|^2 &\leq \xi_\Scal^2 \gamma_{\nabla c} \|\nabla c(x)\|^2, \label{lemma.saa.nabla_c} \\
    \| \nabla^2 f_\Scal(x) - \nabla^2 f(x) \|^2 &\leq \xi_\Scal^2 \gamma_{\nabla^2 f}, \label{lemma.saa.Hess_f} \\ \text{and}\ \ 
    \| \nabla^2 [c_\Scal]_j(x) - \nabla^2 [c]_j(x) \|^2 &\leq \xi_\Scal^2 \gamma_{\nabla^2 c}.\label{lemma.saa.Hess_c}
  \end{align}
  \esubequations
\elemma
\bproof
  Each of the desired bounds can be proved in a similar manner.  We prove the fourth bound, namely,  \eqref{lemma.saa.nabla_c}, with respect to the constraint derivatives. The other bounds follow in a similar manner. First, observe that 
  \begin{align*}
   \nabla c_{\Scal}(x) = \tfrac{1}{|\Scal|} \sum_{i\in\Scal} \nabla c_i(x) &= \tfrac{1}{|\Scal|} \sum_{i\in[N]} \nabla c_i(x) -\tfrac{1}{|\Scal|} \sum_{i\in[N] \setminus \Scal} \nabla c_i(x) \\
   &= \tfrac{N}{|\Scal|} \nabla c(x) - \tfrac{1}{|\Scal|} \sum_{i\in[N] \setminus \Scal} \nabla c_i(x).
  \end{align*}
  Second, observe that for any vector in $v \in \R{N-|\Scal|}$ one has from the Cauchy-Schwarz inequality that $(\ones^Tv)^2 \leq \|\ones\|^2 \|v\|_2^2 = (N - |\Scal|) \|v\|_2^2$.  Consequently, with the triangle inequality and Assumption \ref{ass.bounded.distribute}, one finds that
  \begin{align*}
    \| \nabla c_\Scal(x) - \nabla c(x) \|^2 
    &= \bigg\| \tfrac{N}{|\Scal|} \nabla c(x) - \nabla c(x) - \tfrac{1}{|\Scal|} \sum_{i\in[N] \setminus \Scal} \nabla c_i(x) \bigg\|^2 \\
    &= \tfrac{1}{|\Scal|^2} \bigg\| \sum_{i\in[N] \setminus \Scal} (\nabla c(x) - \nabla c_i(x)) \bigg\|^2 \\
    &\leq \tfrac{1}{|\Scal|^2} \(\sum_{i\in[N] \setminus \Scal} \|\nabla c(x) - \nabla c_i(x) \|\)^2 \\
    &\leq \(\tfrac{N-|\Scal|}{|\Scal|^2}\)\sum_{i\in[N]\setminus\Scal}\left\|\nabla c (x)-\nabla c_i(x) \right\|^2\\
    &\le \(\tfrac{N-|\Scal|}{|\Scal|^2}\)\sum_{i\in[N]}\left\|\nabla c (x)-\nabla c_i(x) \right\|^2 \leq \(\tfrac{N-|\Scal|}{|\Scal|^2}\)N\gamma_{\nabla c}\|\nabla c (x)\|^2,
  \end{align*}
  which gives the desired conclusion.
\eproof

Our third lemma refers to the concept of an \emph{acute perturbation} between two real-valued rectangular matrices.  We provide the following definition, then state our lemma, which shows that if $\Scal$ has sufficiently large cardinality, then $\nabla c_\Scal(x)$ and $\nabla c(x)$ are acute perturbations of each other for all $x \in \R{n}$.

\begin{definition}\label{def.acute.perturb}
  Two matrices $A \in \R{n\times m}$ and $B \in \R{n\times m}$ are acute perturbations of each other if and only if $\rank(AA^\dag BA^\dag A) = \rank(A) = \rank(B)$.  In particular, if $A$ and~$B$ are full-column-rank matrices, then $A^\dag A = I$ and $\rank(A) = \rank(B) = m$, so~$A$ and $B$ are acute perturbations of each other if and only if $\rank(AA^\dag B) = m$.
\end{definition}

\begin{lemma}\label{lemma.acute.perturb}
  If nonempty $\Scal\subseteq[N]$ satisfies $|\Scal| \geq p_1$ and
  \bequation\label{ineq.acute.perturb}
    |\Scal| > \tfrac{2N}{1+\sqrt{1+\tfrac{4\sigma_{\min}^2}{\gamma_{\nabla c}\kappa_{\nabla c}^2}}} \in (0,N],
  \eequation
  then $\nabla c(x)$ and $\nabla c_{\Scal}(x)$ are acute perturbations of each other for all $x \in \R{n}$.
\end{lemma}
\bproof 
  By Assumption~\ref{ass.subproblems}, for all $x \in \R{n}$, the matrices $\nabla c_\Scal(x) \in \R{n \times m}$ and $\nabla c (x) \in \R{n \times m}$ have full-column rank.  All that remains is to show that the rank condition in Definition~\ref{def.acute.perturb} holds for all $x \in \R{n}$.  Observe that, for any $x \in \R{n}$,
  \begin{align}
    \nabla c(x) \nabla c (x)^\dag \nabla c_\Scal(x)
    &= \nabla c(x) \nabla c(x)^\dag (\nabla c(x) + \nabla c_\Scal(x) - \nabla c(x)) \nonumber \\
    &= \nabla c(x) (I + \nabla c(x)^\dag (\nabla c_\Scal(x) - \nabla c(x))). \label{eq.nabla.c.acute.formula}
  \end{align}
  Next observe that, for all $x \in \R{n}$, one has from submultiplicity of the matrix 2-norm, Lemma~\ref{lemma.bound.c.y.L}, Lemma~\ref{lemma.sample.average.result}, and Assumption~\ref{ass.boundness} that
  \begin{align}
    \| \nabla c(x)^\dag (\nabla c_\Scal(x) - \nabla c (x) )\|
    &\leq \|\nabla c(x)^\dag\| \|\nabla c_\Scal(x) - \nabla c(x) \| \nonumber \\
    &\leq \tfrac{\xi_\Scal \sqrt{\gamma_{\nabla c}}}{\sigma_{\min}} \|\nabla c(x)\| \leq \tfrac{\xi_\Scal \sqrt{\gamma_{\nabla c}} \kappa_{\nabla c}}{\sigma_{\min}}. \label{ineq.i}
  \end{align}
  Let us now show that \eqref{ineq.acute.perturb} implies that the right-hand side of \eqref{ineq.i} is strictly less than one.  Define $t := \tfrac{N}{|\Scal|}$ and recall the definition of $\xi_\Scal$ in \eqref{eq.key_value}.  Then, \eqref{ineq.i} being strictly less than one is equivalent to $\tfrac{\kappa_{\nabla c}}{\sigma_{\min}} \sqrt{\gamma_{\nabla c} t(t-1)} < 1$, which in turn is equivalent to $t(t-1)<\tfrac{\sigma_{\min}^2}{\kappa_{\nabla c}^2\gamma_{\nabla c}}$.  The left-hand side of this latter inequality is a quadratic function of~$t$, so it follows that it is satisfied for all positive $t < \thalf + \thalf \sqrt{1+\tfrac{4\sigma_{\min}^2}{\kappa_{\nabla c}^2\gamma_{\nabla c}}}$.  Substituting $\tfrac{N}{|\Scal|}$ for~$t$, one finds after rearrangement that this is equivalent to~\eqref{ineq.acute.perturb}.  Hence, by \eqref{ineq.acute.perturb}, one has that \eqref{ineq.i} yields $\| \nabla c(x)^\dag (\nabla c_\Scal(x) - \nabla c (x) )\| < 1$.  Now observe that from \cite[Theorem 6.6]{Stew1973} (or see \cite[Eq.~(3)]{Stew1979}) that
  \begin{align*}
    &\ \sigma_m (I) - \sigma_m (I + \nabla c(x)^\dag (\nabla c_\Scal(x) - \nabla c (x) )) \\
    \leq&\ | \sigma_m (I) - \sigma_m (I + \nabla c(x)^\dag (\nabla c_\Scal(x) - \nabla c (x) )) | \\
    \leq&\ \|\nabla c(x)^\dag (\nabla c_\Scal(x) - \nabla c (x) )\| < 1,
  \end{align*}
  from which it follows that $\sigma_m (I + \nabla c(x)^\dag (\nabla c_\Scal(x) - \nabla c (x) )) > 0$.  Thus, the matrix $I + \nabla c(x)^\dag (\nabla c_\Scal(x) - \nabla c (x) )$ is nonsingular, which along with the fact that $\nabla c(x)$ has full-column rank means that the matrix in \eqref{eq.nabla.c.acute.formula} has full rank, namely, $m$.
\eproof

The following remark should be emphasized at this point in our analysis.

\bremark\label{rem.loose}
  Condition~\eqref{ineq.acute.perturb} implies that for the result of Lemma~\ref{lemma.acute.perturb} to hold, the sample size $|\Scal|$ needs to be proportional to the full number of samples $N$.  This is a consequence of Lemma~\ref{lemma.sample.average.result}, which in turn relies on the relatively loose Assumption~\ref{ass.bounded.distribute}. However, this fact should not be seen as a limitation of our algorithm and our analysis of it.  One could obtain a looser lower bound than in condition~\eqref{ineq.acute.perturb} for the sample size---to obtain the same conclusion with high probability---if one were to replace Assumption~\ref{ass.bounded.distribute} with a stronger assumption on the distribution of $\omega_c$. A similar comment can be made about our next lemma as well, in terms of condition~\eqref{ineq.y.y-ys} and the distribution of $\omega_f$. No matter what type of trade-off is chosen between the assumption about the random variables and the resulting lower bound for the sample size, our results would still ultimately rely on the full-sampled and sampled derivatives being acute perturbations of each other.  It should also be said that even though a stronger assumption than Assumption~\ref{ass.bounded.distribute} that relies on restrictions on the distributions of $\omega_f$ and $\omega_c$ in order to state high-probability bounds would also be interesting for theoretical purposes, such an assumption would generally not be verifiable in practice.
\eremark

Our next lemma shows that with a similar lower bound on the sample size one can bound differences between sampled and full-sample least-squares multipliers.

\begin{lemma}\label{lemma.bound.y.y_p}
  If nonempty $\Scal\subseteq[N]$ satisfies $|\Scal| \geq p_1$ and
  \bequation\label{ineq.y.y-ys}
    |\Scal| \geq \tfrac{2N}{1+\sqrt{1+\tfrac{4\sigma_{\min}^2}{9\gamma_{\nabla c} \kappa_{\nabla c}^2}}} \in (0,N],
  \eequation
  then for all $x \in \R{n}$ one has $($recall \eqref{eq.lsm}$)$ that
  \bequationNN
    \|y(x)-y_\Scal(x)\| \leq \tfrac{\xi_\Scal (9\kappa_{\nabla f} \sqrt{\gamma_{\nabla c}} \kappa_{\nabla c} + 3 \sqrt{\gamma_{\nabla f}} \sigma_{\min})}{2\sigma_{\min}^2}.
  \eequationNN
 \end{lemma} 
\bproof
  First, let us note that from \eqref{eq.lsm} and the triangle inequality one has
  \begin{align}
    &\ \|y(x) - y_\Scal(x)\| \nonumber \\
    =&\ \|\nabla c(x)^\dag \nabla f(x) - \nabla c_\Scal(x)^\dag \nabla f_\Scal(x) \| \nonumber \\
    \leq&\ \|\nabla c(x)^\dag \nabla f(x) - \nabla c_\Scal(x)^\dag \nabla f(x)\| + \|\nabla c_\Scal(x)^\dag \nabla f(x) - \nabla c_\Scal(x)^\dag \nabla f_\Scal(x)\|. \label{eq.iggy}
  \end{align}
  The desired conclusion will follow by bounding these two right-hand-side terms uniformly over all $x \in \R{n}$. One can accomplish this by employing Theorem~5.2 and Corollary~3.9 from \cite{Stew1977}. Toward this end, let us remark that~\eqref{ineq.y.y-ys} implies that, for all $x \in \R{n}$, the matrices $\nabla c(x)$ and $\nabla c_\Scal(x)$ are acute perturbations of each other. This is straightforward since the right-hand side of \eqref{ineq.y.y-ys} is larger than that of~\eqref{ineq.acute.perturb}.  Hence, when \eqref{ineq.y.y-ys} holds, \eqref{ineq.acute.perturb} also holds, in which case one has from Lemma~\ref{lemma.acute.perturb} that $\nabla c (x)$ and $\nabla c_\Scal(x)$ are acute perturbations of each other for all $x \in \R{n}$, as desired.

  Consider arbitrary $x \in \R{n}$. Let us now introduce the quantities that are needed to apply the results from \cite{Stew1977}. A singular value decomposition of $\nabla c(x)$ yields
  \bequationNN
    \nabla c(x) = U \bbmatrix S \\ 0 \ebmatrix V^T \implies \bbmatrix S \\ 0 \ebmatrix = U^T \nabla c(x) V,
  \eequationNN
  where $U \in \R{n \times n}$ and $V \in \R{m \times m}$ are orthogonal, and where under Assumption~\ref{ass.subproblems} the diagonal matrix $S \in \R{m \times m}$ is positive definite.  Correspondingly, let us define
  \bequationNN
    \bbmatrix E_1 \\ E_2 \ebmatrix := U^T (\nabla c_\Scal(x) - \nabla c(x)) V\ \ \text{and}\ \ \bbmatrix B_1 \\ B_2 \ebmatrix := U^T \nabla c_\Scal(x) V,
  \eequationNN
  where $E_1 \in \R{m \times m}$, $E_2 \in \R{(n-m) \times m}$, $B_1 \in \R{m \times m}$, and $B_2 \in \R{(n-m) \times m}$.  Note that our tuple $(S,E_1,E_2,B_1,B_2)$ corresponds to the tuple $(A_{11},E_{11},E_{21},B_{11},B_{21})$ that is introduced in \cite[page~636]{Stew1977}, and note that since our Assumption~\ref{ass.subproblems} ensures that $\nabla c(x)$ has full-column rank, the matrices $(E_{12},E_{22},B_{12},B_{22})$ that are introduced in \cite[page~636]{Stew1977} are not present in our setting.  Furthermore, for bounding the former term in \eqref{eq.iggy}, let us define for $b_1 \in \R{m}$ and $b_2 \in \R{n-m}$ the matrix
  \bequationNN
    \bbmatrix b_1 \\ b_2 \ebmatrix := U^T(-\nabla f(x)).
  \eequationNN
  
  Let us now bound the first term on the right-hand side of \eqref{eq.iggy}. Recalling that $y(x) = -\nabla c (x)^\dag \nabla f(x)$, one can state that \cite[Theorem 5.2]{Stew1977} in our setting yields
  \bequation\label{ineq.sss}
    \tfrac{\|y(x) + \nabla c_\Scal(x)^\dag \nabla f(x)\|}{\|y(x)\|} \leq \bar{\kappa} \tfrac{\|E_1\|}{\|\nabla c(x)\|} + \bar{\kappa}^2 \tfrac{\|E_2\|}{\|\nabla c(x)\|} \(\eta \tfrac{\|b_2\|}{\|b_1\|} + \tfrac{\|E_2\|}{\|\nabla c(x)\|}\),
  \eequation
  where, as in \cite[Theorem 3.8]{Stew1977} and the displayed equation prior to \cite[Theorem 5.1]{Stew1977},
  \bequationNN
    \bar{\kappa} := \|\nabla c(x)\| \| B_1^{-1} \| = \|\nabla c(x)\| \| (S + E_1)^{-1} \|\ \ \text{and}\ \ \eta := \tfrac{\|b_1\|}{\|\nabla c(x)\|\|y(x)\|}.
  \eequationNN
  We remark in passing that our statement of \eqref{ineq.sss} corrects a typo in the statement of \cite[Theorem 5.2]{Stew1977}.  In particular, in the statement of \cite[Theorem 5.2]{Stew1977}, the latter term on the right-hand side is stated with $\|E_{12}\|$ in the numerator outside of the parentheses.  That is a typo.  One can see through the proof of \cite[Theorem 5.2]{Stew1977} that our statement is correct, where in fact $\|E_{21}\|$ belongs in this numerator.  In our setting, this corresponds to $\|E_2\|$ in this numerator, as we have stated in \eqref{ineq.sss}.

  Our next goal is to prove upper bounds for the terms on the right-hand side of \eqref{ineq.sss}.  First, with respect to $\|b_2\|$, one has from submultiplicity of the matrix 2-norm, $\|U\| = 1$ due to $U$ being orthogonal, and Assumption \ref{ass.boundness} that
  \bequation\label{ineq.b}
    \|b_2\| \leq \|U^T(-\nabla f(x))\| \leq \|U\| \|\nabla f(x)\| \leq \kappa_{\nabla f}.
  \eequation
  Second, with respect to both $\|E_1\|$ and $\|E_2\|$, one has from submultiplicity of the matrix 2-norm and $\|U\| = \|V\| = 1$ due to both $U$ and $V$ being orthogonal that
  \bequation\label{ineq.E}
    \max\{\|E_1\|,\|E_2\|\} \leq \left\| \bbmatrix E_1 \\ E_2 \ebmatrix \right\| = \|U^T(\nabla c_\Scal(x) - \nabla c(x))V\| \leq \|\nabla c_\Scal(x) - \nabla c(x)\|.
  \eequation
  Third, let us bound $\| (S + E_1)^{-1} \|^2$.  For this, let us employ \cite[Theorem~2.2]{Stew1977}, where now our tuple of matrices $(S,E_1,B_1)$ play the role of $(A,E,B)$ in that theorem.  In particular, let us employ the second part of \cite[Theorem 2.2]{Stew1977}, which requires (in our notation) $S$ to be nonsingular and $\|S^{-1}\| \|E_1\| < 1$.  The fact that $S$ is nonsingular has been established earlier through the fact that $S$ is positive definite.  On the other hand, let us now show $\|S^{-1}\| \|E_1\| < 1$.  Toward this end, observe that
  \begin{align*}
    \nabla c(x)^\dag
    &= (\nabla c(x)^T\nabla c(x))^{-1} \nabla c(x)^T \\
    &= \(V \bbmatrix S \\ 0 \ebmatrix^T U^T U \bbmatrix S \\ 0 \ebmatrix V^T\)^{-1} V \bbmatrix S \\ 0 \ebmatrix^T U^T \\
    &= (V S^2 V^T)^{-1} V \bbmatrix S & 0 \ebmatrix U^T \\
    &= (V S^{-2} V^T) V \bbmatrix S & 0 \ebmatrix U^T = V \bbmatrix S^{-1} & 0\ebmatrix U^T.
  \end{align*}
  That is, the right-hand side of the above equation is an SVD for $\nabla c(x)^\dag$, from which it follows that $\|S^{-1}\| = \|\nabla c(x)^\dag\|$.  At the same time, note that using a similar argument as in Lemma~\ref{lemma.acute.perturb}, it follows that \eqref{ineq.y.y-ys} implies that the right-hand side of~\eqref{ineq.i} is less than or equal to $1/3$, which shows that
  \bequation\label{ineq.bound.psedueC.C}
    \| \nabla c(x)^\dag\|\| (\nabla c_\Scal(x) - \nabla c (x) )\| \leq \tfrac{\xi_\Scal \sqrt{\gamma_{\nabla c}} \kappa_{\nabla c}}{\sigma_{\min}} \leq \tfrac{1}{3}.
  \eequation
  Hence, with \eqref{ineq.E} and \eqref{ineq.bound.psedueC.C}, one finds that
  \bequation\label{ineq.1/3}
    \|S^{-1}\| \|E_1\| = \|\nabla c(x)^{\dag}\| \|E_1\| \leq \|\nabla c(x)^{\dag}\| \|\nabla c_\Scal(x) - \nabla c(x)\| \leq \tfrac{1}{3}. 
  \eequation
  Since the requirements of \cite[Theorem 2.2]{Stew1977} thus hold, one obtains with \eqref{ineq.1/3} that
  \bequation\label{ineq.sssss}
    \|(S + E_1)^{-1}\| = \|B_1^{-1}\| \leq \tfrac{\|S^{-1}\|}{1 - \|S^{-1}\| \|E_1\|} \leq \tfrac{\|S^{-1}\|}{1-\tfrac{1}{3}} = \tfrac{3}{2}\|\nabla c (x)^\dag\|.
  \eequation
  Combining \eqref{lemma.bound.nabla_c_dag}, \eqref{ineq.bound.psedueC.C}, \eqref{ineq.sss}, \eqref{ineq.b}, \eqref{ineq.E}, and \eqref{ineq.sssss}, one now obtains
  \begin{align*}
    \tfrac{\|y(x) + \nabla c_\Scal(x)^\dag \nabla f(x)\|}{\|y(x)\|}
    \leq&\ \|E_1\| \|(S + E_1)^{-1}\| + \tfrac{\|(S + E_1)^{-1} \|^2 \|E_2\| \|b_2\|}{\|y(x)\|} \\
    &\ + \|(S + E_1)^{-1}\|^2 \|E_2\|^2 \\
    \leq&\ \(\tfrac{3}{2} \|\nabla c(x)^\dag\| + \tfrac{9}{4} \tfrac{\|\nabla c(x)^\dag\|^2 \kappa_{\nabla f}}{\|y(x)\|}\) \|\nabla c_\Scal(x) - \nabla c(x)\| \\
    &\ + \tfrac{9}{4} \|\nabla c(x)^\dag\|^2 \|\nabla c_\Scal(x) - \nabla c(x)\|^2 \\
    \leq&\ \(\tfrac{9}{4} \|\nabla c(x)^\dag\| + \tfrac{9}{4} \tfrac{\|\nabla c(x)^\dag\|^2 \kappa_{\nabla f}}{\|y(x)\|}\) \|\nabla c_\Scal(x) - \nabla c(x)\| \\
    \leq&\ \(\tfrac{9}{4\sigma_{\min}} +\tfrac{9\kappa_{\nabla f}}{4\sigma_{\min}^2\|y(x)\|}\)\|\nabla c(x)-\nabla c_\Scal(x)\| .
  \end{align*}
  Multiplying $\|y(x)\|$ on both sides of the above inequality, one finds with \eqref{lemma.bound.y_N} that
  \begin{align*}
    \|y(x) + \nabla c_\Scal(x)^\dag \nabla f(x)\|
      &\leq \(\tfrac{9\|y(x)\|}{4\sigma_{\min}} +\tfrac{9\kappa_{\nabla f}}{4\sigma_{\min}^2}\)\|\nabla c_\Scal(x) - \nabla c(x)\| \\
      &\leq \tfrac{9\kappa_{\nabla f}}{2\sigma_{\min}^2}\|\nabla c_\Scal(x) - \nabla c(x)\|.
  \end{align*}
  Further, using Assumptions \ref{ass.boundness} and Lemma~\ref{lemma.sample.average.result}, one finds that
  \bequation\label{eq.iggy2}
    \|y(x) + \nabla c_\Scal(x)^\dag \nabla f(x)\| \leq \tfrac{9\kappa_{\nabla f} \xi_\Scal \sqrt{\gamma_{\nabla c}}}{2\sigma_{\min}^2} \|\nabla c(x)\| \leq \tfrac{9\kappa_{\nabla f} \xi_\Scal \sqrt{\gamma_{\nabla c}} \kappa_{\nabla c}}{2\sigma_{\min}^2}.
  \eequation

  Now let us bound the second term in \eqref{eq.iggy}.  One finds with Lemma~\ref{lemma.sample.average.result} that
  \begin{align*}
    \|\nabla c_\Scal(x)^\dag \nabla f(x) - \nabla c_\Scal(x)^\dag \nabla f_\Scal(x)\| &\leq \|\nabla c_\Scal(x)^\dag\| \|\nabla f(x) - \nabla f_\Scal(x)\| \\
    &\leq \|\nabla c_\Scal(x)^\dag\| \xi_\Scal \sqrt{\gamma_{\nabla f}}.
  \end{align*}
  Furthermore, one has with \eqref{ineq.bound.psedueC.C} that the conditions of \cite[Corollary~3.9]{Stew1977} hold, so
  \bequationNN
    \|\nabla c_\Scal(x)^\dag\| \leq \tfrac{\|\nabla c(x)^\dag\|}{1 - \| \nabla c(x)^\dag \|\|E_1\|} \leq \tfrac32 \|\nabla c(x)^\dag\| \leq \tfrac{3}{2 \sigma_{\min}}.
  \eequationNN
  Consequently, one finds from above that
  \bequationNN
    \|\nabla c_\Scal(x)^\dag \nabla f(x) - \nabla c_\Scal(x)^\dag \nabla f_\Scal(x)\| \leq \tfrac{3 \xi_\Scal \sqrt{\gamma_{\nabla f}}}{2 \sigma_{\min}}
  \eequationNN
  Combining this with \eqref{eq.iggy2}, the proof is complete since both of the derived bounds are independent of $x \in \R{n}$.
\eproof

Our next goal is to prove a bound on the difference of inner products involving the Hessians of the Lagrangians $L_\Scal$ and $L$ when $|\Scal|$ is sufficiently large.  Toward this end, we first present the following preliminary lemma, which for any $x \in \R{n}$ shows a bound on the product between any element of the null space of the sample constraint Jacobian $\nabla c_\Scal(x)$ and a projection matrix that maps vectors onto the span of the columns of the full-sample constraint Jacobian $\nabla c(x)$.

\begin{lemma}\label{lemma.proj.apply}
  If nonempty $\Scal \subseteq [N]$ satisfies $|\Scal| \geq p_1$ and $\xi_\Scal \sqrt{\gamma_{\nabla c}} \kappa_{\nabla c} < \sigma_{\min}$, then for all $x \in \R{n}$ and for all $d_\Scal \in \Null(\nabla c_\Scal(x)^T) \setminus \{0\}$ one has that
  \bequationNN
    \tfrac{\|\Rcal(\nabla c(x)) d_\Scal\|}{\|d_\Scal\|} \leq \tfrac{\xi_\Scal \sqrt{\gamma_{\nabla c}} \kappa_{\nabla c}}{\sigma_{\min}} < 1.
  \eequationNN
\end{lemma}
\bproof
  Consider arbitrary $x \in \R{n}$ and $d_\Scal \in \Null(\nabla c_\Scal(x)^T) \setminus \{0\}$.  It follows from the definitions of $\Rcal$ and $\Ncal$ that $\Rcal(\nabla c_\Scal(x)) d_\Scal = 0$ and $\Ncal(\nabla c_\Scal(x)^T) d_\Scal = d_\Scal$.  Let us now consider $\Rcal(\nabla c(x))d_\Scal$.  By submultiplicity of the matrix 2-norm, one has
  \bequationNN
    \|\Rcal(\nabla c(x)) d_\Scal\| = \|\Rcal(\nabla c(x)) \Ncal(\nabla c_\Scal(x)^T) d_\Scal\| \leq \|\Rcal(\nabla c(x)) \Ncal(\nabla c_\Scal(x)^T)\| \|d_\Scal\|.
  \eequationNN
  Now, since by Assumption \ref{ass.subproblems} the matrices $\nabla c_\Scal(x)$ and $\nabla c(x)$ have the same rank, it follows from \cite[Theorem 2.4]{Stew1977} that
  \begin{align}
    \|\Rcal(\nabla c(x))d_\Scal\|
    &\leq \|\Rcal(\nabla c(x))\Ncal(\nabla c_\Scal(x)^T)\|\|d_\Scal\| \nonumber \\
    &\leq \min\{ \|\nabla c(x)^\dag\|, \|\nabla c_\Scal(x)^\dag\| \} \|\nabla c(x) -\nabla c_\Scal(x)\| \|d_\Scal\| \nonumber \\
    &\leq \|\nabla c(x)^\dag\| \|\nabla c(x) - \nabla c_\Scal(x)\| \|d_\Scal\|. \label{ineq.Ra.Rb}
  \end{align}
  Since $d_\Scal \neq 0 \in \R{n}$, by dividing the expression above by $\|d_\Scal\|$, the resulting expression can be bounded above as in \eqref{ineq.i}, which yields the desired inequality.
\eproof

With the preceding lemma, we can now prove that if the sample set $\Scal$ is sufficiently large in cardinality, then~\eqref{prob.opt.S} is $(\alpha_\Scal,\beta_\Scal)$-strongly Morse with respect to a pair $(\alpha_\Scal,\beta_\Scal) \in \R{2}_{>0}$ defined with respect to $(\alpha,\beta) \in \R{2}_{>0}$ from Assumption~\ref{ass.subproblems}.

\blemma\label{lem: emprical morse}
  Define the tuple $(\kappa_1, \kappa_2, \kappa_3) \in \R{}_{>0} \times \R{}_{>0} \times \R{}_{>0}$ by
  \begin{align}
    &\ \kappa_1 := \tfrac{\sqrt{\gamma_{\nabla c}} \kappa_{\nabla c}}{\sigma_{\min}},\ \ \kappa_2 := \sqrt{(2 \sqrt{\kappa_{\nabla f}} + \kappa_1 \kappa_{\nabla f})^2 + \gamma_c},\ \ \text{and}\ \ \label{eq.kappa} \\
    &\ \kappa_3 := \(3\kappa_{\nabla^2 f}\kappa_1 + \sqrt{\kappa_{\nabla^2 f}} + \tfrac{\sqrt{m} (3 \kappa_{\nabla^2 c} \kappa_1 + \sqrt{\gamma_{\nabla^2 c}}) (5 \kappa_{\nabla f} \kappa_1 + \sqrt{\gamma_{\nabla f}})}{2 \sigma_{\min} \kappa_1}\). \label{eq.def.kappa.3}
  \end{align}
  Then, with $(\alpha,\beta) \in \R{2}_{>0}$ from Assumption~\ref{ass.subproblems}, if nonempty $\Scal \subseteq [N]$ has $|\Scal| \geq p_1$ and 
  \bequation\label{ineq.theorem1.S}
    \xi_\Scal \leq \min\left\{\tfrac{1}{3\kappa_1}, \tfrac{\alpha}{2\kappa_2}, \tfrac{7\beta}{18\kappa_3} \right\},
  \eequation
  then with respect to
  \bequation\label{eq.alpha-beta-s}
    \alpha_\Scal := \alpha - \xi_\Scal \kappa_2 \geq \thalf\alpha > 0\ \ \text{and}\ \ \beta_\Scal := (1-\tfrac{1}{3} \kappa_1 \xi_\Scal) \beta - \kappa_3 \xi_\Scal \geq \thalf\beta > 0
  \eequation
  subproblem \eqref{prob.opt.S} is $(\alpha_\Scal,\beta_\Scal)$-strongly Morse in the sense of Assumption~\ref{ass.subproblems}; in other words, one has for any $x \in \R{n}$ that if $(x,y_\Scal(x))$ satisfies $\|\nabla L(x,y_\Scal(x))\| \leq \alpha_\Scal$, then
  \bequationNN
    |d^T \nabla_{xx}^2 L_\Scal(x,y_\Scal(x)) d| \geq \beta_\Scal \|d\|_2^2\ \ \text{for all}\ \ d \in \Null(\nabla c_\Scal(x)^T).
  \eequationNN
\elemma
\bproof
  Our first aim is to show that with $(\kappa_1,\kappa_2)$ defined in \eqref{eq.kappa} one finds for any point $x \in \R{n}$ and nonempty $\Scal \subseteq [N]$ satisfying $|\Scal| \geq p_1$ that $\|\nabla L_\Scal(x,y_\Scal(x))\| \leq \alpha_\Scal$ implies that $\|\nabla L(x,y(x))\| \leq \alpha$ holds.  Toward this end, consider arbitrary $x \in \R{n}$ and nonempty $\Scal \subseteq [N]$ satisfying $|\Scal| \geq p_1$ with $\|\nabla L_\Scal(x,y_\Scal(x))\| \leq \alpha_\Scal$, where~$(\kappa_1,\kappa_2,\alpha_\Scal)$ is defined by \eqref{eq.kappa} and \eqref{eq.alpha-beta-s}.  Our first aim is to prove a bound on the norm of the difference between the gradient of the Lagrangian with respect to the sample set~$\Scal$ and the gradient of the Lagrangian with respect to the full-sample set.  First, observe that by~\eqref{eq.soc.S.1}, \eqref{eq.lsm}, and the triangle inequality, it follows that
  \begin{align}
    &\ \|\nabla_x L(x,y(x)) - \nabla_x L_\Scal(x,y_{\Scal}(x)) \| \nonumber \\
    \leq&\ \|\nabla f(x) - \nabla f_\Scal(x)\| + \|\nabla c(x) y(x) - \nabla c_\Scal(x) y_\Scal(x) \| \nonumber \\
    =&\ \|\nabla f(x) - \nabla f_\Scal(x)\| + \|\nabla c_\Scal(x) \nabla c_\Scal(x)^\dag \nabla f_\Scal(x) - \nabla c(x) \nabla c(x)^\dag \nabla f(x)\|. \label{ineq.L.grad.x.diff.1}
  \end{align}
  Along with submultiplicity of the matrix 2-norm, the latter term satisfies
  \begin{align*}
    &\ \|\nabla c_\Scal(x) \nabla c_\Scal(x)^\dag \nabla f_\Scal(x) - \nabla c(x) \nabla c(x)^\dag \nabla f(x)\|\\
    =&\ \|\Rcal(\nabla c_\Scal(x)) \nabla f_\Scal(x) - \Rcal(\nabla c(x)) \nabla f(x)\| \\
    =&\ \|\Rcal(\nabla c_\Scal(x)) \nabla f_\Scal(x) - \Rcal(\nabla c_\Scal(x)) \nabla f(x) + \Rcal(\nabla c_\Scal(x)) \nabla f(x) - \Rcal(\nabla c(x)) \nabla f(x)\| \\
    \leq&\ \|\Rcal(\nabla c_\Scal(x))\| \|\nabla f_\Scal(x) - \nabla f(x)\| + \|\Rcal(\nabla c_\Scal(x)) - \Rcal(\nabla c(x))\| \|\nabla f(x)\|.
  \end{align*}
  Combined with \eqref{ineq.L.grad.x.diff.1} and since $\Rcal(\nabla c_\Scal(x))$ is a projection matrix, one obtains that
  \begin{align*}
    &\ \|\nabla_x L(x,y(x)) - \nabla_x L_\Scal(x,y_{\Scal}(x)) \| \nonumber \\
    \leq&\ (1 + \|\Rcal(\nabla c_\Scal(x))\|) \|\nabla f_\Scal(x) - \nabla f(x)\| + \|\Rcal(\nabla c_\Scal(x)) - \Rcal(\nabla c(x))\| \|\nabla f(x)\| \\
    \leq&\ 2 \|\nabla f_\Scal(x) - \nabla f(x)\| + \|\Rcal(\nabla c_\Scal(x)) - \Rcal(\nabla c(x))\| \|\nabla f(x)\|.
  \end{align*}
  Since by Assumption \ref{ass.subproblems} the Jacobians $\nabla c_\Scal(x)$ and $\nabla c(x)$ have full-column rank, by \cite[Theorem 2.3]{Stew1977}, \cite[Theorem 2.4]{Stew1977}, \eqref{lemma.bound.nabla_c_dag}, \eqref{lemma.saa.nabla_f}, and \eqref{lemma.saa.nabla_c} one then has
  \begin{align}
    &\ \|\nabla_x L(x,y(x)) - \nabla_x L_{\Scal}(x,y_{\Scal}(x))\| \nonumber \\
    \leq&\ 2 \|\nabla f_\Scal(x) - \nabla f(x)\| + \|\Rcal(\nabla c(x)) \Ncal(\nabla c_\Scal(x)^T)\| \kappa_{\nabla f} \nonumber \\
    \leq&\ 2 \|\nabla f_\Scal(x) - \nabla f(x)\| + \|\nabla c(x)^\dag\| \|\nabla c_\Scal(x) - \nabla c(x)\| \kappa_{\nabla f} \nonumber \\
    \leq&\ 2 \xi_\Scal \sqrt{\kappa_{\nabla f}} + \tfrac{\xi_\Scal \sqrt{\gamma_{\nabla c}} \kappa_{\nabla c} \kappa_{\nabla f}}{\sigma_{\min}} = \xi_\Scal (2 \sqrt{\kappa_{\nabla f}} + \kappa_1 \kappa_{\nabla f}). \label{ineq.L.grad.x.diff}
  \end{align}
  Second, by \eqref{eq.soc.S.1} and \eqref{lemma.saa.c}, one finds that
  \bequation\label{ineq.L.grad,y.diff}
    \|\nabla_y L(x,y(x)) - \nabla_y L_\Scal(x,y_{\Scal}(x))\| = \|c(x) -  c_\Scal(x)\| \leq \xi_\Scal \sqrt{\gamma_c}.
  \eequation
  Combining \eqref{ineq.L.grad.x.diff} and \eqref{ineq.L.grad,y.diff}, one finds that
  \bequationNN
    \|\nabla L(x,y(x)) - \nabla L_{\Scal}(x,y_{\Scal}(x))\| \leq \xi_\Scal \sqrt{(2 \sqrt{\kappa_{\nabla f}} + \kappa_1 \kappa_{\nabla f})^2 + \gamma_c} = \xi_\Scal \kappa_2.
  \eequationNN
  Thus, since $\|\nabla L_\Scal(x,y_\Scal(x))\| \leq \alpha_\Scal = \alpha - \xi_\Scal \kappa_2$, as desired one obtains
  \begin{align*}
    \|\nabla L(x,y(x))\|
    &= \|\nabla L_\Scal(x,y_\Scal(x)) + \nabla L(x,y(x)) - \nabla L_{\Scal}(x,y_{\Scal}(x))\| \\
    &\leq \|\nabla L_\Scal(x,y_\Scal(x))\| + \|\nabla L(x,y(x)) - \nabla L_{\Scal}(x,y_{\Scal}(x))\| \\
    &\leq \|\nabla L_\Scal(x,y_\Scal(x))\| + \xi_\Scal \kappa_2 \leq \alpha - \xi_\Scal \kappa_2 + \xi_\Scal \kappa_2 = \alpha.
  \end{align*}
  Finally, observe that $\alpha_\Scal := \alpha - \xi_\Scal \kappa_2 \geq \thalf \alpha$ follows since \eqref{ineq.theorem1.S} implies $\xi_\Scal \leq \alpha/(2\kappa_2)$.

  We have shown that for any $x \in \R{n}$ and with $(\kappa_1,\kappa_2)$, nonempty $\Scal \subseteq [N]$ satisfying $|\Scal| \geq p_1$, and $\alpha_\Scal$ satisfying the conditions of the lemma, one finds that $\|\nabla L_\Scal(x,y_\Scal(x))\| \leq \alpha_\Scal$ yields $\|\nabla L(x,y(x))\| \leq \alpha$.  Our next aim is to show that, for such $x$ and $\Scal$ under the conditions of the lemma, one has $|d_\Scal^T \nabla_{xx}^2 L_\Scal(x,y_\Scal(x)) d_\Scal| \geq \beta_\Scal \|d_\Scal\|_2^2$ for all nonzero $d_\Scal \in \Null(\nabla c_\Scal(x)^T)$.  (Observe that $|d_\Scal^T \nabla_{xx}^2 L_\Scal(x,y_\Scal(x)) d_\Scal| \geq \beta_\Scal \|d_\Scal\|_2^2$ holds trivially for $d_\Scal = 0$.)  Consider such $x$ and $\Scal$, and note that by Assumption~\ref{ass.subproblems} it follows that since $\|\nabla L(x,y(x))\| \leq \alpha$ one also has $|d^T\nabla^2_{xx} L(x,y(x))d| \geq \beta \|d\|^2$ for all $d \in \Null(\nabla c(x)^T)$.  Consider any nonzero $d_\Scal \in \Null(\nabla c_\Scal(x)^T)$.  First,
  \begin{align}
    &\ |d_\Scal^T \nabla^2_{xx} L_\Scal(x,y_\Scal(x)) d_\Scal| \nonumber \\
    =&\ |d^T_\Scal \nabla^2_{xx} L(x,y(x)) d_\Scal + d_\Scal^T \nabla^2_{xx} L_\Scal(x,y_\Scal(x)) d_\Scal -d^T_\Scal \nabla^2_{xx} L (x,y(x)) d_\Scal | \nonumber \\
    \geq&\ |d_\Scal^T \nabla^2_{xx} L(x,y(x)) d_\Scal | - |d_\Scal^T \nabla^2_{xx} L_\Scal(x,y_\Scal(x)) d_\Scal - d_\Scal^T \nabla^2_{xx} L(x,y(x)) d_\Scal|. \label{ineq.dSLdS}
  \end{align}
  Let us now derive lower bounds on the two terms on the right-hand side of \eqref{ineq.dSLdS}.  Note that $d_\Scal$ can be decomposed as the sum of two orthogonal vectors, namely, $d_\Scal = \dbar_\Scal + r_\Scal$, where $\dbar_\Scal = \Ncal(\nabla c(x)^T)d_\Scal$ and $r_\Scal = \Rcal(\nabla c(x)) d_\Scal$.  It follows that
  \bequation\label{eq.orth.d.r.}
    \|d_\Scal\|^2 = \|\dbar_\Scal + r_\Scal\|^2 = \|\dbar_\Scal\|^2 + 2\dbar_\Scal^Tr_\Scal + \|r_\Scal\|^2 = \|\dbar_\Scal\|^2 + \|r_\Scal\|^2
  \eequation
  along with $\max\{\|\dbar_\Scal\|, \|r_\Scal\|\} \leq \|d_\Scal\|$.  Now, one finds with the triangle inequality, Cauchy-Schwarz inequality, submultiplicity of the matrix 2-norm, the fact that problem~\eqref{prob.opt.N} is $(\alpha,\beta)$-strongly Morse, and $\max\{\|\dbar_\Scal\|, \|r_\Scal\|\} \leq \|d_\Scal\|$ that
  \begin{align}
    &\ |d_\Scal^T \nabla^2_{xx} L(x,y(x)) d_\Scal| \nonumber \\
    =&\ |(\dbar_\Scal + r_\Scal)^T \nabla^2_{xx} L(x,y(x)) (\dbar_\Scal + r_\Scal) | \nonumber \\
    =&\ |\dbar_\Scal^T \nabla^2_{xx} L(x,y(x)) \dbar_\Scal + 2\dbar_\Scal^T \nabla^2_{xx} L(x,y(x)) r_\Scal + r_\Scal^T \nabla^2_{xx} L(x,y(x)) r_\Scal| \nonumber \\
    \geq&\ |\dbar_\Scal^T \nabla^2_{xx} L(x,y(x)) \dbar_\Scal| - |2\dbar_\Scal^T \nabla^2_{xx} L(x,y(x)) r_\Scal | -|r_\Scal^T \nabla^2_{xx} L(x,y(x)) r_\Scal| \nonumber \\
    \geq&\ |\dbar_\Scal^T \nabla^2_{xx} L(x,y(x)) \dbar_\Scal | - 2\|\nabla^2_{xx} L(x,y(x))\| \|\dbar_\Scal\| \|r_\Scal\| \nonumber \\
    &\ - \|\nabla^2_{xx} L(x,y(x))\| \|r_\Scal\|^2 \nonumber \\
    \geq&\ \beta\|\dbar_\Scal\|^2 - 3\|\nabla^2_{xx} L(x,y(x))\| \|d_\Scal\| \|r_\Scal\|. \label{ineq.dSldS}
  \end{align}
  Next, let us derive an upper bound for $\|r_\Scal\|$ with respect to $\|d_\Scal\|$.  Recall that $r_\Scal = \Rcal(\nabla c(x))d_\Scal$.  Our aim is to employ Lemma~\ref{lemma.proj.apply}, for which $\Scal$ needs to satisfy $\xi_\Scal \sqrt{\gamma_{\nabla c}} \kappa_{\nabla c} < \sigma_{\min}$.  To see that this holds, note that by \eqref{ineq.theorem1.S} one has
  \bequationNN
    \tfrac{\xi_\Scal \sqrt{\gamma_{\nabla c}} \kappa_{\nabla c}}{\sigma_{\min}} = \kappa_1 \xi_\Scal \leq \tfrac{1}{3} < 1.
  \eequationNN
  Now applying the result of Lemma~\ref{lemma.proj.apply} and since $\|d_\Scal\| \neq 0$, one finds that
  \bequationNN
    \tfrac{\|r_\Scal\|}{\|d_\Scal\|}=\tfrac{\|\Rcal(\nabla c(x))d_\Scal\|}{\|d_\Scal\|} \leq  \tfrac{\xi_\Scal \sqrt{\gamma_{\nabla c}} \kappa_{\nabla c}}{\sigma_{\min}} = \kappa_1 \xi_\Scal.
  \eequationNN
  Multiplying by $\|d_\Scal\|$ on both sides yields $\|r_\Scal\| \leq \kappa_1 \xi_\Scal \|d_\Scal\|$.  Combined with \eqref{eq.orth.d.r.}, this yields $\|d_\Scal\|^2 \le \|\dbar_\Scal\|^2 + \kappa_1^2 \xi_\Scal^2 \|d_\Scal\|^2$, which with \eqref{ineq.theorem1.S} implies
  \bequation\label{ineq.bound.dbars}
    \|\dbar_\Scal\|^2 \geq \(1-\kappa_1^2 \xi_\Scal^2\) \|d_\Scal\|^2 \geq \(1 - \tfrac{1}{3} \kappa_1 \xi_\Scal\) \|d_\Scal\|^2.
  \eequation
  Combining \eqref{ineq.bound.dbars}, \eqref{lemma.bound.L_Hess}, $\|r_\Scal\| \leq \kappa_1 \xi_\Scal \|d_\Scal\|$, and \eqref{ineq.dSldS}, it follows that
  \begin{align}
    &\ |d_\Scal^T \nabla^2_{xx} L(x,y(x)) d_\Scal| \nonumber \\
    \geq&\ \(1 - \tfrac{1}{3}\kappa_1 \xi_\Scal\) \beta \|d_\Scal\|^2 - 3 \(\kappa_{\nabla^2 f}+\tfrac{\sqrt{m} \kappa_{\nabla f} \kappa_{\nabla^2 c}}{\sigma_{\min}}\) \kappa_1 \xi_\Scal \|d_\Scal\|^2, \label{ineq.theorem.v1}
  \end{align}
  which gives a lower bound for the first term on the right-hand side of~\eqref{ineq.dSLdS}.  Now let us turn to proving a lower bound for the second term on the right-hand side of \eqref{ineq.dSLdS}, which we derive by proving an upper bound for the absolute difference $|d_\Scal^T \nabla^2_{xx}L_\Scal(x,y_\Scal(x)) d_\Scal - d_\Scal^T \nabla^2_{xx} L (x,y(x)) d_\Scal|$.  First, by the definitions of $L_\Scal$ and $L$,
  \begin{align*}
    &\ \nabla^2_{xx} L_\Scal(x,y_\Scal(x)) - \nabla^2_{xx} L (x,y(x)) \\
    =&\ \nabla^2 f_\Scal(x) - \nabla^2 f(x) + \sum_{j\in[m]} (\nabla^2 [c_\Scal]_j(x) [y_\Scal(x)]_j - \nabla^2 [c]_j(x) [y(x)]_j),
  \end{align*}
  from which it follows by submultiplicity of the matrix 2-norm that
  \begin{align}
    &\ |d_\Scal^T \nabla^2_{xx} L_\Scal(x,y_\Scal(x)) d_\Scal - d_\Scal^T \nabla^2_{xx} L (x,y(x)) d_\Scal| \nonumber \\
    \leq&\ \| \nabla^2 f_\Scal(x) - \nabla^2 f(x)\| \|d_\Scal\|^2 \nonumber \\
    &\ + \left\| \sum_{j\in[m]} (\nabla^2 [c_\Scal]_j(x) [y_\Scal(x)]_j - \nabla^2 [c]_j(x) [y(x)]_j ) \right\| \|d_\Scal\|^2. \label{ineq.usedlater1}
  \end{align}
  With respect to the first norm in \eqref{ineq.usedlater1}, one finds along with the triangle inequality, absolute homogeneity of norms, \eqref{lemma.saa.Hess_c}, and Assumption \ref{ass.subproblems} that
  \begin{align}
    &\ \left\| \sum_{j\in[m]} (\nabla^2 [c_\Scal]_j(x) [y_\Scal(x)]_j - \nabla^2 [c]_j(x) [y(x)]_j ) \right\| \nonumber \\
    \leq&\ \sum_{j\in[m]} \| \nabla^2 [c_\Scal]_j(x) - \nabla^2[c]_j(x) \| |[y_\Scal(x)]_j| \nonumber \\
    &\ + \sum_{j\in[m]} \| \nabla^2[c]_j(x) \| |[y_\Scal(x)]_j - [y(x)]_j | \nonumber \\
    \leq&\ \sum_{j\in[m]} \xi_\Scal \sqrt{\gamma_{\nabla^2 c}} |[y_\Scal(x)]_j| + \sum_{j\in[m]} \| \nabla^2[c]_j(x) \| | [y_\Scal(x)]_j - [y(x)]_j| \nonumber \\
    \leq&\ \xi_\Scal \sqrt{\gamma_{\nabla^2 c}}  \sum_{j\in[m]} |[y_\Scal(x)]_j| +  \kappa_{\nabla^2 c} \sum_{j\in[m]} | [y_\Scal(x)]_j - [y(x)]_j|. \label{ineq.usedlater2}
  \end{align}
  For the second sum in \eqref{ineq.usedlater2}, it follows from the fact that $\|v\|_1 \leq \sqrt{m} \|v\|_2$ for any $v \in \R{m}$, Lemma \ref{lemma.bound.y.y_p}, and $\kappa_1=\tfrac{\sqrt{\gamma_{\nabla c}} \kappa_{\nabla c}}{\sigma_{\min}}$ that
  \begin{align}
    \sum_{j\in[m]} |[y_\Scal(x)]_j - [y(x)]_j| 
    & \leq \sqrt{m} \|y_\Scal(x) - y(x)\| \nonumber \\
    & \leq \tfrac{\xi_\Scal \sqrt{m} (9\kappa_{\nabla f} \sqrt{\gamma_{\nabla c}} \kappa_{\nabla c} + 3 \sqrt{\gamma_{\nabla f}} \sigma_{\min})}{2\sigma_{\min}^2} = \tfrac{3 \sqrt{m} (3\kappa_{\nabla f} \kappa_1 + \sqrt{\gamma_{\nabla f}}) \xi_\Scal}{2\sigma_{\min}}, \label{ineq.usedlater3.1}
  \end{align}
  while for the first sum in \eqref{ineq.usedlater2}, by \eqref{ineq.usedlater3.1}, \eqref{lemma.bound.y_N}, and $\xi_\Scal \leq \tfrac{1}{3\kappa_1}$ one has
  \begin{align}
    \sum_{j\in[m]} | [y_\Scal(x)]_j|
    \leq&\ \sqrt{m} \|y_\Scal(x)\| \leq \sqrt{m} (\|y_\Scal(x) - y(x)\| + \|y(x)\|) \nonumber \\
    \leq&\ \tfrac{\sqrt{m} (3\kappa_{\nabla f} \kappa_1 + \sqrt{\gamma_{\nabla f}})}{2\sigma_{\min} \kappa_1} + \tfrac{\sqrt{m} \kappa_{\nabla f}}{\sigma_{\min}} \nonumber \\
    =&\ \tfrac{\sqrt{m} (5 \kappa_{\nabla f} \kappa_1 + \sqrt{\gamma_{\nabla f}})}{2 \sigma_{\min} \kappa_1}. \label{ineq.usedlater3.2}
  \end{align}  
  Combining \eqref{ineq.usedlater3.1} and \eqref{ineq.usedlater3.2}, one finds that the right-hand side of \eqref{ineq.usedlater2} has
  \begin{align}
    &\ \xi_\Scal \sqrt{\gamma_{\nabla^2 c}} \sum_{j\in[m]} |[y_\Scal(x)]_j| +  \kappa_{\nabla^2 c} \sum_{j\in[m]} | [y_\Scal(x)]_j - [y(x)]_j| \nonumber \\
    \leq&\ \(\tfrac{3 \sqrt{m} \kappa_{\nabla^2 c} \kappa_1 (3\kappa_{\nabla f} \kappa_1 + \sqrt{\gamma_{\nabla f}}) + \sqrt{m} \sqrt{\gamma_{\nabla^2 c}} (5 \kappa_{\nabla f} \kappa_1 + \sqrt{\gamma_{\nabla f}})}{2 \sigma_{\min} \kappa_1}\) \xi_\Scal. \label{ineq.usedlater3}
  \end{align}
  As a result, combining \eqref{lemma.saa.Hess_f}, \eqref{ineq.usedlater1}, \eqref{ineq.usedlater2}, and \eqref{ineq.usedlater3} yields
  \begin{align}
    &\ |d_\Scal^T \nabla^2_{xx} L_\Scal(x,y_\Scal(x)) d_\Scal - d_\Scal^T \nabla^2_{xx} L (x,y(x)) d_\Scal| \nonumber \\
    \leq&\ \(\sqrt{\kappa_{\nabla^2 f}} + \tfrac{3 \sqrt{m} \kappa_{\nabla^2 c} \kappa_1 (3\kappa_{\nabla f} \kappa_1 + \sqrt{\gamma_{\nabla f}}) + \sqrt{m} \sqrt{\gamma_{\nabla^2 c}} (5 \kappa_{\nabla f} \kappa_1 + \sqrt{\gamma_{\nabla f}})}{2 \sigma_{\min} \kappa_1} \) \xi_\Scal \|d_\Scal\|^2. \label{ineq.theorem.v2}
  \end{align}
  Thus, now combining \eqref{ineq.dSLdS}, \eqref{ineq.theorem.v1}, \eqref{ineq.theorem.v2}, and \eqref{eq.def.kappa.3} one finds that
  \begin{align*}
    &\ |d_\Scal^T \nabla^2_{xx}L_\Scal(x,y_\Scal(x)) d_\Scal | \\
    \geq&\ \(1 - \tfrac{1}{3}\kappa_1 \xi_\Scal\) \beta \|d_\Scal\|^2 - 3 \(\kappa_{\nabla^2 f}+\tfrac{\sqrt{m} \kappa_{\nabla f} \kappa_{\nabla^2 c}}{\sigma_{\min}}\) \kappa_1 \xi_\Scal \|d_\Scal\|^2 \nonumber \\
    &\ - \(\sqrt{\kappa_{\nabla^2 f}} + \tfrac{3 \sqrt{m} \kappa_{\nabla^2 c} \kappa_1 (3\kappa_{\nabla f} \kappa_1 + \sqrt{\gamma_{\nabla f}}) + \sqrt{m} \sqrt{\gamma_{\nabla^2 c}} (5 \kappa_{\nabla f} \kappa_1 + \sqrt{\gamma_{\nabla f}})}{2 \sigma_{\min} \kappa_1} \) \xi_\Scal \|d_\Scal\|^2\\
    =&\ \(1 - \tfrac{1}{3}\kappa_1 \xi_\Scal\) \beta \|d_\Scal\|^2 - \kappa_3 \xi_\Scal \|d_\Scal\|^2 = \beta_\Scal \|d_\Scal\|^2.
  \end{align*}
  Finally, observe that $\beta_\Scal := (1 - \tfrac13 \kappa_1 \xi_\Scal) \beta - \kappa_3 \xi_\Scal \geq \thalf \beta$ follows from \eqref{ineq.theorem1.S}, from which one finds $(1 - \tfrac13 \kappa_1 \xi_\Scal) \geq 1 - \tfrac13 \kappa_1 \frac{1}{3\kappa_1} \geq \tfrac89 \beta$, so $(1 - \tfrac13 \kappa_1 \xi_\Scal) \beta - \thalf \beta \geq \tfrac{7}{18} \beta$, and consequently one finds that $\xi_\Scal \leq \tfrac{7\beta}{18 \kappa_3}$ yields $\beta_\Scal \geq \thalf \beta$.
\eproof

We now state our main result as Theorem~\ref{theo.complexity} below, which provides generic worst-case complexity bounds for Algorithm~\ref{alg.psm}.  The theorem states generic properties that a subproblem solver needs to have in order to leverage the strong-Morse properties of the subproblems such that, by solving a sequence of subproblems, the overall sample complexity can be less than that of solving the full-sample problem directly.  (As previously mentioned, the theorem focuses on the behavior of a second-order-type method for solving the subproblems and second-order optimality guarantees.)  A value of $p_1 \in [N]$ that satisfies the requirements of the theorem is not necessarily known in practice, but as shown in our numerical experiments one can obtain strong performance for a wide range of $p_1$ values rather than the full-sample problem directly, even if a chosen $p_1$ value does not satisfy the requirements of our theorem.

\begin{theorem}\label{theo.complexity}
  Suppose that Assumptions~\ref{ass.subproblems}, \ref{ass.boundness}, and \ref{ass.bounded.distribute} hold.  In addition,  suppose that $p_1 \in [N]$ is sufficiently large such that $\Scal_1$ satisfies $|\Scal_1| \geq p_1$ and~\eqref{ineq.theorem1.S} for $\Scal = \Scal_1$ $($where the tuple $($$\kappa_1,\kappa_2,\kappa_3$$)$ is defined by \eqref{eq.kappa} and \eqref{eq.def.kappa.3}$)$, and that, for all $k \in [K]$, the pair of subproblem solver tolerances $(\epsilon_k,\zeta_k) \in \R{}_{>0} \times \R{}_{>0}$ is set as a fraction of $(\alpha_{\Scal_k},\beta_{\Scal_k})$ $($defined by \eqref{eq.alpha-beta-s} with $\Scal = \Scal_k$$)$ such that any $(\epsilon_k,\zeta_k)$-stationary point with respect to~\eqref{prob.opt.S} for $\Scal = \Scal_k$, call it $x_k \in \R{n}$, satisfies
  \begin{align*}
    \|\nabla L_{\Scal_k}(x_k,y_{\Scal_k}(x_k))\| &\leq \epsilon_k \\
    \text{and}\ \ d^T \nabla_{xx}^2 L_{\Scal_k}(x_k,y_{\Scal_k}(x_k)) d &\geq -\zeta_k \|d\|_2^2\ \ \text{for all}\ \ d \in \Null(\nabla c_{\Scal_k}(x_k)^T).
  \end{align*}
  Finally, suppose that the sequence of tolerances $\{(\epsilon_k,\zeta_k)\}$ is monotonically nonincreasing and the subproblem solver employed by Algorithm~\ref{alg.psm} guarantees that:
  \benumerate
    \item[(a)] The number of iterations needed to reach an $(\epsilon_1,\zeta_1)$-stationary point for~\eqref{prob.opt.S} with $\Scal = \Scal_1$ is at most $\Ocal(\max\{\epsilon_1^{-2}, \zeta_1^{-3}\})$.
    \item[(b)] For any $k \geq 2$, the number of iterations needed from $x_{k-1}$ to reach an $(\epsilon_k,\zeta_k)$-stationary point for~\eqref{prob.opt.S} with $\Scal = \Scal_k$ is at most $\Ocal(\log(\tfrac{\epsilon_{k-1}}{\epsilon_k}))$.
  \eenumerate
  Then, the total number of individual objective and constraint gradient evaluations required by a run of Algorithm~\ref{alg.psm} with an initial sample set size $p_1 < N$ is at most
  \bequationNN
    \Ocal\(|\Scal_1| \max\{\epsilon_1^{-2}, \zeta_1^{-3}\} + \sum_{i=2}^K |\Scal_k| \log\(\tfrac{\epsilon_{k-1}}{\epsilon_k}\)\),
  \eequationNN
  whereas the total number of individual objective and constraint gradient evaluations required by a run of Algorithm~\ref{alg.psm} with an initial sample set size $p_1 = N$ is at most
  \bequationNN
    \Ocal(N\max\{\epsilon_1^{-2}, \zeta_1^{-3}\}).
  \eequationNN
\end{theorem}

The following corollary demonstrates a particular example of how Theorem~\ref{theo.complexity} can lead to a specific worst-case sample complexity result when a particular subproblem solver is employed.  In particular, the corollary leverages~\cite[Theorem 3.5]{GoyeEfteBoum2024} to prove that \cite[Algorithm 1]{GoyeEfteBoum2024} possesses a certain worst-case iteration complexity bound when employed as the subproblem solver in our Algorithm~\ref{alg.psm}.  This algorithm involves the minimization, for any $\Scal \subseteq [N]$, of Fletcher's augmented Lagrangian function with respect to our subproblem~\eqref{prob.opt.S}, namely, $F_\Scal : \R{n} \to \R{}$ defined by
\bequation\label{eq.fletcher.al.paper}
  F_\Scal(x) = f(x) + c_\Scal(x)^Ty_\Scal(x) + \rho_\Scal \|c_\Scal(x)\|^2,
\eequation
where $y_\Scal(x)$ is defined as in \eqref{eq.lsm} and $\rho_\Scal \in \R{}_{>0}$ is a penalty parameter. \iftechreport For a proof of the corollary, see Appendix~\ref{sec.specific}. \fi

\bcorollary\label{cor.example}
  Suppose that Assumptions \ref{ass.subproblems}, \ref{ass.boundness}, and \ref{ass.bounded.distribute} hold.  In addition, suppose that for any $\Scal \subseteq [N]$ with $|\Scal| \geq p_1$ and any $\rho_\Scal \in \R{}_{>0}$, Fletcher's Augmented Lagrangian function has a Hessian function $\nabla^2 F_\Scal : \R{n} \to \R{n \times n}$ that is Lipschitz in that there exists $M_\Scal \in \R{}_{>0}$ such that $\|\nabla^2 F_{\Scal}(x) - \nabla^2 F_{\Scal}(\xbar)\| \leq M_\Scal \|x - \xbar\|^2$ for all $(x,\xbar) \in \R{n} \times \R{n}$.  Furthermore, suppose that there exists $R \in \R{}_{>0}$ such that for any $\Scal \subseteq [N]$ with $|\Scal| \geq p_1$ the sublevel set $\Ccal_{\Scal,R} := \{x \in \R{n} : \|c_\Scal(x)\| \leq R\}$ contains $x_0$ and is compact.
  Then, the following statements hold true.
  \benumerate
    \item[(a)] Suppose that with $p_1 = N$ and a penalty parameter that is sufficiently large such that the requirements of \cite[Theorem 3.5]{GoyeEfteBoum2024} hold, \cite[Algorithm 1]{GoyeEfteBoum2024} is employed to solve~\eqref{prob.opt.S} for $\Scal = \Scal_1$.  Then, there exists a pair $(u_{[N],1},u_{[N],2}) \in \R{}_{>0} \times \R{}_{>0}$ such that, for any $(\epsilon,\zeta) \in (0,\tfrac{\sqrt{5}}{2}] \times (0,1]$, the number of individual objective and constraint gradient evaluations that are required until Algorithm~\ref{alg.psm} terminates with an $(\epsilon,\zeta)$ stationary point of~\eqref{prob.opt.N} is
    \bequationNN
      2N \left\lceil \max \left\{u_{[N],1} \epsilon^{-2}, u_{[N],2} \zeta^{-3} \right\} \right\rceil.
    \eequationNN
    \item[(b)]  Suppose that with $p_1 < N$ and a penalty parameter that is sufficiently large such that the requirements of \cite[Theorem 3.5]{GoyeEfteBoum2024} hold, \cite[Algorithm 1]{GoyeEfteBoum2024} is employed to solve~\eqref{prob.opt.S} for $\Scal = \Scal_1$, whereas gradient descent is employed to minimize Fletcher's augmented Lagrangian function for all subsequent $k \in [K]$ under the conditions of \iftechreport Lemma~\ref{lemma.fl.strongly.convex}. \fi  Then, there exists $(u_{\Scal_1,1},u_{\Scal_1,2}) \in \R{}_{>0} \times \R{}_{>0}$ such that, for any $(\epsilon,\zeta) \in (0,\tfrac{\sqrt{5}}{2}] \times (0,1]$ and with $(\tau_1,\kappa_1,\overline\omega)$ defined in \iftechreport Lemma~\ref{lemma.fl.strongly.convex}, \fi the number of individual objective and constraint gradient evaluations that are required until Algorithm~\ref{alg.psm} terminates with an $(\epsilon,\zeta)$ stationary point of~\eqref{prob.opt.N} is
    \begin{align*}
      &\ 2|\Scal_1|\left\lceil\max\left\{u_{\Scal_1,1} \epsilon_1^{-2}, u_{\Scal_1,2} \zeta_1^{-3} \right\}\right\rceil \\
      &\ + 2\(\tfrac{N-\theta|\Scal_1|}{\theta-1}\)\left\lceil\tfrac{1}{2}\log_2 \(45\overline{\omega}^2 \theta^2\(1+\tfrac{(\theta-1)N}{\theta}\)\)\right\rceil + 2N\left\lceil \log_{2}\tfrac{3\sqrt{5}\overline{\omega}\tau_1\sqrt{ \theta^2-\theta }}{\kappa_1\epsilon}\right\rceil.
    \end{align*}
  \eenumerate
\ecorollary

\section{Numerical Results}\label{sec.numerical}

The purpose of our numerical experiments is to demonstrate the performance of Algorithm~\ref{alg.psm} with $p_1 \ll N$ versus a one-shot approach that solves the SAA problem~\eqref{prob.opt.S} with $\Scal = [N]$ directly (i.e., $p_1 = N$).  We present the results of two experiments.  The first experiment involves an artificial two-dimensional problem, for which we present the results obtained using two (sub)problem solvers, namely, \cite[Algorithm~1]{GoyeEfteBoum2024} based on Fletcher's augmented Lagrangian function and \cite[Algorithm~2.2]{berahas2021sequential} based on a sequential quadratic optimization (SQP) approach.  We also use this two-dimensional problem to illustrate experimentally that for a problem with form~\eqref{prob.opt.S} that is strongly Morse when $\Scal = [N]$, a sampled problem with $\Scal \subset [N]$ is also strongly Morse when $|\Scal|$ is large enough relative to $N$.  Our second experiment involves training a neural network to predict the solution of an ordinary differential equation.  All software was written using Matlab R2024b.

\subsection{An Artificial Problem}\label{sec.artificial.prob}

Consider problem~\eqref{prob.expected} with $n = 2$, $m = 1$, and the objective function defined by $f(x) = x_1$, i.e., a deterministic objective.  As for the constraint function $\cbar : \R{n} \to \R{}$, let $\omega$ be a two-dimensional random vector with each component having a uniform distribution over $[-\pi,\pi]$.  Then, for any given $(a,\phi) \in \R{2}_{>0}$, let us define $\Cbar : \R{2} \times [-\pi,\pi]^2 \to \R{}$ according to
\bequationNN
  \Cbar(x,\omega) = x_1 - x_2^2 + a\sin(\phi x_1 + \omega_1) + a\cos(\phi x_2 + \omega_2).
\eequationNN
That is, for any $x \in \R{2}$, one finds that $\cbar(x) = \E[\Cbar(x,\omega)] = x_1 - x_2^2$.  On the other hand, for any positive integer $N$, the corresponding SAA problem has a constraint function that is the average of functions involving $\sin$ and $\cos$ terms.

\begin{figure}[ht]
  \centering
  \includegraphics[width=0.4\textwidth]{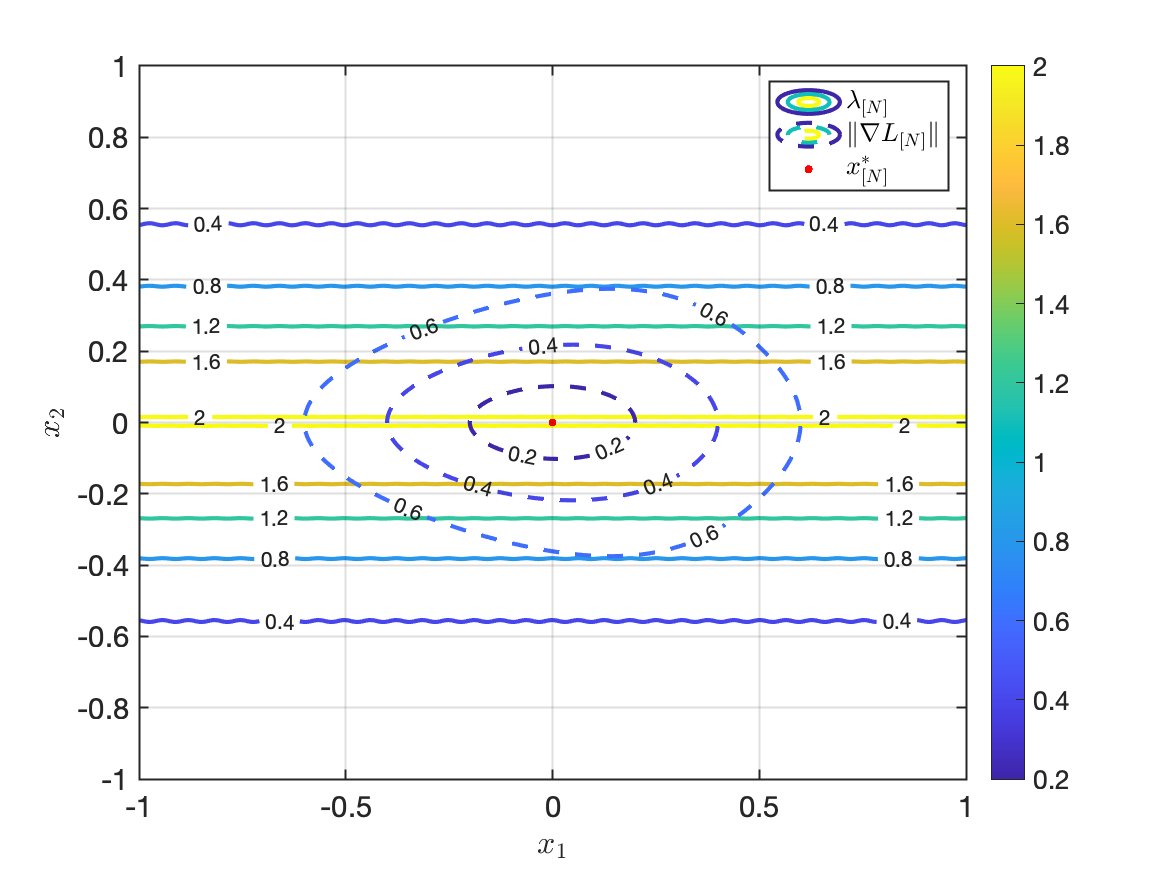}
  \includegraphics[width=0.4\textwidth]{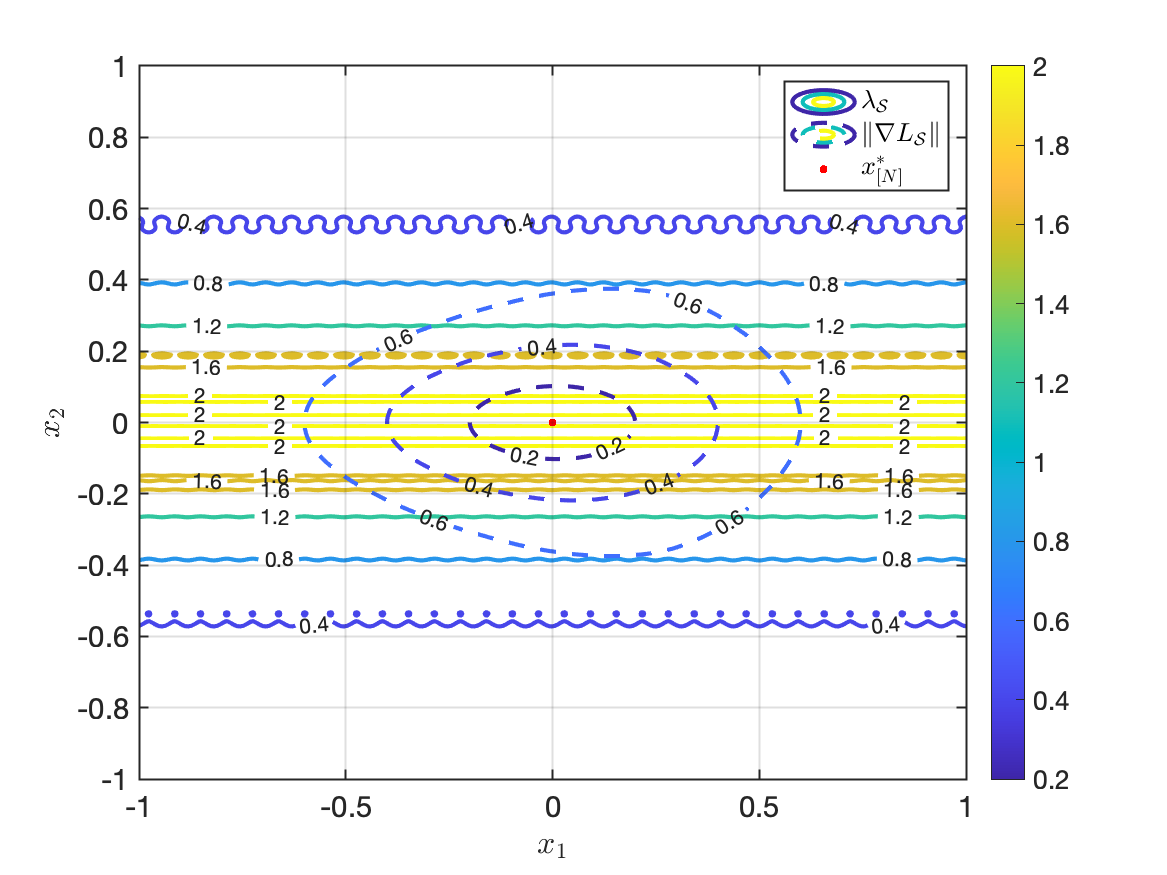}
  \caption{Contour plots for the problem defined in \S\ref{sec.artificial.prob} with different constraint sample sizes.  On the left and the right, $\lambda$ represents the reduced Hessian at $x$ (a real number), $\|\nabla L\|$ represents the norm of the gradient of the Lagrangian at $(x,y(x))$, and $x^*$ is the optimal solution of~\eqref{prob.opt.N}.  On the left, the contours correspond to function values over the full sample set with $N = 2048$. On the right, the contours correspond to values over a sample set with $|\Scal|=64$.}
  \label{fig.artificial.prob.struct}
\end{figure}

For our experiments in this section, we set $a = 10^{-4}$, $\phi = 100$, and $N = 2048$.  Our first aim is to demonstrate numerically that~\eqref{prob.opt.N} (i.e., with $\Scal=[N]$) is strongly Morse over $[-1,1]^2 \subset \R{2}$ (i.e., a box around the origin, the optimal solution of~\eqref{prob.expected}), and that~\eqref{prob.opt.S} is also strongly Morse over this region when $|S|$ is sufficiently large relative to $N$.  This is done through the contour plots in Figure~\ref{fig.artificial.prob.struct}.  These contour plots can be understood as follows.  First, for each sample size we computed the reduced Hessian of the Lagrangian at $x$---i.e., the Hessian of the Lagrangian over the null space of the constraint Jacobian---over the region, which in this setting is a real number $\lambda$ at each point.  Second, for each sample size we computed the norm of the gradient of the Lagrangian evaluated at $(x,y(x))$.  Third, for reference, we plot the point $x_{[N]}^* \approx 0$.  That each problem is strongly Morse can be seen by observing that~$\lambda$ is bounded away from zero when $\|\nabla L\|$ is small.  Specifically, the dashed contour lines for the gradient of the Lagrangian show the boundaries of the region in which $\|\nabla L_{\Scal}\| \leq 0.6$, where over this region one always finds that $\lambda_\Scal \geq 0.8$, which is to say that each problem is $(0.6,0.8)$-strongly Morse in the region $[-1,1]^2$.

Let us now present the results of Algorithm~\ref{alg.psm} with $p_1 = 64 \ll N = 2048$ vs.~$p_1 = N = 2048$.  We present the results with two subproblem solvers.  In each case, our performance criterion is the number of constraint gradient evaluations---i.e., evaluations of $\nabla c_i$ for some $i$---required before an iterate satisfying \eqref{eq.soc.S.approx} was obtained with tolerances $\epsilon = \zeta = 10^{-6}$.  An initial point $x_0$ was chosen at random from $[-1,1]^2$ and used for each run.  For the subproblem tolerances, we set for all $k \in [K]$
\bequation\label{eq.artificial.tolerence}
  \epsilon_k = \zeta_k = 10^{-6} \sqrt{\tfrac{N(N-|\Scal_k|)}{|\Scal_k|^2}+1}.
\eequation
This choice is motivated by \eqref{eq.key_value}, but with the addition of 1 inside the square root to ensure that the tolerance does not vanish at $k = K$.  Rather, the final tolerances for $k = K$ are $\epsilon = \zeta = 10^{-6}$, as stated previously.

Our first comparison considers \cite[Algorithm 1]{GoyeEfteBoum2024} as the subproblem solver, as is considered in Corollary~\ref{cor.example}. This is a second-order method based on minimizing Fletcher's augmented Lagrangian function; recall $F_{\Scal_k} : \R{n} \to \R{}$ defined by~\eqref{eq.fletcher.al.paper}. In short, the algorithm employs gradient descent with a backtracking line search on the augmented Lagrangian function until the norm of gradient of the augmented Lagrangian is below a threshold.  If the resulting point satisfies the second-order termination conditions, then the algorithm terminates; otherwise, a direction of negative curvature is obtained along which a backtracking line search is employed to determine the next iterate.  The algorithm continues until the second-order termination conditions with tolerances in~\eqref{eq.artificial.tolerence} are satisfied.  In our experiment, we set the penalty parameter $\rho_{\Scal_k}$ to be $10$ for all subproblems. With respect to the other parameters in \cite[Algorithm 1]{GoyeEfteBoum2024}, we chose $\alpha_{01} = \alpha_{02} = 1$, $c_1=c_2=10^{-4}$, and $\tau_1 = \tau_2 = \thalf$. 

\begin{figure}[ht]
  \centering
  \includegraphics[width=0.4\textwidth]{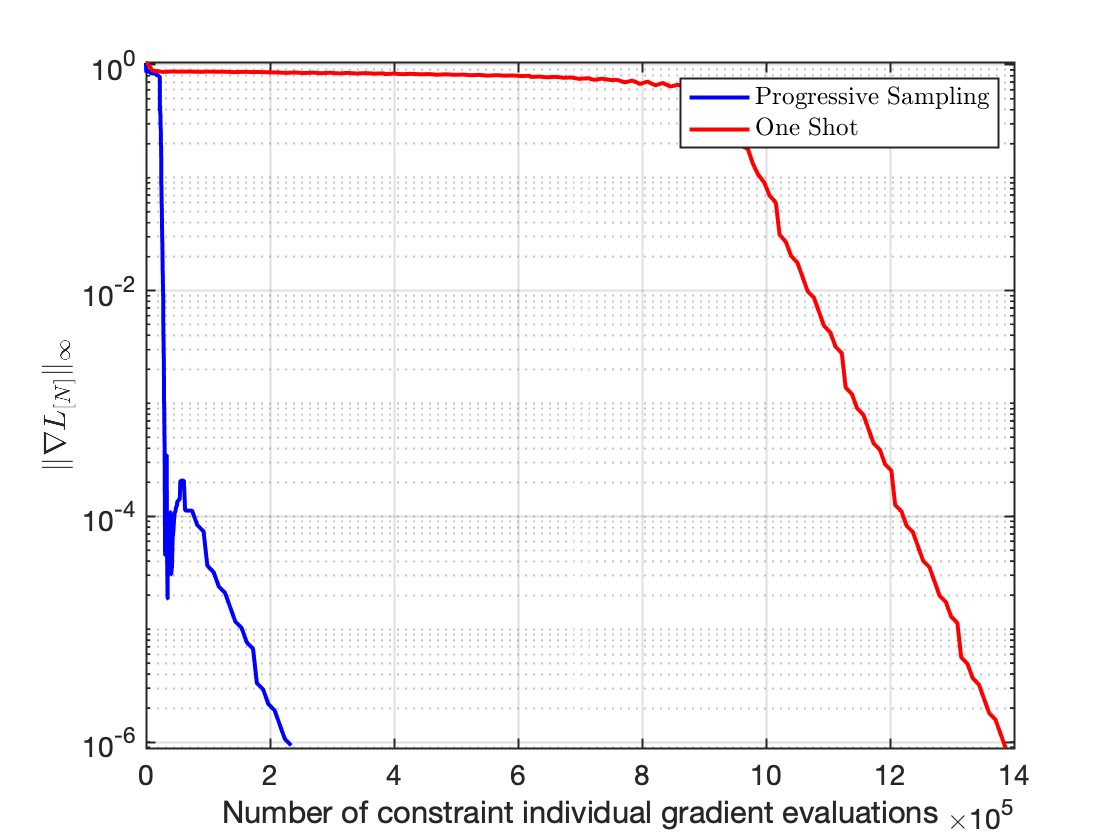}
  \includegraphics[width=0.4\textwidth]{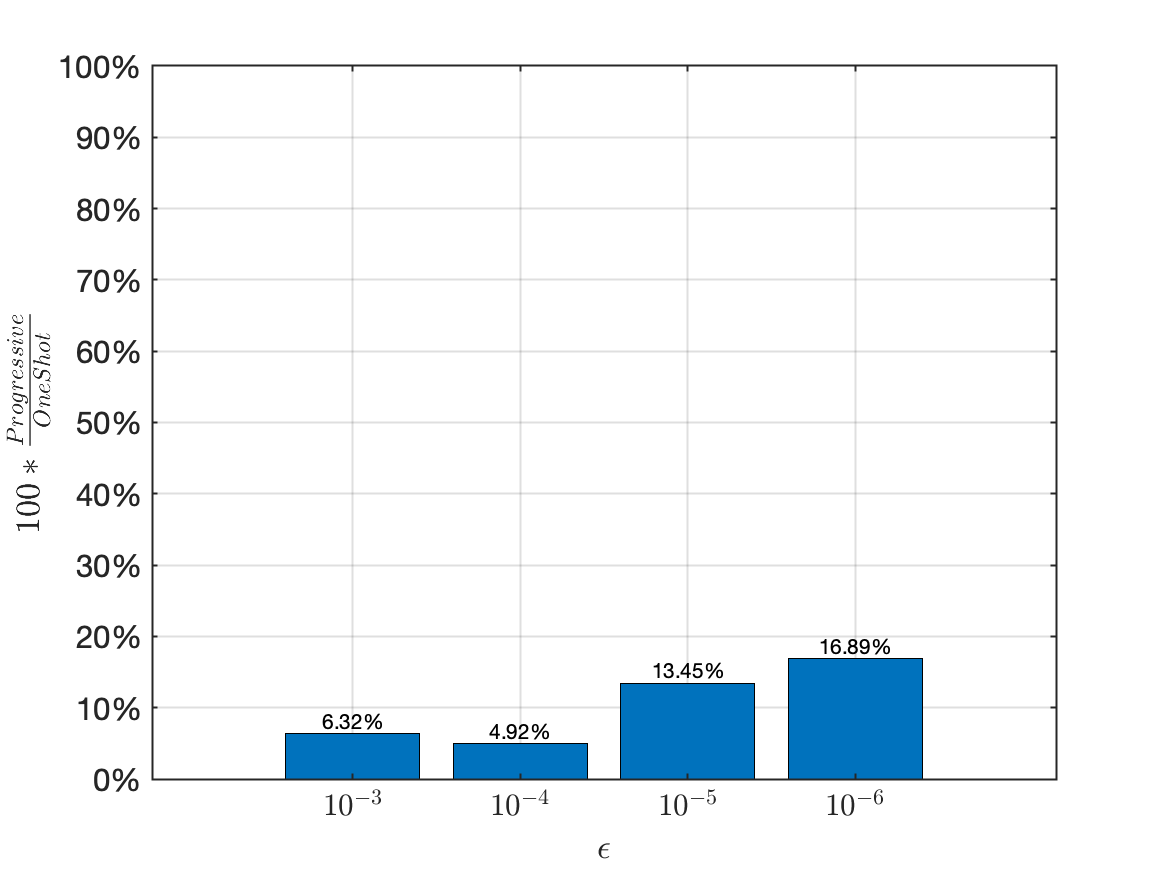}
\caption{Performance of Algorithm~\ref{alg.psm} with $p_1 \ll N$ (``Progressive Sampling'') versus $p_1 = N$ (``One Shot'') when solving the problem in \S\ref{sec.artificial.prob} using a subproblem solver based on minimizing Fletcher's augmented Lagrangian function.  On the left, the norm of the gradient of the Lagrangian with respect to the full-sample problem as a function of the number of constraint gradient evaluations.  On the right, the relative number of constraint gradient evaluations required by the two algorithm instances to obtain solutions satisfying the stated final tolerances.}
\label{fig.artificial.fal}
  \end{figure}

The results are shown in Figure \ref{fig.artificial.fal}.  On the left in the figure, we show the norm of the Lagrangian corresponding to the full-sample problem as a function of the number of constraint gradient evaluations as the algorithm proceeds.  One finds that despite the fact that $p_1 \ll N$, the progressive sampling approach yields a small norm of the Lagrangian after relatively few constraint gradient evaluations.  By contrast, by solving the full-sample problem directly, many more individual constraint gradient evaluations are required before the norm of the gradient of the Lagrangian falls below the desired tolerance.  On the right in the figure, we show the relative number of constraint gradient evaluations required when the final tolerance varies from $10^{-3}$ to $10^{-6}$.  One finds that the progressive sampling approach (with $p_1 \ll N$) requires only a small percentage of the constraint gradient evaluations that are required when solving the full-sample problem directly.

\begin{figure}[ht] 
  \centering
  \includegraphics[width=0.4\textwidth]{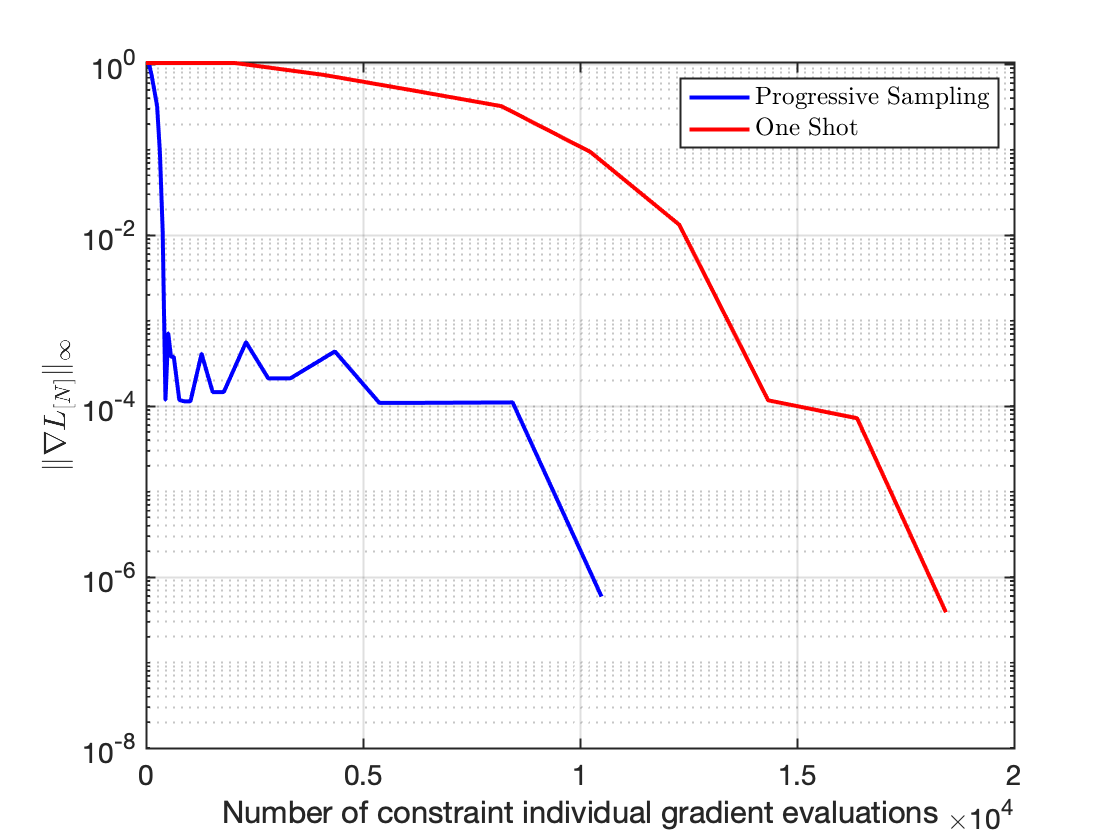}
  \includegraphics[width=0.4\textwidth]{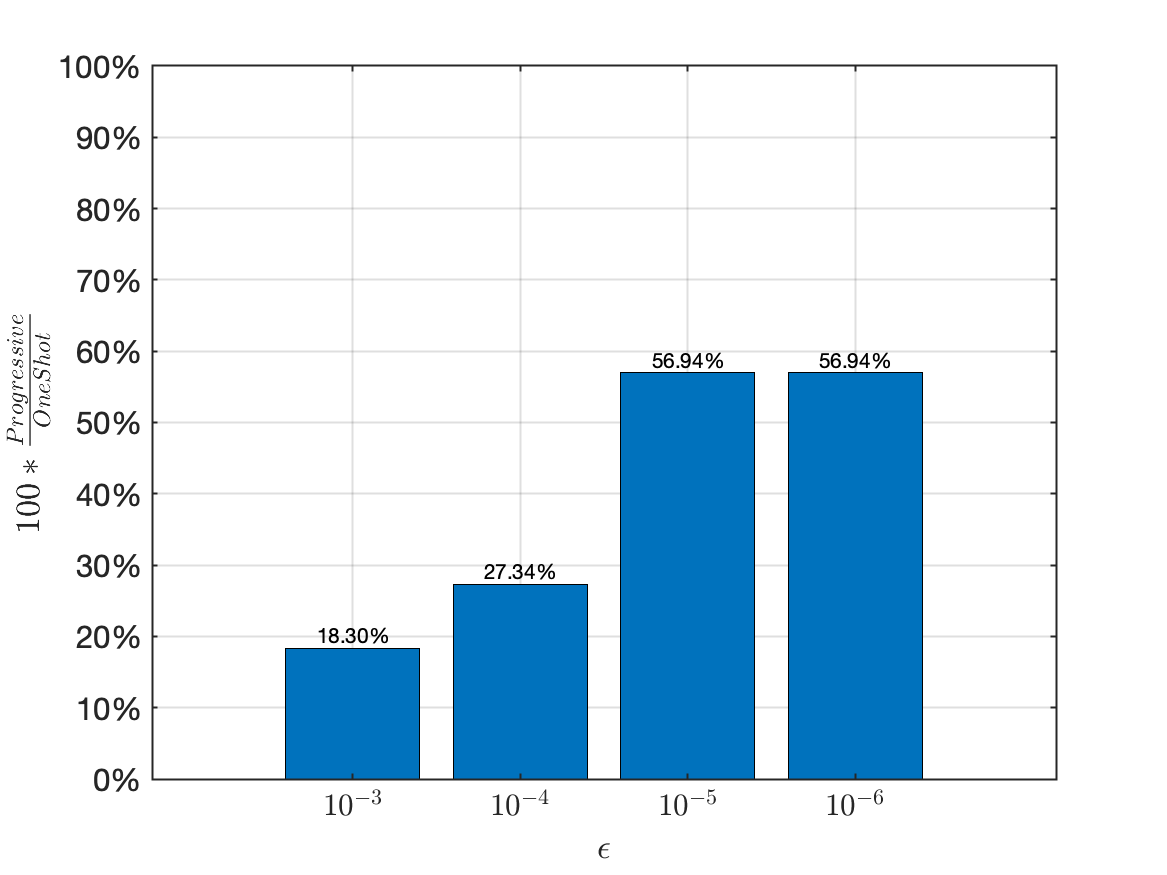}
  \caption{Performance of Algorithm~\ref{alg.psm} with $p_1 \ll N$ (``Progressive Sampling'') versus $p_1 = N$ (``One Shot'') when solving the problem defined in \S\ref{sec.artificial.prob} using a subproblem solver based on a first-order SQP approach.  On the left, the norm of the gradient of the Lagrangian with respect to the full-sample problem as a function of the number of constraint gradient evaluations.  On the right, the relative number of constraint gradient evaluations required by the two algorithm instances to obtain solutions satisfying the stated final tolerances.}
  \label{fig.artificial.sqp}
  \end{figure}

Our second comparison considers using \cite[Algorithm 2.2]{berahas2021sequential} as a subproblem solver, which is an SQP-based method.  This algorithm does not guarantee convergence to an approximate second-order stationary point, yet we still find that progressive sampling yields computational benefits.  For the parameters required by \cite[Algorithm 2.2]{berahas2021sequential}, we chose $(\alpha,\nu,\sigma,\eta,\tau)=(1,0.5,0.5,0.5,1)$. The results are shown in Figure~\ref{fig.artificial.sqp}.  As for the other subproblem solver, one finds that progressive sampling yields better results than the one-shot approach when considering individual constraint gradient evaluations as the performance measure.

\subsection{Training a Physics-Informed Neural Network}\label{sec.neuralODE}

Our second experiment involved training a neural network to predict the solution of an ordinary differential equation (ODE).  In particular, our experiment trained a physics-informed neural network to predict the movement of a damped harmonic (mass-spring) oscillator under the influence of a restoring force and friction~\cite{harmonic.oscillator}.  For the sake of simplicity, we trained the model for known ODE parameters and a single initial condition.  That being said, our algorithm is readily applicable to other such settings of training physics-informed neural networks as well.

The system that we considered is described by a linear, homogeneous, second-order ordinary differential equation with constant coefficients, namely,
\bequationNN
  m \tfrac{d^{2}u(t)}{dt^{2}} + \mu \tfrac{du(t)}{dt} + ku(t) = 0\ \ \text{over}\ \ t \in [0,10]\ \ \text{where}\ \ (m,\mu,k) = (1,0.1,1).
\eequationNN
We let the initial conditions be $u(0) = 1$ and $du(t)/dt = -1$.  To learn the solution of this ODE, we constructed a multilayer perceptron with one node in the input layer (representing time $t$), two hidden layers each with 128 nodes and tanh activation, and an output layer with a single node (representing the predicted height at time~$t$, i.e., $u(t)$).  The trainable parameters of the network are the edge weights as well as bias terms at each of the hidden layer nodes and the output node.  Overall, the network defined a function $\Ncal : \R{n} \times \R{}_{>0} \to \R{}$, where the first argument is the vector of trainable parameters and the second argument is the input (time).

The supervised training process that we considered employed a set of known input-output pairs to define the objective function; namely, for $N_f = 512$, let the pairs be denoted $\{(t_i,u_i)\}_{i\in[N_f]}$, where $u_i$ represents the true height of the spring at time~$t_i$ for all $i \in [N_f]$.  The pairs were determined by solving the ODE with a finite-difference method over 512 evenly spaced points over the time interval $[0,10]$.  We defined the objective function to have two terms, a data-fitting term and a term to minimize the residual of the ODE at the given input times.  Three constraints were included to aid in the training process.  The first two were the initial conditions while the latter required that the average ODE residual at a set of $N_c = 512$ times equals zero.  Overall, the training problem was
\begin{align}
  \min_{x \in \R{n}} &\ \tfrac{1}{N_f} \sum_{i\in[N_f]} \(u_i - \Ncal\(x, \tfrac{10i}{N_f}\)\)^2 + \nonumber \\
  &\ \tfrac{1}{N_c} \sum_{i\in[N_c]} \( \tfrac{\partial^2 \Ncal\(x, \tfrac{10i}{N_c}\)}{\partial t^2} + 0.1 \tfrac{\partial \Ncal\(x, \tfrac{10i}{N_c}\)}{\partial t} + \Ncal\(x, \tfrac{10i}{N_c}\) \)^2 \nonumber \\
  \st &\ \ \Ncal(x,0) = 1,\ \ \tfrac{\partial \Ncal(x,0)}{\partial t} = -1,\ \ \text{and} \nonumber \\
  &\ \ \tfrac{1}{N_c} \sum_{i\in[N_c]} \( \tfrac{\partial^2 \Ncal\(x, \tfrac{10i}{N_c}\)}{\partial t^2} + 0.1 \tfrac{\partial \Ncal\(x, \tfrac{10i}{N_c}\)}{\partial t} + \Ncal\(x, \tfrac{10i}{N_c}\) \) = 0. \label{eq.prob.pinn}
\end{align}

We ran Algorithm~\ref{alg.psm} numerous times with different $(p_{f,1},p_{c,1})$ values.  For the subproblem solver, we employed an implementation of \cite[Algorithm 2.2]{berahas2021sequential}.  This is a first-order method, not a second-order method, yet our results still show a computational benefit of progressive sampling.  The results can be found in Table~\ref{tab:data_accesses_progressive}, where we report the total number of data accesses in each run, which we define to be the total number of times that either an individual objective gradient or individual gradient of constraint~\eqref{eq.prob.pinn} was evaluated. The minimum across each row is written in italicized text and the minimum across each column is written in bold text; when the minimum in a row and column occurs in the same position, we put a box around the number.  The results show that, for any $p_{f,1}$ value, the best result was obtained with $p_{c,1} < 512$, and similarly that for any $p_{c,1}$ value, the best result was obtained with $p_{f,1} < 512$; in other words, progressive sampling was always beneficial compared to solving with the full sample initially.  Overall, the run with $(p_{f,1},p_{c,1}) = (16,32)$ required the fewest number of data accesses.  The required number was approximately only 35\% of the data accesses that were required compared to the run with $(p_{f,1},p_{c,1}) = (512,512)$, which solved the full-sample problem directly.

\begin{table}[ht]
  \tiny
  \centering
  \begin{tabular}{|c|cccccccccc|}
    \multicolumn{1}{c}{} & \multicolumn{10}{c}{$p_{f,1}$} \\
    \cline{2-11}
    \multicolumn{1}{c|}{$p_{c,1}$} & 1 & 2 & 4 & 8 & 16 & 32 & 64 & 128 & 256 & 512 \\
    \hline
    1 & 7.45e7 & 1.03e8 & 7.87e7 & 7.35e7 & 8.38e7 & \textit{1.88e7} & 5.83e7 & 5.96e7 & 5.05e7 & 1.95e8 \\
    2 & \textbf{1.58e7} & \textit{\textbf{9.46e6}} & 2.03e7 & 1.95e7 & 1.74e7 & 1.63e7 & 2.03e7 & 4.17e7 & 4.56e7 & 5.52e7 \\
    4 & 1.92e7 & 3.14e7 & 2.18e7 & 3.50e7 & \textit{1.32e7} & 1.69e7 & 1.61e7 & 2.03e7 & 2.44e7 & 4.40e7 \\
    8 & 2.80e7 & 1.94e7 & \textit{\textbf{8.11e6}} & 2.23e7 & 1.37e7 & \textbf{1.21e7} & 1.17e7 & \textbf{8.90e6} & 1.34e7 & 1.89e7 \\
    16 & 4.40e7 & 2.69e7 & 1.22e7 & \textit{\textbf{9.58e6}} & 1.07e7 & 1.58e7 & 1.46e7 & 1.23e7 & \textbf{1.19e7} & \textbf{1.55e7} \\
    32 & 7.19e7 & 4.78e7 & 1.85e7 & 1.52e7 & \textit{\textbf{8.09e6}} & 1.29e7 & 1.25e7 & 1.34e7 & 1.44e7 & 1.86e7 \\
    64 & 6.93e7 & 7.53e7 & 3.87e7 & 1.89e7 & 1.32e7 & 1.48e7 & 1.42e7 & \textit{1.20e7} & 1.22e7 & 1.57e7 \\
    128 & 8.99e7 & 8.28e7 & 6.06e7 & 3.86e7 & 2.52e7 & 1.54e7 & \textbf{1.07e7} & \textit{9.64e6} & 1.25e7 & 1.75e7 \\
    256 & 8.83e7 & 8.58e7 & 5.47e7 & 8.62e7 & 4.56e7 & 2.77e7 & 1.58e7 & \textit{1.25e7} & 1.43e7 & 2.09e7 \\
    512 & 9.21e7 & 9.29e7 & 5.93e7 & 1.20e8 & 7.35e7 & 4.89e7 & 2.37e7 & \textit{1.60e7} & 1.71e7 & 2.31e7 \\
    \hline
  \end{tabular}
  \caption{Number of data accesses for different $(p_{f,1},p_{c,1})$ pairs for progressive sampling method}
  \label{tab:data_accesses_progressive}
\end{table}

In order to compare the results obtained with Algorithm~\ref{alg.psm} to an alternative type of strategy that has been explored in the literature, we provide the results in Table~\ref{tab:data_accesses_geometric}.  In these runs, rather than solve each subproblem to a desired level of accuracy, the subproblem solver was run $k_{\max}$ iterations, at which time the sample size was increased by a prescribed factor and the merit parameter $\rho$ was increased by a prescribed factor in order to reset the subproblem solver for the subsequent solve.  This strategy is meant to mimic the behavior of adaptive sampling strategies that have been proposed in recent articles, where in some settings one can prove strong worst-case complexity bounds as long as the sample size increases at a geometric rate.  For these runs, we employed the initial sample sizes for the objective and the constraints as the pair $(16,32)$, which gave the best results for progressive sampling.  One finds in Table~\ref{tab:data_accesses_geometric} that the total number of data accesses for all of the tested parameter choices is larger than the number required by Algorithm~\ref{alg.psm} when progressive sampling is employed with $(p_{f,1},p_{c,1}) = (16,32)$.  This shows that, in these experiments, it was better to run the subproblem solver to a prescribed accuracy before increasing the sample size, rather than employ an adaptive-like approach that increases the sample sizes more gradually.  We also ran these experiments with larger $k_{\max}$ values, but in this case the subproblems are solved to high accuracy, meaning that the approach essentially reduces to a progressive sampling strategy like the one that we have proposed.

\begin{table}[ht]
  \tiny
  \centering
  \begin{tabular}{|c|ccc|ccc|ccc|}
    \cline{2-10}
    \multicolumn{1}{c|}{} & \multicolumn{9}{c|}{\textrm{$\rho$ increase factor}} \\
    \multicolumn{1}{c|}{} & \multicolumn{3}{c}{1} & \multicolumn{3}{c}{2} & \multicolumn{3}{c|}{4} \\
    \cline{2-10}
    \multicolumn{1}{c|}{} & \multicolumn{3}{c|}{\textrm{sample increase factor}} & \multicolumn{3}{c|}{\textrm{sample increase factor}} & \multicolumn{3}{c|}{\textrm{sample increase factor}} \\
    \multicolumn{1}{c|}{$k_{\max}$} & 1.1 & 2.0 & 4.0 & 1.1 & 2.0 & 4.0 & 1.1 & 2.0 & 4.0 \\
    \hline
    10 & 1.86e7 & 1.95e7 & 1.97e7 & 1.86e7 & 1.95e7 & 1.99e7 & 1.86e7 & 1.95e7 & 1.99e7 \\
    100 & 1.30e7 & 2.02e7 & 1.95e7 & 1.25e7 & 2.02e7 & 1.98e7 & 1.25e7 & 2.02e7 & 1.98e7 \\
    1000 & 8.87e6 & 8.69e6 & 1.34e7 & 8.91e6 & 8.53e6 & 1.29e7 & 8.91e6 & 8.53e6 & 1.29e7 \\
    \hline
  \end{tabular}
  \caption{Number of data accesses for geometric increase method}
  \label{tab:data_accesses_geometric}
\end{table}

\section{Conclusion}\label{sec.conclusion}

We have proposed, analyzed, and tested an algorithm for solving a sample average approximation \eqref{prob.opt.N} of any problem of the form \eqref{prob.expected}, where the objective and equality-constraint functions are defined by expectations.  The algorithm is based on a progressive sampling strategy, where in each iteration a deterministic algorithm can be employed to solve the resulting equality-constrained problem.  We have shown that the algorithm can achieve improved worst-case sample complexity compared to an approach that solves \eqref{prob.opt.N} directly.  Our numerical experiments have shown that our algorithm can offer computational savings in practice.

%% file: appendix.tex
\section{Subproblem Solver for Algorithm~\ref{alg.psm}}\label{sec.specific}

Our goal in this appendix is to prove Corollary~\ref{cor.example}, which is to say that the purpose is to prove the existence of an algorithm that yields the subproblem requirements of our generic sample-complexity bound in Theorem~\ref{theo.complexity}.  The algorithm that we consider is based on minimization of Fletcher's augmented Lagrangian function, and in particular the strategies proposed and analyzed in \cite{GoyeEfteBoum2024}.  Much of our efforts in this subsection focus on showing that, by employing the algorithm---specifically, \cite[Algorithm 1]{GoyeEfteBoum2024}---we can obtain our desired complexity bounds with respect to our termination conditions in \eqref{eq.soc.S.approx}. We emphasize that this is nontrivial since our termination conditions are different than those in~\cite{GoyeEfteBoum2024}. For simplicity throughout this appendix, we presume that the objective function sampling is initialized with the full sample so that progressive sampling is only conducted with the constraints. This allows us to focus on the particular challenges of constraint sampling. Our analysis is readily extended to the setting of objective function sampling as well.  The main approach of our proofs would be the same; all that objective function sampling would entail are expressions that are more cumbersome.

We begin with a lemma showing that, if the initial sample size is sufficiently large and the subproblem tolerances are set appropriately, then the iterates computed by each subproblem solve are approximate second-order stationary points.

\begin{lemma}\label{lem.start_for_next}
  Suppose that $($$\kappa_1,\kappa_2,\kappa_3$$)$ is defined by \eqref{eq.kappa} and \eqref{eq.def.kappa.3}, and with
  \bequation\label{eq.def.tau2}
    \tau_1 := \kappa_2\ \ \text{and}\ \ \tau_2 := \( \kappa_{\nabla^2 f} + \tfrac{\sqrt{m} \kappa_{\nabla f} \kappa_{\nabla^2 c}}{\sigma_{\min}}\)\kappa_1
  \eequation
  that the subproblem tolerances are set for all $k \in [K]$ $($recall \eqref{eq.key_value}$)$ as
  \bequation\label{eq.tolorence}
    \epsilon_k := \tau_1 \xi_{\Scal_k}\ \ \text{and}\ \ \zeta_k := \tau_2 \xi_{\Scal_k}.
  \eequation
  Then, as long as $p_1 \in [N]$ is sufficiently large such that
  \bequation\label{eq.theorem2.S}
    \sqrt{\tfrac{N(N-p_1)}{p_1^2}} \leq \min \left\{ \tfrac{1}{3\kappa_1}, \tfrac{\alpha}{4\kappa_2}, \tfrac{2\beta}{9\kappa_3} \right\},
  \eequation
  it follows for any $k \in [K]$ that if $x_{\Scal_k} \in \R{n}$ is $(\epsilon_k,\zeta_k)$-stationary with respect to~\eqref{prob.opt.S} for $\Scal = \Scal_k$, then with respect to~\eqref{prob.opt.S} with $\Scal=\Scal_{k+1}$ one has
  \bequation\label{ineq.init.next}
    \begin{aligned}
      \|\nabla L_{\Scal_{k+1}}(x_{\Scal_k}, y_{\Scal_{k+1}}(x_{\Scal_k}))\| &\leq 3\epsilon_k \leq \alpha_{\Scal_{k+1}}\ \ \text{and} \\
      d^T \nabla^2_{xx} L_{\Scal_{k+1}}(x_{\Scal_k}, y_{\Scal_{k+1}}(x_{\Scal_k})) d &\geq \beta_{\Scal_{k+1}} \|d\|^2\ \text{for all}\ d \in \Null(\nabla c_{\Scal_{k+1}}(x_{\Scal_k})^T),
    \end{aligned}
  \eequation
  where $\alpha_{\Scal_{k+1}} \geq \thalf \alpha$ and $\beta_{\Scal_{k+1}} \geq \thalf \beta$ are defined according to \eqref{eq.alpha-beta-s}.
\end{lemma}
\bproof
  Consider arbitrary $k \in [K]$.  Our first aim is to bound the norm of the difference between the gradients of $L_{\Scal_k}$ and $L_{\Scal_{k+1}}$ with respect to the point~$x_{\Scal_k}$.  Toward this end, by \eqref{eq.lsm}, define $y_{\Scal_k} := y_{\Scal_k}(x_{\Scal_k})$ and $y_{\Scal_{k+1}} := y_{\Scal_{k+1}}(x_{\Scal_k})$.  Then, by similar arguments as led to \eqref{ineq.L.grad.x.diff.1} and the triangle inequality, one finds
  \begin{align}
    &\ \|\nabla_x L_{\Scal_{k+1}}(x_{\Scal_k},y_{\Scal_{k+1}}) - \nabla_x L_{\Scal_k}(x_{\Scal_k},y_{\Scal_k})\| \nonumber \\
    \leq&\ \|\Rcal(\nabla c_{\Scal_{k+1}}(x_{\Scal_k})) - \Rcal(\nabla c_{\Scal_k}(x_{\Scal_k}))\| \kappa_{\nabla f} \nonumber \\
    =&\ \|(\Rcal(\nabla c_{\Scal_{k+1}}(x_{\Scal_k})) - \Rcal(\nabla c(x_{\Scal_k}))) + (\Rcal(\nabla c (x_{\Scal_k})) - \Rcal(\nabla c_{\Scal_k}(x_{\Scal_k})))\| \kappa_{\nabla f} \nonumber \\
    \leq&\ (\|\Rcal(\nabla c_{\Scal_{k+1}}(x_{\Scal_k})) - \Rcal(\nabla c(x_{\Scal_k}))\| + \|\Rcal(\nabla c(x_{\Scal_k})) - \Rcal(\nabla c_{\Scal_k}(x_{\Scal_k}))\|) \kappa_{\nabla f}. \label{ineq.theorem2.Lx.0}
  \end{align}
  Next, let us bound the two terms in parentheses on the right-hand side of \eqref{ineq.theorem2.Lx.0}.  Note that by Assumption~\ref{ass.subproblems} the matrices $\nabla c (x_{\Scal_k})$, $\nabla c_{\Scal_k}(x_{\Scal_k})$, and $\nabla c_{\Scal_{k+1}} (x_{\Scal_k})$ all have full-column rank, meaning that they have the same rank, so by \cite[Theorems~2.3--2.4]{Stew1977} and the same arguments as led to \eqref{ineq.L.grad.x.diff}, one finds
  \bsubequations
    \begin{align}
      \|\Rcal(\nabla c_{\Scal_k}(x_{\Scal_k})) - \Rcal(\nabla c(x_{\Scal_k}))\| &\leq \kappa_1 \xi_{\Scal_k} \label{ineq.r_s_k} \\ \text{and}\ \ 
      \|\Rcal(\nabla c_{\Scal_{k+1}}(x_{\Scal_k})) - \Rcal(\nabla c(x_{\Scal_k}))\| &\leq \kappa_1 \xi_{\Scal_{k+1}} \leq \kappa_1 \xi_{\Scal_k}, \label{ineq.r_s_k+1}
    \end{align}
  \esubequations
  where the last inequality follows from \eqref{eq.key_value} and the fact that the algorithm guarantees that $|\Scal_{k+1}| > |\Scal_k|$.  Combining \eqref{ineq.r_s_k} and \eqref{ineq.r_s_k+1} with \eqref{ineq.theorem2.Lx.0}, one obtains
  \bequation\label{ineq.theorem2.Lx}
    \|\nabla_x L_{\Scal_{k+1}}(x_{\Scal_k},y_{\Scal_{k+1}}) - \nabla_x L_{\Scal_k}(x_{\Scal_k},y_{\Scal_k})\| \leq 2 \kappa_{\nabla f} \kappa_1 \xi_{\Scal_k}.
  \eequation
  On the other hand, from the triangle inequality, \eqref{lemma.saa.c}, and $\xi_{\Scal_{k+1}} \leq \xi_{\Scal_k}$, one has
  \begin{align}
    &\ \|\nabla_yL_{\Scal_{k+1}}(x_{\Scal_k},y_{\Scal_{k+1}})-\nabla_yL_{\Scal_k}(x_{\Scal_k},y_{\Scal_k})\| \nonumber \\
    =&\ \|c_{\Scal_{k+1}}(x_{\Scal_k})-c_{\Scal_k}(x_{\Scal_k})\| \nonumber \\
    \leq&\ \|c_{\Scal_{k+1}}(x_{\Scal_k})-c(x_{\Scal_k})\|+\|c(x_{\Scal_k})-c_{\Scal_k}(x_{\Scal_k})\| \nonumber \\
    \leq&\ \sqrt{\gamma_c} \xi_{\Scal_{k+1}} + \sqrt{\gamma_c} \xi_{\Scal_k} \leq 2\sqrt{\gamma_c} \xi_{\Scal_k}. \label{ineq.theorem2.Ly}
  \end{align}
  Combining \eqref{ineq.theorem2.Lx} and \eqref{ineq.theorem2.Ly}, one finds that
  \bequation\label{ineq.LklK1}
    \|\nabla L_{\Scal_{k+1}}(x_{\Scal_k},y_{\Scal_{k+1}}) - \nabla L_{\Scal_k}(x_{\Scal_k},y_{\Scal_k})\|^2 \leq 4 (\gamma_c + \kappa_{\nabla f}^2 \kappa_1^2) \xi_{\Scal_k}^2 = 4 \kappa_2^2 \xi_{\Scal_k}^2,
  \eequation
  which further gives under the conditions of the lemma that
  \begin{align}
    &\ \|\nabla L_{\Scal_{k+1}}(x_{\Scal_k},y_{\Scal_{k+1}})\| \nonumber \\
    \leq&\ \|\nabla L_{\Scal_{k+1}}(x_{\Scal_k},y_{\Scal_{k+1}})-\nabla L_{\Scal_k}(x_{\Scal_k},y_{\Scal_k})\| +\|\nabla L_{\Scal_k}(x_{\Scal_k},y_{\Scal_k})\|\\
    \leq&\ 2 \kappa_2 \xi_{\Scal_k} + \epsilon_k = 3\epsilon_k = 3 \tau_1 \xi_{\Scal_k} \leq \tfrac{3}{4} \alpha. \label{ineq.3/4alpha}
  \end{align}
  At the same time, by \eqref{eq.alpha-beta-s}, $\xi_{\Scal_{k+1}} \leq \xi_{\Scal_k}$, and $\kappa_2 \xi_{\Scal_k} \leq \tfrac{1}{4}\alpha$, one has
  \bequation\label{ineq.alpha.k.k+1}
    \alpha_{\Scal_{k+1}} = \alpha - \kappa_2 \xi_{\Scal_{k+1}} \geq \alpha - \kappa_2 \xi_{\Scal_k} \geq \alpha - \tfrac{1}{4}\alpha = \tfrac{3}{4}\alpha.
  \eequation
  Combining \eqref{ineq.3/4alpha} and \eqref{ineq.alpha.k.k+1}, one obtains the first desired conclusion in \eqref{ineq.init.next} that
  \bequation\label{ineq.S.k+1.grad.cond}
    \|\nabla L_{\Scal_{k+1}}(x_{\Scal_k}, y_{\Scal_{k+1}})\| \leq 3\epsilon_k \leq \alpha_{\Scal_{k+1}}.
  \eequation

  Our next goal is to prove the second inequality in \eqref{ineq.init.next}.  Toward this end, first note that Lemma~\ref{lem: emprical morse} applies here since the right-hand side of \eqref{eq.theorem2.S} is less than or equal to the right-hand side of~\eqref{ineq.theorem1.S}.  Hence, by Lemma~\ref{lem: emprical morse}, it follows that~\eqref{prob.opt.S} with $\Scal = \Scal_{k+1}$ is $(\alpha_{\Scal_{k+1}},\beta_{\Scal_{k+1}})$-strongly Morse.  Combined with the prior conclusion that $x_{\Scal_k}$ satisfies~\eqref{ineq.S.k+1.grad.cond}, one has for all $d_{\Scal_{k+1}}\in\Null(\nabla c_{\Scal_{k+1}}(x_{\Scal_k})^T)$ that
  \bequationNN
    |d_{\Scal_{k+1}}^T \nabla^2_{xx} L_{\Scal_{k+1}}(x_{\Scal_k},y_{\Scal_{k+1}}) d_{\Scal_{k+1}}| \geq \beta_{\Scal_{k+1}} \|d_{\Scal_{k+1}}\|^2.
  \eequationNN
  To prove the second inequality in \eqref{ineq.init.next}, it is necessary to show that, in fact, this inequality always holds without the absolute value on the left-hand side.  Toward showing this, let us define $\dbar_{\Scal_{k+1}} := \Ncal(\nabla c_{\Scal_k}(x_{\Scal_k})^T) d_{\Scal_{k+1}}$.  Then, one finds
  \begin{align}
    &\ d_{\Scal_{k+1}}^T \nabla^2_{xx} L_{\Scal_{k+1}}(x_{\Scal_k},y_{\Scal_{k+1}}) d_{\Scal_{k+1}} \nonumber \\
    =&\ \dbar_{\Scal_{k+1}}^T \nabla^2_{xx} L_{\Scal_k}(x_{\Scal_k},y_{\Scal_k}) \dbar_{\Scal_{k+1}} \nonumber \\
    &\ + (d_{\Scal_{k+1}}^T \nabla^2_{xx} L_{\Scal_{k+1}}(x_{\Scal_k},y_{\Scal_{k+1}}) d_{\Scal_{k+1}} - \dbar_{\Scal_{k+1}}^T \nabla^2_{xx} L_{\Scal_k}(x_{\Scal_k},y_{\Scal_k}) \dbar_{\Scal_{k+1}}) \nonumber \\
    \geq&\ \dbar_{\Scal_{k+1}}^T \nabla^2_{xx} L_{\Scal_k}(x_{\Scal_k},y_{\Scal_k}) \dbar_{\Scal_{k+1}} \nonumber \\
    &\ -|d_{\Scal_{k+1}}^T \nabla^2_{xx} L_{\Scal_{k+1}}(x_{\Scal_k},y_{\Scal_{k+1}}) d_{\Scal_{k+1}} - \dbar_{\Scal_{k+1}}^T \nabla^2_{xx} L_{\Scal_k}(x_{\Scal_k},y_{\Scal_k}) \dbar_{\Scal_{k+1}}|. \label{ineq.dld}
  \end{align}
  Since $\Ncal(\nabla c_{\Scal_k}(x_{\Scal_k})^T)$ is a projection matrix, $\|\Ncal(\nabla c_{\Scal_k}(x_{\Scal_k})^T)\| \leq 1$ and
  \bequationNN
    \|\dbar_{\Scal_{k+1}}\|^2 = \|\Ncal(\nabla c_{\Scal_k}(x_{\Scal_k})^T)  d_{\Scal_{k+1}}\|^2 \leq \|\Ncal(\nabla c_{\Scal_k}(x_{\Scal_k})^T)\|^2 \|d_{\Scal_{k+1}}\|^2 \leq \|d_{\Scal_{k+1}}\|^2.
  \eequationNN
  Along with the fact that $x_{\Scal_k}$ is $(\epsilon_k,\zeta_k)$-stationary, this means that the first-term on the right-hand side of \eqref{ineq.dld} satisfies the inequalities
  \bequation\label{ineq.varepson.cond}
    \dbar_{\Scal_{k+1}}^T \nabla^2_{xx} L_{\Scal_k}(x_{\Scal_k},y_{\Scal_k}) \dbar_{\Scal_{k+1}}
\geq -\zeta_k\|\dbar_{\Scal_{k+1}}\|^2 \geq -\zeta_k \|d_{\Scal_{k+1}}\|^2.
  \eequation
  On the other hand, with respect to the second term on the right-hand side of \eqref{ineq.dld}, observe that with $y(x_{\Scal_k})$ defined by~\eqref{eq.lsm} the triangle inequality yields
  \begin{align}
    &\ |d_{\Scal_{k+1}}^T \nabla^2_{xx} L_{\Scal_{k+1}}(x_{\Scal_k},y_{\Scal_{k+1}}) d_{\Scal_{k+1}} - \dbar_{\Scal_{k+1}}^T \nabla^2_{xx} L_{\Scal_k}(x_{\Scal_k},y_{\Scal_k}) \dbar_{\Scal_{k+1}}| \nonumber \\
    \leq&\ |d_{\Scal_{k+1}}^T \nabla^2_{xx} L_{\Scal_{k+1}}(x_{\Scal_k},y_{\Scal_{k+1}}) d_{\Scal_{k+1}} - d_{\Scal_{k+1}}^T \nabla^2_{xx} L(x_{\Scal_k},y(x_{\Scal_k})) d_{\Scal_{k+1}} | \nonumber \\
    &\ + |d_{\Scal_{k+1}}^T \nabla^2_{xx} L(x_{\Scal_k},y(x_{\Scal_k})) d_{\Scal_{k+1}} - \dbar_{\Scal_{k+1}}^T \nabla^2_{xx} L(x_{\Scal_k},y(x_{\Scal_k})) \dbar_{\Scal_{k+1}}| \nonumber \\
    &\ + |\dbar_{\Scal_{k+1}}^T \nabla^2_{xx} L(x_{\Scal_k},y(x_{\Scal_k})) \dbar_{\Scal_{k+1}} - \dbar_{\Scal_{k+1}}^T \nabla^2_{xx} L_{\Scal_k}(x_{\Scal_k},y_{\Scal_k}) \dbar_{\Scal_{k+1}}|. \label{ineq.expand.3.24.1}
  \end{align}
  Next, one can observe that the first and third terms on the right-hand side of~\eqref{ineq.expand.3.24.1} can be bounded as in~\eqref{ineq.theorem.v2} with respect to $\Scal = \Scal_k$ and $\Scal = \Scal_{k+1}$, respectively.  To see this, note that a bound of the form in~\eqref{ineq.usedlater1} holds since \eqref{ineq.usedlater1} followed using general matrix-norm inequalities.  Furthermore, a bound of the form in~\eqref{ineq.usedlater2} that relies on~\eqref{ineq.usedlater1} and \eqref{lemma.saa.Hess_c} holds, and bounds of the form in~\eqref{ineq.usedlater3.1}--\eqref{ineq.usedlater3.2} hold since these rely on $\xi_{\Scal_k} \leq \tfrac{1}{3\kappa_1}$ and $\xi_{\Scal_{k+1}} \leq \tfrac{1}{3\kappa_1}$, which hold here.  Thus, since \eqref{ineq.theorem.v2} follows from a combination of \eqref{ineq.usedlater1}, \eqref{ineq.usedlater2}, \eqref{ineq.usedlater3.1}, and \eqref{ineq.usedlater3.2}, one obtains that with
  \bequationNN
    \lambda := \sqrt{m} \tfrac{\kappa_{\nabla f}}{2\sigma_{\min}}\(5\sqrt{\gamma_{\nabla^2 c}}    + 9\kappa_1\kappa_{\nabla^2 c} \)
  \eequationNN
  the first and third terms on the right-hand side of~\eqref{ineq.expand.3.24.1} respectively satisfy
  \begin{align}
    &\ |d_{\Scal_{k+1}}^T \nabla^2_{xx} L_{\Scal_{k+1}}(x_{\Scal_k},z_{\Scal_{k+1}}) d_{\Scal_{k+1}} - d_{\Scal_{k+1}}^T \nabla^2_{xx} L(x_{\Scal_k},y(x_{\Scal_k})) d_{\Scal_{k+1}}| \nonumber \\
    \leq&\ \lambda \xi_{\Scal_{k+1}} \|d_{\Scal_{k+1}}\|^2 \leq \lambda \xi_{\Scal_k} \|d_{\Scal_{k+1}}\|^2 \label{ineq.3.24.1} \\
    \text{and}\ \ 
    &\ |\dbar_{\Scal_{k+1}}^T \nabla^2_{xx} L(x_{\Scal_k},y(x_{\Scal_k})) \dbar_{\Scal_{k+1}} - \dbar_{\Scal_{k+1}}^T \nabla^2_{xx} L_{\Scal_k}(x_{\Scal_k},y_{\Scal_k}) \dbar_{\Scal_{k+1}}| \nonumber \\
    \leq&\ \lambda \xi_{\Scal_k} \|\dbar_{\Scal_{k+1}}\|^2 \leq \lambda \xi_{\Scal_k} \|d_{\Scal_{k+1}}\|^2, \label{ineq.3.24.2}
  \end{align}
  As for the second term on the right of \eqref{ineq.expand.3.24.1}, it follows from the submultiplicity of the matrix 2-norm, the triangle inequality, $\|\dbar_{\Scal_{k+1}}\|\le\|d_{\Scal_{k+1}}\|$ and \eqref{lemma.bound.L_Hess} that
  \begin{align}
    &\ |d_{\Scal_{k+1}}^T \nabla^2_{xx} L(x_{\Scal_k},y(x_{\Scal_k})) d_{\Scal_{k+1}} - \dbar_{\Scal_{k+1}}^T \nabla^2_{xx} L(x_{\Scal_k},y(x_{\Scal_k})) \dbar_{\Scal_{k+1}}| \nonumber \\
    =&\ |(d_{\Scal_{k+1}} - \dbar_{\Scal_{k+1}})^T \nabla^2_{xx} L(x_{\Scal_k},y(x_{\Scal_k})) (d_{\Scal_{k+1}} + \dbar_{\Scal_{k+1}})| \nonumber \\
    \leq&\ \|\nabla^2_{xx} L(x_{\Scal_k},y(x_{\Scal_k}))\| \|d_{\Scal_{k+1}} - \dbar_{\Scal_{k+1}}\|\|d_{\Scal_{k+1}} + \dbar_{\Scal_{k+1}}\| \nonumber \\
    \leq&\ 2\(\kappa_{\nabla^2 f} + \tfrac{\sqrt{m} \kappa_{\nabla f} \kappa_{\nabla^2 c}}{\sigma_{\min}}\)\|d_{\Scal_{k+1}}\| \|d_{\Scal_{k+1}} - \dbar_{\Scal_{k+1}}\| \nonumber\\
    =&\ 2\tfrac{\tau_2}{\kappa_1}\|d_{\Scal_{k+1}}\| \|d_{\Scal_{k+1}} - \dbar_{\Scal_{k+1}}\|. \label{ineq.theorem2.iv2}
  \end{align}
  Now with respect to $\|d_{\Scal_{k+1}}- \dbar_{\Scal_{k+1}}\|$ one finds from the triangle inequality, submultiplicity of the matrix 2-norm, Lemma~\ref{lemma.proj.apply}, and \eqref{ineq.r_s_k} that
  \begin{align}\nonumber
    &\ \|d_{\Scal_{k+1}} - \dbar_{\Scal_{k+1}}\| \\
    =&\ \|\Rcal(\nabla c_{\Scal_k}(x_{\Scal_k})) d_{\Scal_{k+1}}\|\nonumber \\
    \leq&\ \|\Rcal(\nabla c(x_{\Scal_k})) d_{\Scal_{k+1}}\| + \|(\Rcal(\nabla c(x_{\Scal_k})) - \Rcal(\nabla c_{\Scal_k}(x_{\Scal_k})))\| \|d_{\Scal_{k+1}}\|\nonumber \\
    \leq&\ \kappa_1 \xi_{\Scal_{k+1}} \|d_{\Scal_{k+1}}\| + \kappa_1 \xi_{\Scal_k} \|d_{\Scal_{k+1}}\| \leq 2 \kappa_1 \xi_{\Scal_k} \|d_{\Scal_{k+1}}\|.\label{ineq.theorem2.d.dbar}
  \end{align}
  Combining \eqref{ineq.dld}--\eqref{ineq.theorem2.d.dbar} with \eqref{eq.def.kappa.3}, \eqref{eq.def.tau2}, and $\zeta_k = \tau_2 \xi_{\Scal_k}$, one finds
  \begin{align*}
    &\ d_{\Scal_{k+1}}^T \nabla^2_{xx} L_{\Scal_{k+1}}(x_{\Scal_k}, y_{\Scal_{k+1}}) d_{\Scal_{k+1}} \nonumber \\
    \geq&\ -\tau_2 \xi_{\Scal_k} \|d_{\Scal_{k+1}}\|^2 - 2 \lambda \xi_{\Scal_k} \|d_{\Scal_{k+1}}\|^2 - 4\tau_2 \xi_{\Scal_k} \|d_{\Scal_{k+1}}\|^2 \\
    =&\ \(-5\( \kappa_{\nabla^2 f} + \tfrac{\sqrt{m} \kappa_{\nabla f} \kappa_{\nabla^2 c}}{\sigma_{\min}}\)\kappa_1-2 \tfrac{\sqrt{m}\kappa_{\nabla f}}{2\sigma_{\min}}\(5\sqrt{\gamma_{\nabla^2 c}}    + 9\kappa_1 \kappa_{\nabla^2 c}\)\)\xi_{\Scal_k}\|d_{\Scal_{k+1}}\|^2\\
    \geq&\ - 2\(3\kappa_{\nabla^2 f}\kappa_1 + \(\tfrac{5\sqrt{m}\kappa_{\nabla f}}{2\sigma_{\min}}\) (3\kappa_1\kappa_{\nabla^2 c} + \sqrt{\gamma_{\nabla^2 c}})\) \xi_{\Scal_k} \|d_{\Scal_{k+1}}\|^2\\
    =&\ - 2 \kappa_3 \xi_{\Scal_k} \|d_{\Scal_{k+1}}\|^2.
  \end{align*}
  Since \eqref{eq.theorem2.S} requires $\xi_{\Scal_k} \leq \tfrac{2\beta}{9\kappa_3}$ it follows that
  \bequation\label{ineq.var.k+1}
    d_{\Scal_{k+1}}^T \nabla^2_{xx} L_{\Scal_{k+1}}(x_{\Scal_k},y_{\Scal_{k+1}}) d_{\Scal_{k+1}} \geq -2 \kappa_3 \xi_{\Scal_k} \|d_{\Scal_{k+1}}\|^2 \geq -\tfrac{4}{9} \beta \|d_{\Scal_{k+1}}\|^2.
  \eequation
  On the other hand, recall from \eqref{eq.alpha-beta-s} that $\beta_{\Scal_{k+1}} = \(1-\tfrac{1}{3} \kappa_1 \xi_{\Scal_{k+1}}\)\beta -\kappa_3 \xi_{\Scal_{k+1}}$, which along with $\xi_{\Scal_{k+1}} \leq \xi_{\Scal_k}$ and \eqref{eq.theorem2.S} (specifically, $\xi_{\Scal_k} \leq \min\{\tfrac{1}{3\kappa_1}, \tfrac{2\beta}{9\kappa_3}\})$ implies
  \bequationNN
    \beta_{\Scal_{k+1}} \geq \(1-\tfrac{1}{3} \kappa_1 \xi_{\Scal_k}\) \beta - \kappa_3 \xi_{\Scal_k} \geq \beta - \tfrac19 \beta - \tfrac29 \beta = \tfrac{2}{3} \beta.
  \eequationNN
  Since subproblem~\eqref{prob.opt.S} with $\Scal=\Scal_{k+1}$ is $(\alpha_{k+1},\beta_{k+1})$-strongly morse, it follows that
  \bequationNN
    |d_{\Scal_{k+1}}^T \nabla^2_{xx} L_{\Scal_{k+1}}(x_{\Scal_k},y_{\Scal_{k+1}}) d_{\Scal_{k+1}}| \geq \beta_{\Scal_{k+1}} \|d_{\Scal_{k+1}}\|^2 \geq \tfrac{2}{3}\beta \|d_{\Scal_{k+1}}\|^2.
  \eequationNN
  That said, since \eqref{ineq.var.k+1} holds and $-\tfrac49 > -\tfrac23$, it must hold that
  \bequationNN
    d_{\Scal_{k+1}}^T \nabla^2_{xx} L_{\Scal_{k+1}}(x_{\Scal_k},y_{\Scal_{k+1}}) d_{\Scal_{k+1}} \geq \beta_{\Scal_{k+1}} \|d_{\Scal_{k+1}}\|^2,
  \eequationNN
  which completes the proof.
\eproof

Our next aim is to employ~\cite[Theorem 3.5]{GoyeEfteBoum2024} to prove that \cite[Algorithm 1]{GoyeEfteBoum2024} possesses a certain worst-case iteration complexity bound when employed as the subproblem solver in our Algorithm~\ref{alg.psm}.  For reference in our subsequent analysis, we introduce, for any $\Scal \subseteq [N]$, Fletcher's augmented Lagrangian function with respect to our subproblem~\eqref{prob.opt.S}, namely, $F_\Scal : \R{n} \to \R{}$ defined by
\bequation\label{eq.fletcher.al}
  F_\Scal(x) = f(x) + c_\Scal(x)^Ty_\Scal(x) + \rho_\Scal \|c_\Scal(x)\|^2,
\eequation
where $y_\Scal(x)$ is defined as in \eqref{eq.lsm} and $\rho_\Scal \in \R{}_{>0}$ is a penalty parameter.

\begin{lemma}\label{lemma.fletures.lag.result}
  Suppose there exists $R \in \R{}_{>0}$ such that for any $\Scal \subseteq [N]$ with $|\Scal| \geq p_1$ the set $\Ccal_{\Scal,R} := \{x \in \R{n} : \|c_\Scal(x)\| \leq R\}$ contains $x_0$ and is compact.  Then, for any such $\Scal$, there exists $\hat\rho_\Scal \in \R{}_{>0}$ such that the requirements of \cite[Theorem 3.5]{GoyeEfteBoum2024} hold for all $\rho_\Scal \geq \hat\rho_\Scal$.  Thus, for any such $\rho_\Scal$ there exists $(u_{\Scal,1},u_{\Scal,2}) \in \R{}_{>0} \times \R{}_{>0}$ such that for any $(\epsilon,\zeta) \in (0,\tfrac{\sqrt{5}}{2}R] \times (0,1]$ one has that~\cite[Algorithm 1]{GoyeEfteBoum2024} with starting point $x_0$ locates an $(\epsilon,\zeta)$-stationary point $($see~\eqref{eq.soc.S.approx}$)$ in a number of iterations that is at most
  \bequation\label{eq.fleture.complexity.second}
    T_\Scal = \max\{u_{\Scal,1} \epsilon^{-2}, u_{\Scal,2} \zeta^{-3}\}.
  \eequation
\end{lemma}
\bproof
  To prove the lemma, it suffices to prove that the requirements of \cite[Theorem 3.5]{GoyeEfteBoum2024} hold in our present setting and that the worst-case complexity bound from that theorem holds with respect to our termination condition~\eqref{eq.soc.S.approx}, which is different from the termination condition for~\cite[Algorithm 1]{GoyeEfteBoum2024}.  Toward these ends, let us begin by considering an arbitrary sample set $\Scal \subseteq [N]$ with $|\Scal| \geq p_1$.

  Assumption~A1 in \cite{GoyeEfteBoum2024} requires the existence of $\Rhat \in \R{}_{>0}$ and $\hat\sigma \in \R{}_{>0}$ such that for all $x \in \R{n}$ with $\|c_\Scal(x)\| \leq \Rhat$ one has $\sigma_{\min}(\nabla c_\Scal(x)^T) \geq \hat\sigma$.  This holds in our present setting under the conditions of this lemma and Assumption~\ref{ass.subproblems}.  Assumption~A2 in \cite{GoyeEfteBoum2024} requires the existence of $R$ as stated in the conditions of this lemma. Assumption~A3 in \cite{GoyeEfteBoum2024} requires the existence of $\eta_\Scal \in \R{}_{>0}$ such that
  \bequation\label{ineq.A3}
    \|c_\Scal(x + d) - c_\Scal(x) - \nabla c_\Scal(x)^Td\| \leq \eta_\Scal \|d\|^2\ \ \text{for all}\ \ (x,d) \in \Ccal_{\Scal,R} \times \R{n}.
  \eequation
  To show that this holds in the present setting, observe that for all $(x,j) \in \R{n} \times [m]$ it follows from Assumption~\ref{ass.bounded.distribute}, \eqref{lemma.saa.Hess_c}, and the triangle inequality that
  \begin{align}
    \|\nabla^2 [c_\Scal]_j(x)\|
      &=    \|\nabla^2 [c]_j(x) + \nabla^2 [c_\Scal]_j(x) - \nabla^2 [c]_j(x)\| \nonumber \\
      &\leq \|\nabla^2 [c]_j(x)\| + \|\nabla^2 [c_\Scal]_j(x) - \nabla^2 [c]_j(x)\| \nonumber \\
      &\leq \|\nabla^2 [c]_j(x)\| + \xi_\Scal \sqrt{\gamma_{\nabla^2 c}} \nonumber \\
      &\leq \kappa_{\nabla^2 c} + \sqrt{\tfrac{N(N-p_1)\gamma_{\nabla^2c}}{p_1^2}} =: \kappa_{p_1}. \label{ineq.used}
  \end{align}
  Now consider arbitrary $(x,d) \in \Ccal_{\Scal,R} \times \R{n}$ and observe that, by Taylor's theorem and \eqref{ineq.used}, it follows that for all $j \in [m]$ there exists a point $\xtilde_j \in \R{n}$ such that
  \begin{align}
    &\ |[c_\Scal]_j(x + d) - [c_\Scal]_j(x) - \nabla [c_\Scal]_j(x)^Td| \nonumber \\
    =&\ |\thalf d^T \nabla^2 [c_\Scal]_j(\xtilde_j)d| \leq \thalf \|\nabla^2 [c_\Scal]_j(\xtilde_j)\| \|d\|^2 \leq \thalf \kappa_{p_1} \|d\|^2. \label{ineq.taylor}
  \end{align}
  Consequently, one finds that
  \bequationNN
    \|c_\Scal(x + d) - c_\Scal(x) - \nabla c_\Scal(x)^Td\| \leq \sqrt{\sum_{j \in [m]} \tfrac{1}{4} \kappa_{p_1}^2 \|d\|^4} = \thalf \sqrt{m} \kappa_{p_1} \|d\|^2,
  \eequationNN
  from which it follows that \eqref{ineq.A3} holds for any $\eta_\Scal = \thalf \sqrt{m} \kappa_{p_1}$.  Assumption~A4 in~\cite{GoyeEfteBoum2024} requires that $x_0 \in \Ccal_{\Scal,R}$, as is required in the conditions of this lemma.  Finally, Assumption~A5 in \cite{GoyeEfteBoum2024} requires that the penalty parameter is chosen greater than
  \bequationNN
    \max_{x\in\Ccal_{\Scal,R}} \max \left\{ \tfrac{\sigma_{\max}(\nabla c_\Scal(x)) \sigma_{\max}(\nabla y_\Scal(x))}{2\sigma_{\min}(\nabla c_\Scal(x))^2},
\tfrac{\sigma_{\max}(\nabla y_\Scal(x))}{\sigma_{\min}(\nabla c_\Scal(x))},
\tfrac{1}{\sigma_{\min}(\nabla c_\Scal(x))} \right\}.
  \eequationNN
  Since $\Ccal_{\Scal,R}$ is compact, it follows under Assumption~\ref{ass.subproblems} and the extreme value theorem that there exists $\rho_\Scal$ as stated in the lemma.  All together, one can conclude that the requirements of \cite[Theorem 3.5]{GoyeEfteBoum2024} hold in our present setting, as desired.

  All that remains is to prove that the worst-case iteration complexity bound from \cite[Theorem 3.5]{GoyeEfteBoum2024} yields the desired conclusion of the lemma for our setting.  Toward this end, let us introduce $\tau_\Scal \in \R{}_{>0}$ as equal to the positive real number ``$C$'' introduced in \cite[Corollary~2.7]{GoyeEfteBoum2024} (dependent on second-order derivatives of the element functions $\{[c_{\Scal}]_j\}_{j\in[m]}$ and $\{[y_\Scal]_j\}_{j\in[m]}$ over $\Ccal_{\Scal,R}$).  One can now state by \cite[Theorem 3.5]{GoyeEfteBoum2024} that with respect to \cite[Algorithm 1]{GoyeEfteBoum2024} employed to solve subproblem~\eqref{prob.opt.S} with $x_0 \in \Ccal_{\Scal,R}$ that there exists $(\ubar_{\Scal,1}, \ubar_{\Scal,2}, \ubar_{\Scal,3}) \in \R{}_{>0} \times \R{}_{>0} \times \R{}_{>0}$ such that for any $(\epsilon,\zeta) \in (0,\tfrac{\sqrt{5}}{2}R] \times (0,1]$ and corresponding $\epsilon_{F_\Scal} := \min\{\tfrac{\sqrt{5}\epsilon}{5}, \tfrac{\zeta}{100\tau_\Scal}\}$ and $\zeta_{F_\Scal} := \tfrac{99\zeta}{100} + \tau_\Scal \epsilon_{F_\Scal}$ one finds that in a number of iterations that is at most
  \bequation\label{eq.TS}
    \overline{T}_\Scal := \max \left\{ \ubar_{\Scal,1} \epsilon_{F_\Scal}^{-2}, \ubar_{\Scal,2} \zeta_{F_\Scal}^{-1}, \ubar_{\Scal,3} \zeta_{F_\Scal}^{-3} \right\}
  \eequation
  the algorithm produces a point $\xbar \in \R{n}$ satisfying
  \bequation\label{ineq.fletcher.termin.second}
    \|\nabla F_\Scal(\xbar)\| \leq \epsilon_{F_\Scal}\ \ \text{and}\ \ d^T\nabla^2 F_\Scal(\xbar) d \geq -\zeta_{F_\Scal} \|d\|^2\ \ \text{for all}\ \ d \in \R{n}.
  \eequation
  (Note that to apply~\cite[Theorem 3.5]{GoyeEfteBoum2024} it has been observed that the former tolerance satisfies $\epsilon_{F_\Scal} \leq \tfrac{\sqrt{5}}{5} \epsilon \leq \thalf R$.)  In addition, it also follows from \cite[Theorem 3.5]{GoyeEfteBoum2024} that the point $\xbar$ satisfies $\|\nabla_x L(\xbar,y(\xbar))\| \leq 2\epsilon_{F_\Scal}$, $\|\nabla_y L(\xbar,y(\xbar))\| \leq \epsilon_{F_\Scal}$, and
  \bequationNN
    d^T\nabla_{xx}^2 L(\xbar,y(\xbar))d \geq -\zeta_{F_\Scal} \|d\|^2\ \ \text{for all}\ \ d \in \Null(\nabla c_\Scal(\xbar)^T).
  \eequationNN
  Observing that this means $\|\nabla L(\xbar,y(\xbar))\| \leq \sqrt{5} \epsilon_{F_\Scal} \leq \epsilon$, and observing that $\zeta_{F_\Scal} = \tfrac{99\zeta}{100} + \tau_\Scal \epsilon_{F_\Scal} \leq \zeta$, it follows that \eqref{eq.soc.S.approx} holds, as desired.  Finally, to show that this is achieved in a number of iterations that is at most $T_\Scal$ of the form in~\eqref{eq.fleture.complexity.second}, one only needs to plug in the definitions of $\epsilon_{F_\Scal}$ and $\zeta_{F_\Scal}$ into $\overline{T}_\Scal$ in \eqref{eq.TS}.  This yields 
  \begin{align*}
    \overline{T}_\Scal
      &\leq \max \{5 \ubar_{\Scal,1} \epsilon^{-2}, (100 \tau_\Scal)^2 \ubar_{\Scal,1} \zeta^{-2}, \ubar_{\Scal,2} (\tfrac{99\zeta}{100} + \tau_{\Scal} \tfrac{\sqrt{5}\epsilon}{5})^{-1}, \ubar_{\Scal,2} \zeta^{-1}, \\
      &\qquad\qquad \ubar_{\Scal,3} (\tfrac{99\zeta}{100} + \tau_\Scal \tfrac{\sqrt{5}\epsilon}{5})^{-3}, \ubar_{\Scal,3} \zeta^{-3} \} \\
      &\leq \max \{5 \ubar_{\Scal,1} \epsilon^{-2}, (100 \tau_\Scal)^2 \ubar_{\Scal,1} \zeta^{-2}, \ubar_{\Scal,2} (\tfrac{99\zeta}{100})^{-1}, \ubar_{\Scal,3} (\tfrac{99\zeta}{100})^{-3} \}.
  \end{align*}
  Thus, since $\zeta \in (0,1]$ implies $\max \{\zeta^{-1}, \zeta^{-2}\} \leq \zeta^{-3}$, the conclusion holds with $u_{\Scal,1} = 5 \ubar_{\Scal,1}$ and $u_{\Scal,2} = \max\{ (100 \tau_\Scal)^2 \ubar_{\Scal,1}, \tfrac{100}{99} \ubar_{\Scal,2},(\tfrac{100}{99})^3 \ubar_{\Scal,3} \}$.
\eproof

Since \cite{GoyeEfteBoum2024} is based on Fletcher's Augmented Lagrangian function, we now show that \eqref{ineq.init.next} offers bounds on first- and second-order derivatives of \eqref{eq.fletcher.al}.

\begin{lemma}\label{lemma.ourtermin.fletchertermin}
  Suppose there exists $R \in \R{}_{>0}$ such that for any $\Scal \subseteq [N]$ with $|\Scal| \geq p_1$ the sublevel set $\Ccal_{\Scal,R} := \{x \in \R{n} : \|c_\Scal(x)\| \leq R\}$ contains $x_0$ and is compact.  In addition, for any such $\Scal$ and any given $\beta_\Scal \in \R{}_{>0}$, define  
  \bequation\label{eq.def.kappa.scal}
    \baligned
      \kappa_{\Scal} &:= \max_{x \in \Ccal_{\Scal,R}} \bigg\{ \|\nabla f(x)\|, \|y_\Scal(x)\|, \|\nabla y_\Scal(x)\|, \|\nabla c_\Scal(x)\|, \\
      &\qquad\qquad\qquad \max_{j\in[m]} \left\{ \|\nabla^2 [y_\Scal]_j(x)\|, \|\nabla^2 [c_\Scal]_j(x)\| \right\} \bigg\}, \\
      \eta_{\Scal,1} &:= \(\tfrac{2\sqrt{m}}{\sigma_{\min}} + m\) \kappa_{\Scal}, \\
      \eta_{\Scal,2} &:= 2 m \kappa_\Scal,\\ 
      \eta_{\Scal,3} &:= \tfrac{\thalf \beta_\Scal \sigma_{\min}^2}{\eta_{\Scal,1} \sigma_{\min}^2 + \thalf \beta_\Scal \eta_{\Scal,2}}, \\
      \eta_{\Scal,4} &:= \eta_{\Scal,1}\eta_{\Scal,3} + \(1 + m \kappa_\Scal\) \kappa_\Scal \\ \text{and}\ \ 
      \overline\epsilon_\Scal &:= \tfrac{\thalf \beta_\Scal \sigma_{\min}^2}{\eta_{\Scal,1} \sigma_{\min}^2 + \thalf \beta_\Scal \eta_{\Scal,2} + \eta_{\Scal,2} \eta_{\Scal,4}}
    \ealigned
  \eequation
  along with the functions $\underline{\rho_\Scal} : (0,\overline\epsilon_\Scal) \to \R{}_{>0}$ and $\overline{\rho_\Scal} : (0,\overline\epsilon_\Scal) \to \R{}_{>0}$ defined by
  \bequation\label{eq.def.rho.func}
    \underline{\rho_\Scal}(\epsilon) = \tfrac{\thalf \beta_\Scal + \eta_{\Scal,4}}{2\sigma^2_{\min} - \eta_{\Scal,2} \epsilon}\ \ \text{and}\ \ \overline{\rho_\Scal}(\epsilon) = \tfrac{\thalf \beta_\Scal - \eta_{\Scal,1} \epsilon}{\eta_{\Scal,2} \epsilon}.
  \eequation
  Then, for any such $\Scal$ and any $\beta_\Scal \in \R{}_{>0}$, the following hold.
  \benumerate
    \item[(a)] for any $\epsilon \in (0,\overline\epsilon_\Scal)$ one has $0 < \underline{\rho_\Scal}(\epsilon) < \overline{\rho_\Scal}(\epsilon)$;
    \item[(b)] for any $\epsilon \in (0,\min\{R,\overline\epsilon_\Scal\})$ and $\rho_\Scal \in (\underline{\rho_\Scal}(\epsilon), \overline{\rho_\Scal}(\epsilon))$, if $x$ yields
    \bequation\label{ineq.init.simple}
      \begin{aligned}
        \|\nabla L_{\Scal}(x, y_{\Scal}(x))\| &\leq \epsilon \\
        \text{and}\ \ d^T \nabla^2_{xx} L_{\Scal}(x,y_{\Scal}(x))d &\geq \beta_{\Scal} \|d\|^2\ \ \text{for all}\ \ d \in \Null(\nabla c_{\Scal}(x)^T)
      \end{aligned}
    \eequation
    then
    \bequation\label{ineq.flecturer.condi}
      \|c_\Scal(x)\| \leq R,\ \ \|\nabla F_\Scal(x)\| \leq (1 + \kappa_{\Scal} + 2 \rho_\Scal \kappa_{\Scal}) \epsilon,\ \ \text{and}\ \ \nabla^2 F_\Scal(x) \succeq \thalf \beta_\Scal I.
    \eequation
  \eenumerate
\end{lemma}
\bproof
  Consider arbitrary $\beta_\Scal \in \R{}_{>0}$ and $\Scal \subseteq [N]$ with $|\Scal| \geq p_1$.  Part~(a) follows from the fact that, for any $\epsilon \in (0,\overline\epsilon_\Scal)$, one finds $\epsilon < \overline\epsilon_\Scal \leq \tfrac{\sigma_{\min}^2}{\eta_{\Scal,2}}$ and
  \begin{align*}
    \left\{\overline\epsilon_\Scal \leq \tfrac{\thalf \beta_\Scal \sigma_{\min}^2}{\eta_{\Scal,1} \sigma_{\min}^2 + \thalf \beta_\Scal \eta_{\Scal,2} + \eta_{\Scal,2} \eta_{\Scal,4}}\right\} \iff& \left\{ \tfrac{\thalf \beta_\Scal + \eta_{\Scal,4}}{\sigma_{\min}^2} \leq \tfrac{\thalf \beta_\Scal - \eta_{\Scal,1} \overline\epsilon_\Scal}{\eta_{\Scal,2} \overline\epsilon_\Scal } \right\} \\
    \Longrightarrow \left\{ \tfrac{ \thalf \beta_\Scal + \eta_{\Scal,4}}{2 \sigma_{\min}^2 - \eta_{\Scal,2} \epsilon} < \tfrac{\thalf \beta_\Scal - \eta_{\Scal,1} \epsilon}{\eta_{\Scal,2} \epsilon} \right\} \iff& \left\{\underline{\rho_\Scal}(\epsilon) < \overline{\rho_\Scal}(\epsilon)\right\}.
  \end{align*}
  Let us now prove part (b).  Toward this end, consider arbitrary $\epsilon \in (0,\min\{R,\overline\epsilon_\Scal\})$, $\rho_\Scal \in (\underline{\rho_\Scal}(\epsilon), \overline{\rho_\Scal}(\epsilon))$, and $x$ satisfying \eqref{ineq.init.simple}.  Thus,
  \bequation\label{ineq.nablax.nablay}
    \max\{\|\nabla f(x) + \nabla c_\Scal(x) y_\Scal(x)\|, \|c_\Scal(x)\|\} \leq  \|\nabla L_{\Scal}(x, y_{\Scal}(x))\| \leq \epsilon \leq R,
  \eequation
  which among other things gives the first inequality in \eqref{ineq.flecturer.condi}.  Consequently, $x \in \Ccal_{\Scal,R}$, so it follows from \eqref{eq.fletcher.al}, the chain rule, the triangle inequality, and submultiplicity of the matrix 2-norm that
  \begin{align}
    \|\nabla F_\Scal(x)\|
      &\leq \|\nabla f(x) + \nabla c_\Scal(x) y_\Scal(x)\| + (\|\nabla y_\Scal(x)\| + 2 \rho_\Scal \|\nabla c_\Scal(x)\|) \|c_\Scal(x)\| \nonumber \\
      &\leq \|\nabla f(x) + \nabla c_\Scal(x) y_\Scal(x)\| + (1 + 2 \rho_\Scal) \kappa_{\Scal} \|c_\Scal(x)\|. \label{eq.fletcher.al.grad}
  \end{align}
  Combining \eqref{eq.fletcher.al.grad} with \eqref{ineq.nablax.nablay}, one obtains the second inequality in \eqref{ineq.flecturer.condi}.  Our final aim is to prove the third inequality in \eqref{ineq.flecturer.condi}.  Toward this end, first note that
  \begin{align}
    &\ \nabla^2 F_\Scal(x) \nonumber \\
    =&\ \nabla^2_{xx} L_\Scal(x, y_\Scal(x)) + \nabla y_\Scal(x) \nabla c_\Scal(x)^T + \nabla c_\Scal(x) \nabla y_\Scal(x)^T \nonumber \\
    &\ + 2\rho_\Scal \nabla c_\Scal(x)^T \nabla c_\Scal(x) + \sum_{j\in[m]} (\nabla^2 [y_\Scal]_j(x) + 2 \rho_\Scal \nabla^2 [c_\Scal]_j(x)) [c_\Scal]_j(x). \label{eq.hess.F}
  \end{align}
  Let us now express this Hessian in an equivalent form involving a decomposition into orthogonal spaces defined by the constraint Jacobian at $x$.  Let us first derive an expression for $\nabla y_\Scal(x)\nabla c_\Scal(x)^T$ by differentiating the linear system that defines $y_\Scal(x)$ through \eqref{eq.lsm}.  Specifically, by \eqref{eq.lsm}, one finds that
  \bequation\label{eq.grad.y}
    \nabla c_\Scal(x)^T\nabla c_\Scal(x) y_\Scal(x) = -\nabla c_\Scal(x)^T \nabla f(x).
  \eequation
  Abbreviating notation and differentiating the left-hand side yields
  \begin{align*}
    &\ \nabla(\nabla c_\Scal^T\nabla c_\Scal y_\Scal) |_x \\
    =&\ \nabla ( (\nabla c_\Scal^T|_x \nabla c_\Scal|_x ) y_\Scal) |_x + \nabla (\nabla c_\Scal^T \nabla c_\Scal|_x y_\Scal|_x ) |_x + \nabla (\nabla c_\Scal^T|_x \nabla c_\Scal y_\Scal |_x) |_x \\
    =&\ \nabla y_\Scal|_x \nabla c_\Scal^T |_x \nabla c_\Scal|_x + \bbmatrix \nabla^2 [c_\Scal]_1|_x \nabla c_\Scal|_x y_\Scal|_x & \cdots
    & \nabla^2 [c_\Scal]_m|_x \nabla c_\Scal|_x y_\Scal|_x
    \ebmatrix \\
    &\ + \(\sum_{j\in[m]} \nabla^2 [c_\Scal]_j|_x [y_\Scal]_j|_x \) \nabla c_\Scal|_x,
  \end{align*}
  while at the same time differentiating the right-hand side yields
  \begin{align*}
    \nabla ( -\nabla c_\Scal^T\nabla f) |_x
    =&\ -\nabla (\nabla c_\Scal^T|_x \nabla f)|_x - \nabla (\nabla c_\Scal^T\nabla f|_x) |_x \\
    =&\ -\nabla^2 f|_x \nabla c_\Scal|_x -\bbmatrix \nabla^2 [c_\Scal]_1|_x \nabla f|_x & \cdots & \nabla^2 [c_\Scal]_m|_x \nabla f|_x \ebmatrix.
  \end{align*}
  Combining these derivations and rearranging yields
  \begin{align}
    &\ \nabla y_\Scal(x) \nabla c_\Scal^T(x) \nabla c_\Scal(x) \nonumber \\
    =&\ -\nabla^2_{xx} L_\Scal(x,y_\Scal(x)) \nabla c_\Scal(x) \nonumber \\
    &\ - \underbrace{\bbmatrix \nabla^2 [c_\Scal]_1(x) \nabla_x L_\Scal(x,y_\Scal(x)) & \cdots & \nabla^2 [c_\Scal]_m(x) \nabla_x L_\Scal(x,y_\Scal(x)) \ebmatrix}_{=: \Ecal(x)}. \label{ineq.important}
  \end{align}
  Multiplying \eqref{ineq.important} on the right by $(\nabla c_\Scal(x)^T \nabla c_\Scal(x))^{-1} \nabla c_\Scal(x)^T$ yields
  \begin{align}
    &\ \nabla y_\Scal(x) \nabla c_\Scal(x)^T \nonumber \\
    =&\ -(\nabla^2_{xx} L_\Scal(x,y_\Scal(x)) \nabla c_\Scal(x) + \Ecal(x)) (\nabla c_\Scal(x)^T \nabla c_\Scal(x))^{-1} \nabla c_\Scal(x)^T \nonumber \\
    =&\ -\nabla^2_{xx} L_\Scal(x,y_\Scal(x)) \Rcal(\nabla c_\Scal(x)) -\Ecal(x) \nabla c_\Scal(x)^\dag. \label{eq.hess.e}
  \end{align}
  On the other hand, denoting $\Rcal(x) := \Rcal(\nabla c_\Scal(x))$, $\Ncal(x) := \Ncal(\nabla c_\Scal(x)^T)$, and $H_{xx} := \nabla^2_{xx} L_\Scal(x,y_\Scal(x))$ for the sake of notational simplicity and using the fact that $\Rcal(x) + \Ncal(x) = I$, one finds that
  \begin{align}
    H_{xx} - H_{xx}\Rcal(x) - \Rcal (x)H_{xx}
    &= H_{xx} \Ncal(x) - \Rcal(x) H_{xx} \nonumber \\
    &=(\Ncal(x) + \Rcal(x)) H_{xx} \Ncal(x) - \Rcal(x) H_{xx} \nonumber \\
    &= \Ncal(x) H_{xx} \Ncal(x) - \Rcal(x) H_{xx} \Rcal(x). \label{eq.hess.simp}
  \end{align}
  Combining \eqref{eq.hess.F}, \eqref{eq.hess.e}, and \eqref{eq.hess.simp} now yields the expression for the Hessian as
  \begin{align}
    \nabla^2 F_\Scal(x)
    =&\ \Ncal(x) H_{xx} \Ncal(x) -\Rcal(x) H_{xx} \Rcal(x) + 2 \rho_\Scal \nabla c_\Scal(x)^T \nabla c_\Scal(x) \nonumber \\
    &\ -\Ecal(x) \nabla c_\Scal(x)^\dag -(\nabla c_\Scal(x)^\dag)^T \Ecal(x)^T \nonumber \\
    &\ + \sum_{j\in[m]} (\nabla^2 [y_\Scal]_j(x) + 2 \rho_\Scal \nabla^2 [c_\Scal]_j(x)) [c_\Scal]_j(x). \label{eq.fl.hess.rewrite}
  \end{align}
  Let us now observe that from the definition of $\Ecal(x)$, norm inequalities, submultiplicity of the matrix 2-norm, the fact that $x \in \Ccal_{\Scal,R}$, \eqref{eq.def.kappa.scal}, and \eqref{ineq.nablax.nablay} that
  \begin{align}
    \|\Ecal(x)\|
    &= \|\Ecal(x)^T\| \nonumber \\
    &= \max_{d\in\R{n} \st \|d\|=1} \|\Ecal(x)^Td\| \nonumber \\
    &= \max_{d\in\R{n} \st \|d\|=1} \left\| \bbmatrix \nabla_x L_\Scal(x,y_\Scal(x))^T \nabla^2 [c_\Scal]_1(x)^Td \\ \vdots \\ \nabla_x L_\Scal(x,y_\Scal(x))^T \nabla^2 [c_\Scal]_m(x)^Td \ebmatrix \right\| \nonumber \\
    &= \max_{d\in\R{n} \st \|d\|=1} \sqrt{\sum_{j\in[m]} \( \nabla_x L_\Scal(x,y_\Scal(x))^T \nabla^2 [c_\Scal]_j(x)^Td\)^2} \nonumber \\
    &\leq \sqrt{m} \max_{d\in\R{n} \st \|d\|=1} \(\max_{j\in[m]} |\nabla_x L_\Scal(x,y_\Scal(x))^T \nabla^2 [c_\Scal]_j(x)^Td| \) \nonumber \\
    &\leq \sqrt{m} \max_{d\in\R{n} \st \|d\|=1} \(\max_{j\in[m]} \|\nabla_x L_\Scal(x,y_\Scal(x))\| \|\nabla^2 [c_\Scal]_j(x)\| \|d\|\) \nonumber \\
    &\leq \sqrt{m} \max_{d\in\R{n} \st \|d\|=1} (\epsilon\kappa_\Scal\|d\|) = \sqrt{m} \kappa_\Scal \epsilon. \label{ineq.norm.Ecal}
  \end{align}
  Now consider arbitrary $d \in \R{n} \setminus \{0\}$ decomposed as $d \equiv d_\Ncal + d_\Rcal$, where $d_\Ncal \in \Null(\nabla c_\Scal(x)^T)$ and $d_\Rcal \in \Range(\nabla c_\Scal(x))$.  One has from \eqref{eq.fl.hess.rewrite} that
  \begin{align}
    &\ d^T\nabla^2 F_\Scal(x) d \nonumber \\
    =&\ d^T\Ncal(x) \nabla^2_{xx} L_\Scal(x,y_\Scal(x)) \Ncal(x) d - d^T \Rcal(x) \nabla^2_{xx} L_\Scal(x,y_\Scal(x)) \Rcal(x) d \nonumber \\
    &\ + 2\rho_\Scal d^T \nabla c_\Scal(x)^T\nabla c_\Scal(x) d - 2 d^T \Ecal(x) \nabla c_\Scal(x)^\dag d \nonumber \\
    &\ + d^T\(\sum_{j\in[m]} \(\nabla^2 [y_\Scal]_j(x) + 2\rho_\Scal \nabla^2 [c_\Scal]_j(x)\) [c_\Scal]_j(x)\)d \nonumber \\
    =&\ d_\Ncal^T \nabla^2_{xx} L_\Scal(x,y_\Scal(x)) d_\Ncal - d_\Rcal^T \nabla^2_{xx} L_\Scal(x,y_\Scal(x)) d_\Rcal \nonumber \\
    &\ + 2\rho_\Scal \|\nabla c_\Scal(x) d_\Rcal\|^2 - 2 d^T\Ecal(x) \nabla c_\Scal(x)^\dag d \nonumber \\
    &\ + d^T \(\sum_{j\in[m]} \(\nabla^2 [y_\Scal]_j(x) + 2\rho_\Scal \nabla^2 [c_\Scal]_j(x)\) [c_\Scal]_j(x)\) d, \label{ineq.hess.dd}
  \end{align}
  where for the latter two terms one has from \eqref{eq.def.kappa.scal}, \eqref{ineq.nablax.nablay}, and \eqref{ineq.norm.Ecal} that
  \begin{align}
    &\ -2d^T \Ecal(x) \nabla c_\Scal(x)^\dag d + d^T \(\sum_{j\in[m]} \(\nabla^2 [y_\Scal]_j(x) + 2\rho_\Scal \nabla^2 [c_\Scal]_j(x)\) [c_\Scal]_j(x)\)d \nonumber \\
    \geq&\ -2 \|\Ecal(x)\| \|\nabla c_\Scal(x)^\dag\| \|d\|^2 - \left\| \sum_{j\in[m]} (\nabla^2 [y_\Scal]_j(x) + 2\rho_\Scal \nabla^2 [c_\Scal]_j(x)) [ c_\Scal]_j(x) \right\| \|d\|^2 \nonumber \\
    \geq&\ -\tfrac{2\sqrt{m} \kappa_\Scal}{\sigma_{\min}} \|d\|^2 \epsilon -\sum_{j\in[m]} ( \|\nabla^2 [y_\Scal]_j(x)\| + 2\rho_\Scal \|\nabla^2 [c_\Scal]_j(x)\|) |[c_\Scal]_j(x)| \|d\|^2 \nonumber \\ 
    \geq&\ -\tfrac{2\sqrt{m}\kappa_\Scal}{\sigma_{\min}} \|d\|^2 \epsilon -\sum_{j\in[m]} ( \|\nabla^2 [y_\Scal]_j(x)\| + 2\rho_\Scal \|\nabla^2 [c_\Scal]_j(x)\| ) \|c_\Scal(x)\| \|d\|^2 \nonumber \\
    \geq&\ -\tfrac{2\sqrt{m}\kappa_\Scal}{\sigma_{\min}} \|d\|^2 \epsilon -m\kappa_{\Scal} (1 + 2\rho_\Scal) \epsilon \|d\|^2. \label{ineq.bound.sum}
  \end{align}
  Now since $d_\Ncal \in \Null(\nabla c_\Scal(x)^T)$, \eqref{ineq.init.simple} gives $d_\Ncal^T \nabla_{xx}^2 L_\Scal(x,y_\Scal(x)) d_\Ncal \geq \beta_\Scal \|d_\Ncal\|^2$.  Along with \eqref{eq.def.kappa.scal}, the triangle inequality, and submultiplicity of the matrix 2-norm,
  \begin{align}
    \|\nabla^2_{xx} L_\Scal(x,y_\Scal(x))\|
      &= \left\|\nabla^2f(x) + \sum_{j\in[m]} \nabla^2 [c_\Scal]_j(x)[y_\Scal]_j(x) \right\| \nonumber \\
      &\leq \|\nabla^2 f(x)\| + \sum_{j\in[m]} \|\nabla^2 [c_\Scal]_j(x)\| \|y_\Scal(x)\| \leq \kappa_\Scal + m \kappa_\Scal^2. \label{ineq.hess.lag}
  \end{align}
  At the same time, let us derive a lower bound for $\|\nabla c_\Scal(x)^Td_\Rcal\|$.  Let $\nabla c_\Scal(x)$ have the singular value decomposition $\sum_{i\in[m]} u_i \sigma_i v_i^T$ where $(u_i,\sigma_i,v_i) \in \R{n} \times \R{} \times \R{m}$ for all $i\in[m]$ with $\{u_i\}_{i\in[m]}$ and $\{v_i\}_{i\in[m]}$ being sets of orthonormal vectors. Then, by Assumption~\ref{ass.subproblems} and the fact that $d_\Rcal\in\Range(\nabla c_\Scal(x))$, one finds that
  \begin{align}
    \|\nabla c_\Scal(x)^Td_\Rcal\|^2
    &= \left\| \sum_{i\in[m]} v_i \sigma_i u_i^T d_\Rcal \right\|^2 = \sum_{i\in[m]} \|v_i \sigma_i u_i^T d_\Rcal\|^2 = \sum_{i\in[m]} \sigma_i^2 (u_i^T d_\Rcal)^2 \nonumber \\
&\ge\sigma_{\min}^2\sum_{i=1}^m(u_i^Td_\Rcal)^2=\sigma_{\min}^2\|d_\Rcal\|^2. \label{ineq.nabla cd}
  \end{align}
  Combining \eqref{ineq.hess.dd}, \eqref{ineq.bound.sum}, \eqref{ineq.hess.lag}, \eqref{ineq.nabla cd}, and the decomposition $d = d_\Ncal + d_\Rcal$,
  \begin{align}
    &\ d^T\nabla^2 F_\Scal(x) d \nonumber \\
    \geq&\ \beta_\Scal \|d_\Ncal\|^2 + (2\rho_\Scal\sigma_{\min}^2 - (\kappa_\Scal + m \kappa_\Scal^2 )) \|d_\Rcal\|^2 \nonumber \\
    &\ - \(\tfrac{2 \sqrt{m}}{\sigma_{\min}} + m (1 + 2\rho_\Scal ) \) \kappa_{\Scal} \|d\|^2 \epsilon \nonumber \\
    =&\ \(\beta_\Scal - \(\tfrac{2\sqrt{m}}{\sigma_{\min}} + m (1 + 2\rho_\Scal) \) \kappa_{\Scal} \epsilon\) \|d_\Ncal\|^2 \nonumber \\
    &\ + \(2\rho_\Scal\sigma_{\min}^2 - (1 + m \kappa_\Scal) \kappa_\Scal - \(\tfrac{2\sqrt{m}}{\sigma_{\min}} + m (1 + 2 \rho_\Scal) \) \kappa_{\Scal} \epsilon\) \|d_\Rcal\|^2. \label{ineq.final}
  \end{align}
  We now claim that our desired conclusion follows as long as $(\rho_\Scal,\epsilon)$ yields
  \bequation\label{ineq.whateever}
    \baligned
      \beta_\Scal - \(\tfrac{2\sqrt{m}}{\sigma_{\min}} + m (1 + 2\rho_\Scal) \) \kappa_{\Scal} \epsilon &\geq \thalf \beta_\Scal\ \ \text{and} \\
      2\rho_\Scal\sigma_{\min}^2 - (1 + m \kappa_\Scal) \kappa_\Scal - \(\tfrac{2\sqrt{m}}{\sigma_{\min}} + m (1 + 2 \rho_\Scal) \) \kappa_{\Scal} \epsilon &\geq \thalf\beta_\Scal,
    \ealigned
  \eequation
  since these inequalities along with \eqref{ineq.final} would yield the desired fact that
  \bequationNN
    d^T \nabla^2 F_\Scal(x) d \geq \thalf \beta_\Scal \|d_\Ncal\|^2 + \thalf\beta_\Scal \|d_\Rcal\|^2 = \thalf\beta_\Scal\|d\|^2.
  \eequationNN
  Indeed, since $\epsilon \in (0,\overline\epsilon_\Scal)$ and $\rho_\Scal \in (\underline{\rho_\Scal}(\epsilon), \overline{\rho_\Scal}(\epsilon))$, one finds that
  \begin{align*}
    \beta_\Scal - \(\tfrac{2\sqrt{m}}{\sigma_{\min}} + m (1 + 2 \rho_\Scal) \) \kappa_{\Scal} \epsilon
    &= \beta_\Scal - (\eta_{\Scal,1} \epsilon + \eta_{\Scal,2} \rho_\Scal \epsilon) \\
    &\geq \beta_\Scal - (\eta_{\Scal,1} \epsilon + \eta_{\Scal,2} \overline{\rho_\Scal}(\epsilon) \epsilon) \\
    &= \beta_\Scal - (\eta_{\Scal,1} \epsilon + \thalf \beta_\Scal-\eta_{\Scal,1} \epsilon ) = \thalf\beta_\Scal
  \end{align*}
  as well as
  \begin{align*}
    &\ 2 \rho_\Scal \sigma_{\min}^2 - (1 + m \kappa_\Scal) \kappa_\Scal - \(\tfrac{2\sqrt{m}}{\sigma_{\min}} + m (1 + 2\rho_\Scal)\) \kappa_{\Scal} \epsilon \\
    =&\ 2 \rho_\Scal \sigma_{\min}^2 - \(\eta_{\Scal,1}\epsilon+(1 + m \kappa_\Scal) \kappa_\Scal\) - \eta_{\Scal,2} \rho_\Scal \epsilon \\
    \geq&\ 2 \rho_\Scal \sigma_{\min}^2 - \(\eta_{\Scal,1}\overline\epsilon_\Scal+(1 + m \kappa_\Scal) \kappa_\Scal\) - \eta_{\Scal,2} \rho_\Scal \epsilon \\
    \ge&\ 2 \rho_\Scal \sigma_{\min}^2 - \(\eta_{\Scal,1}\eta_{\Scal,3}+(1 + m \kappa_\Scal) \kappa_\Scal\) - \eta_{\Scal,2} \rho_\Scal \epsilon \\
    =&\ 2 \rho_\Scal \sigma_{\min}^2 - \eta_{\Scal,4} - \eta_{\Scal,2} \rho_\Scal \epsilon \\
    \geq&\ 2 \underline{\rho_\Scal}(\epsilon) \sigma_{\min}^2 - \eta_{\Scal,4} - \eta_{\Scal,2} \overline{\rho_\Scal}(\epsilon) \epsilon \\
    =&\ \tfrac{(\beta_\Scal + 2\eta_{\Scal,4}) \sigma_{\min}^2}{2 \sigma_{\min}^2 - \eta_{\Scal,2} \epsilon} - \eta_{\Scal,4} - \(\thalf \beta_\Scal - \eta_{\Scal,1} \epsilon\) \\
    \geq&\ \tfrac{\(\beta_\Scal + 2\eta_{\Scal,4}\) \sigma_{\min}^2}{\sigma_{\min}^2} - 2 \eta_{\Scal,4} - \thalf\beta_\Scal = \thalf\beta_\Scal.
  \end{align*}
  For the reasons stated previously, the proof is complete.
\eproof

To obtain our desired complexity result, we need one more assumption.

\bassumption\label{ass.Hessian.Lip.F}
  For any $\Scal \subseteq [N]$ with $|\Scal| \geq p_1$ and any $\rho_\Scal \in \R{}_{>0}$, Fletcher's Augmented Lagrangian function has a Hessian function $\nabla^2 F_\Scal : \R{n} \to \R{n \times n}$ that is Lipschitz continuous in the sense that there exists $M_\Scal \in \R{}_{>0}$ such that
  \bequationNN
    \|\nabla^2 F_{\Scal}(x) - \nabla^2 F_{\Scal}(\xbar)\| \leq M_\Scal \|x - \xbar\|^2\ \ \text{for all}\ \ (x,\xbar) \in \R{n} \times \R{n}.
  \eequationNN
\eassumption

\begin{lemma}\label{lemma.fl.strongly.convex}
  Suppose that the conditions of Lemma~\ref{lem.start_for_next} hold in the sense that $(\kappa_1,\kappa_2,\kappa_3)$ is defined by \eqref{eq.kappa} and \eqref{eq.def.kappa.3}, $(\tau_1,\tau_2)$ is defined by \eqref{eq.def.tau2}, for all $k \in [K]$ the pair of subproblem tolerances $(\epsilon_k,\zeta_k)$ is defined by \eqref{eq.tolorence}, $p_1$ is sufficiently large such that \eqref{eq.theorem2.S} holds, and for all $k \in [K]$ the point $x_{\Scal_k} \in \R{n}$ is $(\epsilon_k,\zeta_k)$-stationary with respect to~\eqref{prob.opt.S} for $\Scal = \Scal_k$.  In addition, suppose that the conditions of Lemma~\ref{lemma.ourtermin.fletchertermin} hold that there exists $R \in \R{}_{>0}$ such that for any $\Scal \subseteq [N]$ with $|\Scal| \geq p_1$ the set $\Ccal_{\Scal,R} := \{x \in \R{n} : \|c_\Scal(x)\| \leq R\}$ contains $x_0$ and is compact. Further, with \eqref{eq.def.kappa.scal} and \eqref{eq.def.rho.func} for all $\beta_\Scal \geq \thalf \beta$ and $\Scal \subset [N]$ with $|\Scal| \geq p_1$, suppose with
  \bequation\label{eq.def.delta_s}
    \baligned
      \underline{\delta} &:= \min_{\Scal \subseteq [N] \st \eqref{eq.theorem2.S}} \left\{ R, \overline\epsilon_\Scal, \tfrac{\thalf \beta}{2\eta_{\Scal,1} + 3 \eta_{\Scal,2} \tfrac{\max\{1,\kappa_\Scal\}}{\sigma_{\min}}} \right\}, \\
      \overline\omega &:= \max_{\Scal \subseteq [N] \st \eqref{eq.theorem2.S}} 1 + \kappa_\Scal + 2 \overline{\rho_\Scal}(\underline\delta) \kappa_\Scal, \\ \text{and}\ \ 
      \overline{M} &:= \max_{\Scal \subseteq [N] \st \eqref{eq.theorem2.S}} M_{\Scal}\ \text{(see Assumption~\ref{ass.Hessian.Lip.F})}
    \ealigned
  \eequation
  that the sample set $\Scal_k$ for each $k \in [K]$ yields
  \bequation\label{ineq.S.lemma}
    \sqrt{\tfrac{N (N-|\Scal_k|)}{|\Scal_k|^2}} \leq \tfrac{1}{\kappa_2} \min\left\{\tfrac{\underline{\delta}}{3},
    \tfrac{\beta}{18 \overline{\omega} ^2},\tfrac{\beta^2}{432 \overline{\omega} \overline{M}},\tfrac{R}{3\overline{\omega}}\right\}.
  \eequation
  Finally, suppose that with $\rho_{\Scal_{k+1}} \in \R{}_{>0}$ and $t_{\Scal_{k+1}} \in \R{}_{>0}$, gradient descent with constant step size $t_{\Scal_{k+1}}$ is employed to minimize $F_{\Scal_{k+1}}$ with initial point~$x_{\Scal_k}$.  Then, there exists a positive interval and a positive upper bound such that if $\rho_{\Scal_{k+1}}$ is within the positive interval and $t_{\Scal_{k+1}}$ is below the upper bound, in at most
  \bequationNN
    \left \lceil \log_2 \tfrac{3 \sqrt{5} \overline{\omega} \epsilon_k}{\epsilon_{k+1}} \right \rceil\ \text{iterations}
  \eequationNN
  the method gives $x_{\Scal_{k+1}}$ that is $(\epsilon_{k+1}, \zeta_{k+1})$-stationary for \eqref{prob.opt.S} with $\Scal = \Scal_{k+1}$.
\end{lemma}
\bproof
  Consider arbitrary $k \in [K]$.  By Lemma~\ref{lem.start_for_next}, it follows that \eqref{ineq.init.next} holds, where $\alpha_{\Scal_{k+1}} \geq \thalf \alpha$ and $\beta_{\Scal_{k+1}} \geq \thalf \beta$ are defined according to \eqref{eq.alpha-beta-s}.  Consequently, one finds that the result of Lemma~\ref{lemma.ourtermin.fletchertermin}(b) holds with $\Scal = \Scal_{k+1}$, $\epsilon = 3\epsilon_k$, $\beta_{\Scal_{k+1}} \geq \thalf \beta$ defined by \eqref{eq.alpha-beta-s}, $x = x_{\Scal_k}$, and $\rho_{\Scal_{k+1}}=\overline{\rho_{\Scal_{k+1}}}(\underline\delta)$.  Indeed, with \eqref{eq.def.tau2}, \eqref{eq.tolorence}, the first term of the minimum in \eqref{ineq.S.lemma}, and \eqref{eq.def.delta_s}, one finds that $3 \epsilon_k = 3 \kappa_2 \xi_{\Scal_k} \leq \underline{\delta} \leq \min\{R,\overline\epsilon_{\Scal_{k+1}}\}$ and $\rho_{\Scal_{k+1}}\ge\underline{\rho_{\Scal_{k+1}}}(\underline\delta)\ge\underline{\rho_{\Scal_{k+1}}}(3\epsilon_k)$.  Thus,
  \bequation\label{ineq.i-0}
    x_{\Scal_k} \in \Ccal_{\Scal_{k+1},R},\ \|\nabla F_{\Scal_{k+1}}(x_{\Scal_k})\|
\leq 3 \overline{\omega} \epsilon_k\ \text{and}\ \nabla^2 F_{\Scal_{k+1}}(x_{\Scal_k})\succeq \thalf\beta_{\Scal_{k+1}}I.
  \eequation
  Before proceeding, since $\Ccal_{\Scal_{k+1},R} := \{x \in \R{n} : \|c_{\Scal_{k+1}}(x)\| \leq R \}$ is compact, let
  \bequationNN
    \sigma_{\nabla^2 F_{\Scal_{k+1}}} := \max_{x \in \Ccal_{\Scal_{k+1},R}} \{\sigma_{\max}(\nabla^2 F_{\Scal_{k+1}}(x))\} < \infty.
  \eequationNN
  
  Let us now consider the behavior of gradient descent employed with penalty parameter $\rho_{\Scal_{k+1}} \in \R{}_{>0}$ and step size $t_{\Scal_{k+1}} \in \R{}_{>0}$; that is, with $x_{\Scal_k}^0 \gets x_{\Scal_k}$, consider the iterative method defined for all $i = 0,1,2,\dots$ by
  \begin{align}
    x_{\Scal_k}^{i+1} \gets&\ x_{\Scal_k}^i - t_{\Scal_{k+1}} \nabla F_{\Scal_{k+1}}(x_{\Scal_k}^i) \label{eq.fl.gd} \\ \text{where}\ \ 
    t_{\Scal_{k+1}} \gets&\ \min\left\{ \tfrac{\beta_{\Scal_{k+1}}}{9 \overline{\omega} \overline{M} \epsilon_k}, \tfrac{1}{\sigma_{\nabla^2 F_{\Scal_{k+1}}}}, \tfrac{3}{\beta_{\Scal_{k+1}}} \right\}. \label{eq.stepsize}
  \end{align}
  Our next aim is to show that, for any such $i \geq 1$, one has that
  \bequation\label{ineq.induct.end}
    x_{\Scal_k}^i \in \Ccal_{\Scal_{k+1},R},\ \ \|\nabla F_{\Scal_{k+1}}(x_{\Scal_k}^i)\| \leq \(1 - \tfrac{1}{6} \beta_{\Scal_{k+1}} t_{\Scal_{k+1}}\)^i \|\nabla F_{\Scal_{k+1}}(x_{\Scal_k}^{0})\|
  \eequation
  and $\sigma_{\min}(\nabla^2 F_{\Scal_{k+1}}(x_{\Scal_k}^0)) \geq \tfrac12 \beta_{\Scal_{k+1}}$ while for $i \geq 1$ one has
  \bequation\label{ineq.induct.end.2}
    \baligned
      &\ \sigma_{\min}(\nabla^2 F_{\Scal_{k+1}}(x_{\Scal_k}^{i})) \\
      \geq&\ \tfrac{1}{2} \beta_{\Scal_{k+1}} - \overline{M} \|\nabla F_{\Scal_{k+1}}(x_{\Scal_k}^{0})\| t_{\Scal_{k+1}} \sum_{l=0}^{i-1} \(1 - \tfrac{1}{6}\beta_{\Scal_{k+1}} t_{\Scal_{k+1}})\)^l \geq \tfrac{1}{3} \beta_{\Scal_{k+1}}.
    \ealigned
  \eequation
  Toward this end, first observe that by defining $G_i := \nabla F_{\Scal_{k+1}}(x_{\Scal_k}^i)$ for all $i$ one has from Taylor's theorem and submultiplicity that
  \begin{align}
    \|G_{i+1}\|
      &\leq \left\|G_i + \nabla^2 F_{\Scal_{k+1}}(x_{\Scal_k}^i) (x_{\Scal_k}^{i+1} - x_{\Scal_k}^i) \right\| + \tfrac{1}{2} \overline{M} \|x_{\Scal_k}^{i+1}-x_{\Scal_k}^{i}\|^2 \nonumber \\
      &= \|G_i - t_{\Scal_{k+1}} \nabla^2 F_{\Scal_{k+1}}(x_{\Scal_k}^i)G_i\| + \tfrac{1}{2} \overline{M}t_{\Scal_{k+1}}^2 \|G_i\|^2 \nonumber \\
  &\le \|I - t_{\Scal_{k+1}} \nabla^2 F_{\Scal_{k+1}}(x_{\Scal_k}^i)\| \|G_i\| + \tfrac{1}{2} \overline{M}t_{\Scal_{k+1}}^2 \|G_i\|^2. \label{ineq.taylor.exp}
  \end{align}
  Let us now employ the aforementioned induction.  For $i=0$, the base case holds since $x_{\Scal_k}^0 = x_{\Scal_k}$ satisfies \eqref{ineq.i-0}.  Now consider an arbitrary positive integer $i$ and suppose that $x_{\Scal_k}^i$ satisfies \eqref{ineq.induct.end}--\eqref{ineq.induct.end.2}.  By \eqref{eq.stepsize}, one has
  \bequationNN
    t_{\Scal_{k+1}} \sigma_{\max}(\nabla^2 F_{\Scal_{k+1}}(x_{\Scal_k}^i)) \leq t_{\Scal_{k+1}} \sigma_{\nabla^2 F_{\Scal_{k+1}}} \leq 1,
  \eequationNN
  which in turn along with \eqref{ineq.induct.end.2} shows that
  \begin{align}
    &\ \|I - t_{\Scal_{k+1}} \nabla^2 F_{\Scal_{k+1}}(x_{\Scal_k}^i)\| \nonumber \\
    =&\ \max_{v\in\R{n}\st\|v\|=1} \left|v^Tv-t_{\Scal_{k+1}} v^T\nabla^2  F_{\Scal_{k+1}}(x_{\Scal_k}^i)v\right| \nonumber \\
    =&\ \max\left\{ |1 - t_{\Scal_{k+1}} \sigma_{\min}(\nabla^2 F_{\Scal_{k+1}}(x_{\Scal_k}^i))|, |1 - t_{\Scal_{k+1}} \sigma_{\max}(\nabla^2 F_{\Scal_{k+1}}(x_{\Scal_k}^i)) | \right\} \nonumber \\
    =&\ 1 - t_{\Scal_{k+1}} \sigma_{\min}( \nabla^2  F_{\Scal_{k+1}}(x_{\Scal_k}^i)) \nonumber \\
    \leq&\ 1 - \tfrac{1}{3} \beta_{\Scal_{k+1}} t_{\Scal_{k+1}}. \label{ineq.bound.1}
  \end{align}
  At the same time, combining \eqref{ineq.i-0}, \eqref{eq.stepsize}, and \eqref{ineq.induct.end}, one finds that
  \bequation\label{ineq.g0}
    \overline{M} t_{\Scal_{k+1}} \|G_i\| \leq \overline{M} t_{\Scal_{k+1}} \|G_0\| \leq 3 \overline{M} t_{\Scal_{k+1}} \overline{\omega} \epsilon_k \leq \tfrac{1}{3}\beta_{\Scal_{k+1}}.
  \eequation
  Now combining \eqref{ineq.taylor.exp}, \eqref{ineq.bound.1}, and \eqref{ineq.g0}, one obtains that
  \begin{align}
    \|G_{i+1}\|
      &\leq \(1-\tfrac{1}{3} \beta_{\Scal_{k+1}} t_{\Scal_{k+1}}\) \|G_i\| + \tfrac16 \beta_{\Scal_{k+1}} t_{\Scal_{k+1}} \|G_i\| \nonumber \\
      &= \(1 - \tfrac{1}{6} \beta_{\Scal_{k+1}} t_{\Scal_{k+1}}\) \|G_i\| \leq \(1 - \tfrac{1}{6} \beta_{\Scal_{k+1}} t_{\Scal_{k+1}}\)^{i+1} \|G_0\|, \label{ineq.G1.G0}
  \end{align}
  which shows the inequality in~\eqref{ineq.induct.end}, and in fact shows that $\|\nabla F_{\Scal_{k+1}}(\cdot)\|$ decreases with each iteration.  Also, from \cite[Theorem 6.6]{Stew1973} (or see \cite[Eq.~(3)]{Stew1979}), one finds
  \begin{align*}
    &\ |\sigma_{\min}(\nabla^2 F_{\Scal_{k+1}}(x_{\Scal_k}^{i+1}))-\sigma_{\min}(\nabla^2 F_{\Scal_{k+1}}(x_{\Scal_k}^i))| \\
    \leq&\ \|\nabla^2 F_{\Scal_{k+1}}(x_{\Scal_k}^{i+1}))- \nabla^2 F_{\Scal_{k+1}}(x_{\Scal_k}^i)\| \leq \overline{M}\|x_{\Scal_k}^{i+1}-x_{\Scal_k}^i\|= \overline{M}t_{\Scal_{k+1}}\|G_i\|.
  \end{align*}
  Combined with \eqref{eq.stepsize}, \eqref{ineq.induct.end.2}, and \eqref{ineq.G1.G0}, one now obtains that
  \begin{align}
    &\ \sigma_{\min}(\nabla^2 F_{\Scal_{k+1}}(x_{\Scal_k}^{i+1})) \nonumber \\
    \geq&\ \sigma_{\min}(\nabla^2 F_{\Scal_{k+1}}(x_{\Scal_k}^i))-|\sigma_{\min}(\nabla^2 F_{\Scal_{k+1}}(x_{\Scal_k}^{i+1}))-\sigma_{\min}(\nabla^2 F_{\Scal_{k+1}}(x_{\Scal_k}^i))| \nonumber \\
    \geq&\ \tfrac{1}{2}\beta_{\Scal_{k+1}}-\overline{M}\|G_0\|t_{\Scal_{k+1}}\sum_{l=0}^{i-1}\(1-\tfrac{1}{6}\beta_{\Scal_{k+1}}t_{\Scal_{k+1}}\)^{l}-\overline{M}t_{\Scal_{k+1}}\|G_i\| \nonumber \\
    \geq&\ \tfrac{1}{2}\beta_{\Scal_{k+1}}-\overline{M}\|G_0\|t_{\Scal_{k+1}}\sum_{l=0}^i\(1-\tfrac{1}{6}\beta_{\Scal_{k+1}}t_{\Scal_{k+1}}\)^{l} \nonumber \\
    \geq&\ \tfrac{1}{2}\beta_{\Scal_{k+1}}-\overline{M}\|G_0\|t_{\Scal_{k+1}}\sum_{l=0}^\infty\(1-\tfrac{1}{6}\beta_{\Scal_{k+1}}t_{\Scal_{k+1}}\)^{l} \nonumber \\
    =&\ \tfrac{1}{2}\beta_{\Scal_{k+1}}-\overline{M}\|G_0\|t_{\Scal_{k+1}} \tfrac{6}{\beta_{\Scal_{k+1}}t_{\Scal_{k+1}}} \nonumber \\
    \geq&\ \tfrac{1}{2}\beta_{\Scal_{k+1}}-\tfrac{18 \overline{M}\overline{\omega}\epsilon_k}{\beta_{\Scal_{k+1}}} \ge \tfrac{1}{3}\beta_{\Scal_{k+1}}. \label{ineq.strong.convex.2}
  \end{align}
  Here, for the last inequality, note that $\epsilon_k=\tau_1\xi_{\Scal_k}=\kappa_2\xi_{\Scal_k}$ by~\eqref{eq.tolorence}. Then, since $\Scal_k$ satisfies \eqref{ineq.S.lemma}, one finds $\epsilon_k\le \tfrac{\beta^2}{432\overline{M}\overline{\omega}}$, which further gives $\tfrac{18 \overline{M}\overline{\omega}\epsilon_k}{\beta_{\Scal_{k+1}}}\le \tfrac{\beta^2}{24\beta_{\Scal_{k+1}}}$. Now since $\beta_{\Scal_{k+1}}\ge\thalf\beta$ by \eqref{eq.alpha-beta-s}, it follows that $\tfrac{18 \overline{M}\overline{\omega}\epsilon_{\Scal_k}}{\beta_{\Scal_{k+1}}} \le  \tfrac{\beta_{k+1}}{6}$. As a result, one has $\tfrac{1}{2}\beta_{\Scal_{k+1}}-\tfrac{18 \overline{M}\overline{\omega}\epsilon_k}{\beta_{\Scal_{k+1}}} \ge \tfrac{1}{2}\beta_{\Scal_{k+1}}-\tfrac{1}{6}\beta_{\Scal_{k+1}}=\tfrac{1}{3}\beta_{\Scal_{k+1}}$, which gives the last inequality.

  All that remains to complete the induction is to show that $x_{\Scal_k}^{i+1}\in\Ccal_{\Scal_{k+1},R}$.  Toward this end, let us employ \cite[Corollary 2.7]{GoyeEfteBoum2024}.  This requires that the singular values of $\nabla c_{\Scal_{k+1}}$ are bounded away from zero and that $\Ccal_{\Scal_{k+1},R}$ is compact, both of which are assumed here.  It also requires that
  \bequation\label{ineq.ensreu.tau}
    \rho_{\Scal_{k+1}} > \max_{x\in\Ccal_{\Scal_{k+1},R}} \tfrac{ \max\{1,\|\nabla y_{\Scal_{k+1}}(x)\|\}}{\sigma_{\min}(\nabla c_{\Scal_{k+1}}(x)) }.
  \eequation
  Let us show that the choice $\rho_{\Scal_{k+1}} = \overline{\rho_{\Scal_{k+1}}}(\underline\delta)$, which is reflected in \eqref{eq.def.delta_s}, yields this lower bound.  Indeed, since $\sigma_{\min}(\nabla c_{\Scal_{k+1}}(x))\ge\sigma_{\min}$ by Assumption \ref{ass.subproblems} and $\|\nabla y_{\Scal_{k+1}}(x)\|\le\kappa_{\Scal_{k+1}}$ by the definition of $\kappa_{\Scal_{k+1}}$ in \eqref{eq.def.kappa.scal}, one sees \eqref{ineq.ensreu.tau} holds if
  \bequationNN
    \rho_{\Scal_{k+1}}=\overline{\rho_{\Scal_{k+1}}}(\underline{\delta})\ge \tfrac{3\max\{1,\kappa_{\Scal_{k+1}}\}}{2\sigma_{\min}}.
  \eequationNN
  For any $\delta>0$, the definition $\overline{\rho_{\Scal_{k+1}}}(\delta) = \tfrac{\thalf \beta_{\Scal_{k+1}} - \eta_{\Scal_{k+1},1} \delta}{\eta_{\Scal_{k+1},2} \delta}$ yields
  \begin{align*}
    & \left\{\overline{\rho_{\Scal_{k+1}}}(\delta) \geq \tfrac{3 \max\{1,\kappa_{\Scal_{k+1}}\}}{2\sigma_{\min}} \right\} \iff \left\{ \tfrac{\thalf\beta_{\Scal_{k+1}} - \eta_{\Scal_{k+1},1} \delta}{\eta_{\Scal_{k+1},2} \delta} \geq \tfrac{3\max\{1,\kappa_{\Scal_{k+1}}\}}{2\sigma_{\min}}\right\} \\
    &\iff \left\{\beta_{\Scal_{k+1}}\sigma_{\min} - 2 \eta_{\Scal_{k+1},1} \delta \sigma_{\min} \geq 3 \eta_{\Scal_{k+1},2} \max\{1,\kappa_{\Scal_{k+1}}\}\delta\right\} \\
    &\iff \left\{ \delta \leq \tfrac{\beta_{\Scal_{k+1}} \sigma_{\min}}{2 \eta_{\Scal_{k+1},1} \sigma_{\min} + 3 \eta_{\Scal_{k+1},2} \max\{1,\kappa_{\Scal_{k+1}}\}}\right\},
  \end{align*}
  where the last inequality is ensured to hold with $\underline{\delta}$ defined in \eqref{eq.def.delta_s}. Now combining inequalities \eqref{ineq.G1.G0} and \eqref{ineq.strong.convex.2} to write
  \bequationNN
    \|G_{i+1}\| \leq \tilde{\epsilon} := \(1 - \tfrac{1}{6} \beta_{\Scal_{k+1}} t_{\Scal_{k+1}}\)^{i+1} \|G_0\|\ \ \text{and}\ \ \nabla^2 F_{\Scal_{k+1}}(x_{\Scal_k}^{i+1}) \succeq \tfrac{1}{3} \beta_{\Scal_{k+1}} I,
  \eequationNN
  it follows from \cite[Corollary 2.7]{GoyeEfteBoum2024} that the iterate $x_{\Scal_k}^{i+1}$ is an $(\tilde{\epsilon},2\tilde{\epsilon},-\tfrac{\beta_{\Scal_{k+1}}}{3}+(1+\tfrac{2}{\sigma_{\min}})\kappa_{\Scal_{k+1}}\tilde{\epsilon})$ stationary point of \eqref{prob.opt.S}, in the sense that 
  \begin{align}
    &\ \|c_{\Scal_{k+1}}(x_{\Scal_k}^{i+1})\| \leq  \tilde{\epsilon},\ \ \|\nabla_x L_{\Scal}(x_{\Scal_k}^{i+1},y_{\Scal_{k+1}}(x_{\Scal_k}^{i+1}))\| \leq 2\tilde{\epsilon}, \nonumber \\
    \text{and} &\ d^T \nabla_{xx}^2 L_{\Scal_{k+1}}(x_{\Scal_k}^{i+1},y_{\Scal_{k+1}}(x_{\Scal_k}^{i+1}))d \geq \( \tfrac{1}{3} \beta_{\Scal_{k+1}} - \(1+\tfrac{2}{\sigma_{\min}}\) \kappa_{\Scal_{k+1}} \tilde{\epsilon}\) \|d\|^2 \nonumber \\
    \text{for all} &\ d \in \Null(\nabla c_{\Scal_{k+1}}(x_{\Scal_{k+1}}^{i+1})^T). \label{ineq.rela1}
  \end{align}
  Note that $\(1-\tfrac{1}{6}\beta_{\Scal_{k+1}}t_{\Scal_{k+1}}\) < 1$, so $\tilde{\epsilon}\le \|G_0\|$. Combined with the fact that $\|G_0\| \leq 3\overline{\omega}\epsilon_{\Scal_k}\le R$ by \eqref{ineq.S.lemma}, one finds that $\|c_{\Scal_{k+1}}(x_{\Scal_k}^{i+1})\|\le \tilde{\epsilon}\le R$. Thus, the desired conclusion that $x_{\Scal_{k}}^{i+1}\in\Ccal_{\Scal_{k+1},R}$ has been proved.
  
  Through the induction it has been shown that \eqref{ineq.rela1} holds for all $i \geq 0$.  Let us now show that that right-hand side of the second line of \eqref{ineq.rela1} is nonnegative for all $i \geq 0$.  Consider arbitrary such $i$. From $\beta_{\Scal_{k+1}}\ge\thalf\beta$ and $\tilde{\epsilon}\le\|G_0\|\le3\overline{\omega}\epsilon_k $,
  \bequationNN
    \tfrac{1}{3} \beta_{\Scal_{k+1}} - \(1 + \tfrac{2}{\sigma_{\min}}\) \kappa_{\Scal_{k+1}} \tilde{\epsilon} \geq \tfrac{1}{6} \beta - \(1+\tfrac{2}{\sigma_{\min}}\)\kappa_{\Scal_{k+1}}3\overline{\omega}\epsilon_k.
  \eequationNN
  Next, from \eqref{eq.def.delta_s} and \eqref{ineq.ensreu.tau}, it follows that
  \begin{align*}
    \overline{\omega}
    &= \max_{\Scal \subseteq [N] \st \eqref{eq.theorem2.S}} 1+ \kappa_\Scal + 2 \overline{\rho_\Scal}(\underline\delta) \kappa_\Scal\\
    &\geq \kappa_{\Scal_{k+1}} + 2 \overline{\rho_{\Scal_{k+1}}}(\underline\delta) \kappa_{\Scal_{k+1}} \geq \kappa_{\Scal_{k+1}} + \tfrac{3\max\{1,\kappa_{\Scal_{k+1}}\}}{\sigma_{\min}} \kappa_{\Scal_{k+1}} \geq \(1 + \tfrac{2}{\sigma_{\min}}\)\kappa_{\Scal_{k+1}}.
  \end{align*}
  Combined with the previous displayed inequality, one finds $\tfrac{1}{3}\beta_{\Scal_{k+1}} - (1+\tfrac{2}{\sigma_{\min}}) \kappa_{\Scal_{k+1}}\tilde{\epsilon} \geq \tfrac{1}{6} \beta - 3\overline{\omega}^2\epsilon_k$. Now combined with \eqref{ineq.S.lemma} and $\epsilon_k=\kappa_2\xi_{\Scal_k}$, one has
  \begin{align*}
    \tfrac{1}{3} \beta_{\Scal_{k+1}} - \(1+\tfrac{2}{\sigma_{\min}}\) \kappa_{\Scal_{k+1}} \tilde{\epsilon} &\geq \tfrac{1}{6} \beta - 3\overline{\omega}^2\epsilon_k\ge0.
  \end{align*}
  Thus, the right-hand side of the second line of \eqref{ineq.rela1} is nonnegative for all $i \geq 0$.
  
  Let us now bound the number of iterations until gradient descent terminates, which due to the result in the previous paragraph can be determined by the number of iterations until the first-order condition is satisfied.  Recall the fact that $\nabla_y L_{\Scal_{k+1}}(x_{\Scal_k}^{i+1},y_{\Scal_{k+1}}(x_{\Scal_k}^{i+1}))=c_{\Scal_{k+1}}(x_{\Scal_k}^{i+1})$.  By \eqref{ineq.rela1}, one finds for all $i \geq 0$ that
  \begin{align}
    &\ \|\nabla L_{\Scal_{k+1}}(x_{\Scal_k}^{i+1},y_{\Scal_{k+1}}(x_{\Scal_k}^{i+1}))\| \nonumber \\
    =&\ \sqrt{\|\nabla_x L_{\Scal_{k+1}}(x_{\Scal_k}^{i+1},y_{\Scal_{k+1}}(x_{\Scal_k}^{i+1}))\|^2+\|\nabla_y L_{\Scal_{k+1}}(x_{\Scal_k}^{i+1},y_{\Scal_{k+1}}(x_{\Scal_k}^{i+1}))\|^2} \nonumber \\
    \leq&\ \sqrt{5}\tilde{\epsilon}= \sqrt{5}\(1-\tfrac{1}{6}\beta_{\Scal_{k+1}}t_{\Scal_{k+1}}\)^{i+1}\|G_0\|. \label{ineq.rela2}
  \end{align}
  By the step-size rule in \eqref{eq.stepsize}, the right-hand side of this expression is monotonically decreasing as $i$ increases.  Letting $T$ denote the iteration index at which the first-order condition is satisfied, one finds from \eqref{ineq.rela2} that
  \bequationNN
    \sqrt{5}\(1-\tfrac{1}{6}\beta_{\Scal_{k+1}} t_{\Scal_{k+1}}\)^T \|G_0\| \leq \epsilon_{k+1} \Longrightarrow T \leq \left\lceil \log_{1-\tfrac{1}{6}\beta_{\Scal_{k+1}}t_{\Scal_{k+1}}}\tfrac{\epsilon_{k+1}}{\sqrt{5}\|G_0\|} \right\rceil.
  \eequationNN
  Moreover, from \eqref{ineq.i-0} and \eqref{eq.stepsize}, one finds for the log term that
  \begin{align*}
    \log_{1-\tfrac{1}{6}\beta_{\Scal_{k+1}}t_{\Scal_{k+1}}}\tfrac{\epsilon_{k+1}}{\sqrt{5} \|G_0\|}
    &= \log_{\tfrac{1}{1-\tfrac{1}{6}\beta_{\Scal_{k+1}}t_{\Scal_{k+1}}}}\tfrac{\sqrt{5}\|G_0\|}{\epsilon_{k+1}} \leq \log_{2}\tfrac{\sqrt{5}G_0}{\epsilon_{k+1}} \leq \log_{2}\tfrac{3\sqrt{5}\overline{\omega}\epsilon_k}{\epsilon_{k+1}},
  \end{align*}
  which completes the proof.
\eproof

\begin{theorem}\label{theo.complexity.OLD}
  Suppose that Assumptions \ref{ass.subproblems}, \ref{ass.boundness}, \ref{ass.bounded.distribute}, and \ref{ass.Hessian.Lip.F} hold.
  \benumerate
    \item[(a)] Suppose that with $p_1 = N$ the conditions of Lemma~\ref{lemma.fletures.lag.result} hold and that \cite[Algorithm 1]{GoyeEfteBoum2024} is employed to solve subproblem~\eqref{prob.opt.S} for $\Scal = \Scal_1$.  Then, there exists $(u_{[N],1},u_{[N],2}) \in \R{}_{>0} \times \R{}_{>0}$ such that, for any $(\epsilon,\zeta) \in (0,\tfrac{\sqrt{5}}{2}] \times (0,1]$, the number of constraint gradient evaluations that are required until Algorithm~\ref{alg.psm} terminates with an $(\epsilon,\zeta)$ stationary point of~\eqref{prob.opt.N} is
    \bequation\label{eq.total.grad.eval.dir}
      N \left\lceil \max \left\{u_{[N],1} \epsilon^{-2}, u_{[N],2} \zeta^{-3} \right\} \right\rceil.
    \eequation
    \item[(b)]  Suppose that with $p_1 < N$ the conditions of Lemma~\ref{lemma.fletures.lag.result} hold and that \cite[Algorithm 1]{GoyeEfteBoum2024} is employed to solve subproblem~\eqref{prob.opt.S} for $\Scal = \Scal_1$, whereas gradient descent is employed to minimize Fletcher's augmented Lagrangian function for all subsequent $k \in [K]$ under the conditions of Lemma~\ref{lemma.fl.strongly.convex}.  Then, there exists $(u_{\Scal_1,1},u_{\Scal_1,2}) \in \R{}_{>0} \times \R{}_{>0}$ such that, for any $(\epsilon,\zeta) \in (0,\tfrac{\sqrt{5}}{2}] \times (0,1]$ and with $(\tau_1,\kappa_1,\overline\omega)$ defined in Lemma~\ref{lemma.fl.strongly.convex}, the number of constraint gradient evaluations that are required until Algorithm~\ref{alg.psm} terminates with an $(\epsilon,\zeta)$ stationary point of~\eqref{prob.opt.N} is
    \begin{multline}\label{eq.total.grad.eval.our}
      |\Scal_1|\left\lceil\max\left\{u_{\Scal_1,1} \epsilon_1^{-2}, u_{\Scal_1,2} \zeta_1^{-3} \right\}\right\rceil \\ + \(\tfrac{N-p|\Scal_1|}{p-1}\)\left\lceil\tfrac{1}{2}\log_2 \(45\overline{\omega}^2p^2\(1+\tfrac{(p-1)N}{p}\)\)\right\rceil  + N\left\lceil \log_{2}\tfrac{3\sqrt{5}\overline{\omega}\tau_1\sqrt{ p^2-p }}{\kappa_1\epsilon}\right\rceil.
    \end{multline}
  \eenumerate
\end{theorem}
\bproof
  Part (a) of the theorem and the aspect of part (b) pertaining to the first subproblem follow from Lemma~\ref{lemma.fletures.lag.result}, \cite[Theorem 3.5]{GoyeEfteBoum2024}, and the fact that each iteration of \cite[Algorithm~1]{GoyeEfteBoum2024} requires $\Ocal(p_1)$ constraint gradient evaluations.

  Let us now proceed to analyze the aspect of part (b) pertaining to the subproblem solves after the initial one.  Consider two phases: (1) all subproblems in iterations prior to the last one, and (2) the final subproblem in iteration $K$.

  Consider iterations $k \in \{2,\dots,K-1\}$ of Algorithm~\ref{alg.psm}.  The aim of each iteration is to find an $(\epsilon_k,\zeta_k)$-stationary point using the final iterate from the previous subproblem solve as a starting point.  By Lemma~\ref{lemma.fl.strongly.convex}, for all $k \in \{2,\dots,K-1\}$ it follows that the number of iterations required by gradient descent is at most
  \bequationNN
    \left\lceil \log_{2}\tfrac{3\sqrt{5}\overline{\omega}\epsilon_{k-1}}{\epsilon_k}\right\rceil.
  \eequationNN
  Recalling that $|\Scal_k| = p|\Scal_{k-1}| < N$ for all such $k$, one finds that
  \begin{align*}
    \tfrac{\epsilon_{k-1}}{\epsilon_k}
    &=\tfrac{\tau_1\xi_{\Scal_{k-1}}}{\tau_1\xi_{\Scal_k}}=\sqrt{\tfrac{(N-|\Scal_{k-1}|)|\Scal_k|^2}{|\Scal_{k-1}|^2(N-|\Scal_k|)}}\\
    &=\sqrt{\tfrac{p^2(N-\tfrac{1}{p}|\Scal_k|)}{(N-|\Scal_k|)}}=p \sqrt{1+\tfrac{(p-1)|\Scal_k|}{p(N-|\Scal_k|)}} \leq p \sqrt{1+\tfrac{(p-1)N}{p}}.
  \end{align*}
  Here, the last inequality follows from the fact that $x/(N-x)$ increases as $x$ increases over $(0,N)$, and since $|\Scal_k| \in \N{}$ and $|\Scal_k|<N$ imply $|\Scal_k|\leq N-1$.  Since each iteration of gradient descent requires $|\Scal_k|$ constraint gradient evaluations, it follows that the total number of constraint gradient evaluations in iteration $k$ is at most
  \bequationNN
    |\Scal_k| \left\lceil \log_2 \(3\sqrt{5}\overline{\omega}p\sqrt{1 + \tfrac{(p-1)N}{p}} \) \right\rceil = |\Scal_k|\left\lceil\tfrac{1}{2}\log_2 \(45\overline{\omega}^2p^2\(1+\tfrac{(p-1)N}{p}\)\)\right\rceil.
  \eequationNN
  Thus, the total number of constraint gradient evaluations required to solve all of the subproblems in iterations $k \in \{2,\dots,K-1\}$ is at most
  \begin{align*}
    &\ \sum_{k=2}^{K-1}|\Scal_k|\left\lceil\tfrac{1}{2}\log_2 \(45\overline{\omega}^2p^2\(1+\tfrac{(p-1)N}{p}\)\)\right\rceil \\
    &= \left\lceil\tfrac{1}{2}\log_2 \(45\overline{\omega}^2p^2\(1+\tfrac{(p-1)N}{p}\)\)\right\rceil \sum_{i=1}^{\lceil\log_pN/|\Scal_1|\rceil-2}|\Scal_1|p^i\\
    &= \left\lceil\tfrac{1}{2}\log_2 \(45\overline{\omega}^2p^2\(1+\tfrac{(p-1)N}{p}\)\)\right\rceil \(p |\Scal_1| \tfrac{p^{\lceil\log_pN/|\Scal_1|\rceil-2}-1}{p-1}\)\\
    &\le \left\lceil\tfrac{1}{2}\log_2 \(45\overline{\omega}^2p^2\(1+\tfrac{(p-1)N}{p}\)\)\right\rceil \(p|\Scal_1| \tfrac{p^{\log_pN/|\Scal_1|-1}-1}{p-1}\)\\
    &=\left\lceil\tfrac{1}{2}\log_2 \(45\overline{\omega}^2p^2\(1+\tfrac{(p-1)N}{p}\)\)\right\rceil \(\tfrac{N-p|\Scal_1|}{p-1}\).
  \end{align*}
  Now consider iteration $k = K$, the aim of which is to find an $(\epsilon,\zeta)$-stationary point by gradient descent using the final iterate from the subproblem solve in iteration $K-1$ as the starting point.  From the definition that $\epsilon_{K-1}
      =\tau_1\xi_{\Scal_{K-1}}$ in \eqref{eq.tolorence}, the function $\xi_\Scal$ is decreasing when $|\Scal|$ increasing and that $|\Scal_{K-1}|\ge\tfrac{N}{p}$, we have
  \begin{equation}
      \label{ineq.last.tolorence}
      \epsilon_{K-1}
      =\tau_1\sqrt{\tfrac{(N-|\Scal_{K-1}|)N}{|\Scal_{K-1}|^2}}\le\tau_1\sqrt{\tfrac{(N-N/p)N}{(N/p)^2}}=\tau_1\sqrt{ p^2-p }.
  \end{equation}
  Combining Lemma~\ref{lemma.fl.strongly.convex}, the inequality \eqref{ineq.last.tolorence}, and the fact that each iteration involves $N$ individual constraint gradient evaluations, it follows that the number of constraint gradient evaluations that are required is at most
  \bequationNN
    N\left\lceil \log_{2}\tfrac{3\sqrt{5}\overline{\omega}\epsilon_{K-1}}{\epsilon}\right\rceil\le N\left\lceil \log_{2}\tfrac{3\sqrt{5}\overline{\omega}\tau_1\sqrt{ p^2-p }}{\kappa_1\epsilon}\right\rceil.
  \eequationNN
  Combining the bounds that have been proved, the proof is complete.
\eproof